\documentclass{amsart}

\usepackage[utf8]{inputenc}
\usepackage{amsfonts}
\usepackage{amsmath}
\usepackage{amsthm}
\usepackage{tikz-cd}
\usepackage{hyperref}
\usepackage{mathtools}
\usepackage{cleveref}
\usepackage{xcolor}

\usetikzlibrary{tikzmark,calc}
\newcommand{\tetrahedron}{\tikzmarknode[text opacity=0,path picture={%
\draw let \p1=($(path picture bounding box.north)-(path picture bounding
box.south)$),\p2=($(path picture bounding box.east)-(path picture bounding
box.west)$),\n1={0.5*min(\x2,\y1)},\n2={0.42*max(\x2,\y1)-\n1}
in ([yshift=-\n2]path picture bounding box.center) coordinate(aux0)
-- ++ (200:\n1) --
($(aux0)+(80:\n1)$)
-- (aux0) -- ++ (320:\n1)
-- ($(aux0)+(80:\n1)$)
($(aux0)+(320:\n1)$) --
($(aux0)+(200:\n1)$);
}]{aux}{W}}

\usepackage{lmodern}

\usepackage{soul}

\usepackage{quiver}
\usetikzlibrary{nfold}

\usepackage[all]{xy}

\usepackage[T1]{fontenc}
\usepackage{microtype}

\usepackage{geometry}
\title[Multiplicative Serre spectral sequences in HoTT]
{Multiplicativity of the Serre Spectral Sequence in Homotopy Type Theory}
\author{Arseniy Kryazhev}
\address{University of California San Diego, La Jolla, CA 92093, USA}
\email{akryazhev@ucsd.edu}

\newtheorem{theorem}{Theorem}[section]
\theoremstyle{definition}
\newtheorem{thm}[theorem]{Theorem}
\newtheorem{prp}[theorem]{Proposition}
\newtheorem{de}[theorem]{Definition}
\newtheorem{re}[theorem]{Remark}
\newtheorem{exa}[theorem]{Example}
\newtheorem{cor}[theorem]{Corollary}
\newtheorem{lem}[theorem]{Lemma}

\newcommand{\id}{\operatorname{id}}

\newcommand{\ZZ}{\mathbb{Z}}
\newcommand{\NN}{\mathbb{N}}
\newcommand{\CC}{\mathbb{C}}
\newcommand{\QQ}{\mathbb{Q}}
\newcommand{\ap}{\operatorname{ap}}
\newcommand{\refl}{\operatorname{refl}}

\newcommand{\inl}{\operatorname{inl}}
\newcommand{\inr}{\operatorname{inr}}
\newcommand{\glue}{\operatorname{glue}}
\newcommand{\gluel}{\operatorname{gluel}}
\newcommand{\gluer}{\operatorname{gluer}}
\newcommand{\auxl}{\operatorname{auxl}}
\newcommand{\auxr}{\operatorname{auxr}}

\newcommand{\loo}{\operatorname{loop}}
\newcommand{\merid}{\operatorname{merid}}
\newcommand{\mmerid}{\overline{{\operatorname{merid}}}}
\newcommand{\swap}{\bar{\operatorname{swap}}}
\newcommand{\fib}{\operatorname{fib}}
\newcommand{\im}{\operatorname{im}}
\newcommand{\tot}[2]{\sum\limits_{#1}^{}{#2}}
\newcommand{\depfun}[2]{\prod\limits_{#1}^{}{#2}}
\newcommand{\bs}[1]{\mathrm{pt}_{#1}}
\newcommand{\funext}{\operatorname{funext}}

\newcommand{\ptd}[1]{\operatorname{pted}_{#1}}
\newcommand{\pto}{\to_*}
\newcommand{\bsd}[1]{\operatorname{based}_{#1}}
\newcommand{\pin}{\operatorname{pinch}}

\newcommand{\adj}[1]{\operatorname{adj}^{#1}}
\newcommand{\oout}{\operatorname{\Omega-out}}

\newcommand{\siin}{\operatorname{\Sigma-in}}
\newcommand{\sout}{\operatorname{\Sigma-out}}
\newcommand{\sgtn}{\mathbf{1}}
\newcommand{\s}{\mathbf{2}}

\newcommand{\ctos}{\operatorname{ctos}}
\newcommand{\jtos}{\operatorname{jtos}}
\newcommand{\whisk}{\raisebox{0.5ex}{$\displaystyle\centerdot$}}
\newcommand{\join}{\star}
\newcommand{\lps}[2]{(#1, #2 \rangle}
\newcommand{\lpsl}[2]{[#1, #2 \rangle}
\newcommand{\rps}[2]{\langle #1, #2 )}
\newcommand{\rpsl}[2]{\langle #1, #2 ]}
\newcommand{\ooutk}[1]{\operatorname{\Omega-out}^{#1}}
\newcommand{\onum}[1]{\stackrel{#1}{\Omega}}
\newcommand{\snum}[1]{\stackrel{#1}{\Sigma}}

\newcommand{\cone}{\vartriangleright}
\newcommand{\mir}{\operatorname{mirror}}

\tikzset{%
    symbol/.style={%
        draw=none,
        every to/.append style={%
            edge node={node [sloped, allow upside down, auto=false]{$#1$}}}
    }
}

\begin{document}

\begin{abstract}
Working in homotopy type theory, we define relative homotopy groups and give three equivalent definitions of the product on these groups. Using this product, we establish the multiplicativity of the cohomological Serre spectral sequence. We show a connection between the homotopy product and the cup-product in cohomology, which leads to the classically familiar description of the starting page of the spectral sequence as a tensor product of graded rings.
As applications, we compute cohomology rings of several classifying spaces and
loop spaces of spheres.
\end{abstract}

\maketitle

\tableofcontents

\section{Introduction}

Algebraic topology is primarily concerned with algebraic characteristics of topological spaces that are invariant under continuous deformation, i.e. homotopy. Classically, these topological spaces are seen as sets with certain structure, or, more abstractly, as objects of a certain category. Homotopy type theory (HoTT) recasts this study in a type-theoretic light: now topological spaces manifest as types, and
homotopy is encoded by equalities between points (terms) of types. This perspective, unifying topology and type theory, is promising for several reasons. On the one hand, the notion of homotopy is now built in, which often leads to significant simplifications and deeper insights. On the other hand, this approach is a generalization of the classical one, by which we mean that homotopy type theory is modeled by any infinity-topos (\cite{Shu}). Moreover, homotopy type theory lends itself to computer formalization as it is implemented in modern proof assistants such as Coq, Agda and Rzk.

Among the most important algebraic characteristics studied by algebraic topology are homology $H_\bullet(-)$ and cohomology $H^\bullet(-)$, which are graded groups assigned to a space. Classically, one of the main tools for computing these invariants is the Serre spectral sequence.
It arises whenever one has a fibration (which can be thought of as a ``twisted direct product'') and encapsulates the relationship between the (co)homology of the three spaces involved. More concretely, given a fibration $E \xrightarrow{F} B$, there is a trigraded sequence of groups $E^{p,q}_r$ and homomorphisms $d^{p,q}_r: E^{p,q}_r \to E^{p+r,q-r+1}_r$ such that 1) $E^{p,q}_2 = H^p(B; H^q(F))$; 2) $E^{p,q}_{r+1} = \ker{d_r^{p,q}}/\im d_r^{p-r,q+r-1}$; 3) $H^{p+q}(E)$ is obtained by iterated extensions from the groups $E^{p,q}_r$ for large $r$. This is usually expressed by writing
$$ H^p(B; H^q(F)) \implies H^{p+q}(E) .$$

In HoTT, the cohomological spectral sequence has recently been constructed in \cite{Flo}. The merit of this result comes from the fact that certain things are challenging to construct in HoTT, which stems from two of its aspects. On the one hand, the logic involved has to be \textit{constructive}. This makes things like homological algebra much less straightforward. On the other hand, since ``equal'' means ``homotopic'', we cannot distinguish between homotopy equivalent spaces; thus, only homotopy-invariant notions can be used.

An important feature of cohomology is the existence of a multiplicative (``cup-product'') structure, which makes cohomology a graded ring rather than just a group. A natural question is then whether the Serre spectral sequence can keep track of this additional structure. Classically, the answer is positive: each $E^{p,q}_r$ can be converted into a ring in a way compatible with the differentials, the input and the output of the spectral sequence.

The present paper imports this result into HoTT -- see the statement of Theorem \ref{th:main-result}. Classically, one deals with cohomology on the level of cochains and then defines all the products essentially via a multiplication of ring-valued functions (which cochains are). However, in HoTT it is challenging to implement this approach as it is unclear what the domain of a cochain should be. It cannot be the set of singular simplices because the latter are not homotopy-invariant (any simplex is contractible). It could be cellular chains, but this restricts the theory to cellular complexes, while also raising new technical issues.
Thus, neither the spectral sequence nor the cohomological product is defined in HoTT in the ``classic textbook'' way via cochains. Instead, one uses the definition of cohomology $H^n(-,R)$ as (path components of) mapping spaces into classifying spaces $K_n R$; that is, one sets $H^n(X,R) : \equiv ||X \to K_n R||_0$. Then one can leverage the ring-space structure of $K_n R$ to give a Yoneda-style definition of the product\footnote{See \cite{LjuMor} for a comprehensive exposition.}.

The difficulty then is that the constructions of the spectral sequence and the cup-product become much more distant than classically, and it is no longer so apparent how to connect them. Our approach is guided by \cite{Dug1}, where the product arises from the pairing on relative homotopy groups
$$\pi_n(B,A) \otimes \pi_m(D,C) \to \pi_{n+m}(B \wedge D, (A \wedge D) \coprod^{A \wedge C}(B \wedge C) ).$$
We notice that the starting term of the spectral sequence of \cite{Flo} can be given in the language of relative homotopy groups, whereafter the product above is carried to the spectral sequence. The main challenge then is to overcome the technical difficulties in proving that everything works as expected. For example, the pairing above involves a three-dimensional diagram with smash-products at the vertices, so proving that it is distributive over addition will necessitate 5-fillers for 4-dimensional diagrams. To deal with these problems, we provide two additional perspectives on the pairing: one with the cross-product instead of the smash-product, and an inductive one. We prove that they are all equivalent, and then we use the appropriate one for each property of the product we try to prove.

Some more work is required to prove that the multiplication we introduce specializes to the cup-product pairing. We prove it via induction on the dimensions.

Another major challenge throughout our work is getting the signs right. Roughly, switching between different perspectives forces us to permute certain functors, which introduces numerous minus signs. For example, the product on Eilenberg--Mac Lane spaces $K_n G \times K_m G \to K_{n+m}G$ can be alternatively seen as $\Omega K_{n+1} G \times \Omega K_{m+1} G \to \Omega^2 K_{n+m+2}G$ and we'd like to pair loop functors $\Omega$ with each other first, for which we'd need to ``permute'' $\Omega$ and $K_{n+1}$ (see Lemma \ref{le:pairing-gives-inner-prod}).

Once the multiplicative spectral sequence is obtained, one can reap numerous fruits. We illustrate this by performing several standard computations, such as those of cohomology rings of some classifying spaces and loop spaces of spheres. One of the promising avenues that this opens is towards homotopy group calculations via Serre killing, but the lack of a universal coefficient theorem in any form remains the major obstacle.

The results are being formalized using the Coq HoTT library.

The structure of the paper is as follows.

In Section \ref{se:prelim}, we recall the basics about pointed types and give some necessary definitions.

In Section \ref{se:pairings}, we define the pairing operation: first on homotopy groups and then on spectral sequences. To prove the desired properties of the pairing -- naturality, distributivity and the Leibniz relation -- we give three alternative definitions of the pairing: an inductive one, one using pushout products, and one using pushout smashes.

Then, in Section \ref{se:SSS}, we show how this construction specializes to the classical statement of the multiplicativity of the Serre spectral sequence.

Finally, in Section \ref{se:appls}, we show how this additional structure can be used to compute cohomology rings of certain types.

Most technical lemmas are aggregated in Section \ref{se:lemmas} to make it easier to trace the general logic in the main sections.

\subsection*{AI usage statement}

The mathematics of this text is due to the author. AI was used only for the final edits: formatting, typo correction, reference list, and to fill in several short paragraphs with (straightforward) formulas.

\subsection*{Acknowledgements}

I am grateful to Eva Belmont for pointing me to \cite{Dug1}, to Dan Christensen for feedback on the early version of the text, and to Michael Shulman for suggestions and feedback.

\section{Preliminaries}\label{se:prelim}

We use ``book HoTT'' (\cite{UF}): Martin-L\"of type theory with univalence and higher inductive types.

All of our notation and conventions are from \cite{UF} with a couple of exceptions. We use $\sgtn$ for the singleton type and $\s$ for the type with two elements $0_{\s}$ and $1_\s$. For path concatenation we use $\cdot$, which we can sometimes suppress.
The universe is denoted by $U$.
\\
We write $cancel$ for the (essentially unique) paths witnessing obvious path algebra equalities such as $p \cdot p^{-1} = \refl$.

\subsection{The universe of pointed types}
\begin{de}
    Let $U_* :\equiv \tot{X:U}{X}$. An element $(X, \bs{X}) : U_*$ is called a \textit{pointed type} and $\bs{X} : X$ is called the \textit{basepoint}. A map of pointed types $\phi: (X, \bs{X}) \to (Y, \bs{Y})$  (a ``\textit{pointed map}'') is a map $f: X \to Y$ together with a witness of pointedness $\ptd{\phi}: f(\bs{X}) = \bs{Y}$. The type of these will be denoted by $(X, \bs{X}) \pto (Y, \bs{Y})$, or $X \pto Y$ when the basepoints are clear from the context.
    We will also abuse notation by using the same symbol for $\phi$ and the underlying map $\pi_1 \phi$.
\end{de}

\begin{re}
    Recall that in HoTT a \textit{wild category} is defined via a propositions-as-types translation of the classical definition of a category (\cite{SphereSyms}, Definition B.1). Pointed types and pointed maps form a wild category. In particular, one can compose pointed maps by $(f_1, \ptd{1}) \circ (f_2, \ptd{2}) : \equiv (f_1 \circ f_2, \ap_{f_1}(\ptd{2}) \cdot \ptd{1})$.
\end{re}

Given a homotopy (i.e. an equality) $\alpha$ between functions $f,g: A \to B$, we will write $\alpha(a)$ for $\operatorname{happly}(\alpha)(a): f(a) =_B g(a)$. By function extensionality, $\alpha$ is determined by $\lambda {a}.{\alpha(a)}$.
\begin{de}
    A \textit{based homotopy} between two pointed functions $f,g: X \pto Y$ is a homotopy between the underlying functions $h : f =_{X \to Y} g$ together with a witness that it is based $\bsd{h}: h(\bs{X}) = \ptd{f} \cdot \ptd{g}^{-1} $.
\end{de}

\begin{re}
    Unfolding the definitions, one sees that showing an equality of two pointed maps $f,g : (X, \bs{X}) \to (Y, \bs{Y})$ amounts to exhibiting a based homotopy between them. We will use the notation $=_*$ for such equalities.
\end{re}

\begin{de}
    When we refer to a \textit{diagram in $U_*$} (or $U$), we mean a diagram in the wild category $U_*$ (or $U$). Thus, a 2-dimensional diagram in $U_*$ is a diagram with objects pointed, maps pointed and 2-cells (homotopies) based. A 1-dimensional diagram in $U_*$ determining parallel arrows $f,g: a\to b$ \textit{commutes} if it has a based 2-dimensional filler $h: f =_{a\pto b} g$. For a 2-dimensional diagram in $U_*$, the condition of commutativity is quite involved if spelled out explicitly; fortunately, we only need simple special cases such as Remark \ref{re:equality-of-squares-simple-version}.
\end{de}

\begin{de}\label{de:pointed-equiv}
    A \textit{pointed equivalence} between pointed types is a pointed map whose underlying function is an equivalence of types. By \cite{Flo}, 2.2.7, this is equivalent to being an invertible arrow in the wild category $U_*$.
\end{de}

\begin{re}
    Generalizing the previous remark, from univalence it follows that an equality (in $\tot{X,Y}{X \pto Y}$) between $f : (X, \bs{X}) \to (Y, \bs{Y})$ and $f' : (X', \bs{X'}) \to (Y', \bs{Y'})$ is equivalently two pointed equivalences $g_X: X \pto X', g_Y: Y \pto Y'$ such that the resulting square commutes in $U_*$ (so that $g_Y \circ f =_{X \pto Y'} f' \circ g_X$).
\end{re}

\begin{de}
    By a \textit{pointed pair} we will understand any pointed map $i: A \pto B$. A map of pointed pairs is a commutative (in $U_*$) square. We will use the notation $(B, A)$ when $i$ is clear from the context.
\end{de}

\begin{re}
    A pointed type $(X, \bs{X})$ is equivalently the pointed pair $\sgtn \xrightarrow{\bs{X}} X$. A map of pointed types $(X, \bs{X}) \pto (Y, \bs{Y})$ is equivalently a map of pointed pairs $(X,\sgtn) \pto (Y,\sgtn)$.
\end{re}

\begin{de}
    Suppose given arrows $f,g: A \to B$ and $h: B \to C$ in a wild category $\mathbb{C}$, and an equality $\alpha: f =_{\operatorname{Hom}_{\mathbb{C}}(A,B)} g$. Define the \textit{whiskered homotopy} $h \whisk \alpha: h \circ f =_{\operatorname{Hom}_{\mathbb{C}}(A,C)} h \circ g$ as
    $$h \whisk \alpha : \equiv \ap_{h \circ -}(\alpha).$$
    We will be interested in the cases $\mathbb{C} \equiv U$ and $\mathbb{C} \equiv U_*$. Explicitly, in the former case, the resulting homotopy is (essentially uniquely) determined by the relation $(h \whisk \alpha)(a) = \ap_h (\alpha(a))$. In the latter case, the underlying unbased homotopy is the unbased whisker, and it is based by
    $(h \whisk \alpha)(\bs{A}) = \ap_h(\alpha(\bs{A})) =
    \ap_h(\ptd{f} \cdot \ptd{g}^{-1}) =
    \ap_h(\ptd{f}) \cdot \ap_h(\ptd{g})^{-1} =
    \ap_h(\ptd{f}) \cdot \ptd{h} \cdot \ptd{h}^{-1} \cdot \ap_h(\ptd{g})^{-1} =
    \ptd{h \circ f} \ptd{h \circ g}^{-1}.
    $
\end{de}

Whiskering $\whisk$ has obvious algebraic properties such as $h \whisk (\alpha \cdot \beta) = (h \whisk \alpha) \cdot (h \whisk \beta)$ and $(h' \circ h) \whisk \alpha = h' \whisk (h \whisk \alpha)$.

Another manifestation of the same motif is the whisker of a 1-dimensional equality with a 2-dimensional equality. Suppose given a path $r: x=_X y$ in some type $X$, and a 2-cell $\eta:p=_{y=z}q$ between two other paths $p,q: y=_Xz$:
\[\begin{tikzcd}
	x & y & z
	\arrow["r", squiggly, from=1-1, to=1-2]
	\arrow[""{name=0, anchor=center, inner sep=0}, "p", curve={height=-6pt}, squiggly, from=1-2, to=1-3]
	\arrow[""{name=1, anchor=center, inner sep=0}, "q"', curve={height=6pt}, squiggly, from=1-2, to=1-3]
	\arrow["\eta"', shorten <=2pt, shorten >=2pt, Rightarrow, from=0, to=1]
\end{tikzcd}\]
Then there is a path in the type $r \cdot p=_{x=z} r \cdot q$, denoted by $r \whisk \eta $. (It equals $\eta$ for $r \equiv \refl_y$ and extends to arbitrary $r$ by path induction.)
Similarly, given $r': z =_Xt$, we can whisker $\eta$ from the right to get $\eta \whisk r': p \cdot r'=_{y=t}  q \cdot r'$.

\subsection{Constructions with types}

Given a span $C \xleftarrow{f} A \xrightarrow{g} B$ in $U_*$, the pushout $C \coprod^A B$ will be pointed by $\inl(\bs{C}).$

\begin{lem}\label{pushout}
    The pointed pushout $(C \coprod^A B, \inl(\bs{C}))$ is the pushout of the span $C \xleftarrow{f} A \xrightarrow{g} B$ in $U_*$.
\end{lem}
See Lemma \ref{pushout-proof} for the detailed formulation and the proof.

\begin{de}
    For $A :U$, define its \textit{suspension} $\Sigma A: U_*$ as the pushout of $\sgtn \xleftarrow{} A \xrightarrow{} \sgtn$ and denote $N :\equiv \inl(*), S :\equiv \inr(*), \merid_a :\equiv \glue(a): N =_{\Sigma A} S$. If $A$ is pointed, denote $ \mmerid_a :\equiv \merid_a \cdot \merid_{\bs{A}}^{-1} : \Omega \Sigma A$.
\end{de}

\begin{de}
    For $A : U$, define recursively the \textit{iterated suspension} $\Sigma^{\bullet} A : \NN \to U$ as
    \begin{align*}
        \Sigma^0 A : \equiv A ;\\
        \Sigma^{n+1} A : \equiv \Sigma \Sigma^n A.
    \end{align*}
\end{de}
For $n \geq 1$ we can think of $\Sigma^n A$ as a pointed type with basepoint $N$.

\begin{de}
    For $A : U_*$, define recursively the \textit{iterated loop space} $\Omega^{\bullet} A : \NN \to U_*$ as
    \begin{align*}
        \Omega^0 A : \equiv A; \\
        \Omega^{n+1} A : \equiv \Omega^n \Omega A.
    \end{align*}
\end{de}

Recall that in $U_*$ there is an equivalence of types $(\Sigma A \pto B) \to (A \pto \Omega B)$ given by $f \mapsto \lambda a.\ptd{f}^{-1}\cdot \ap_f\mmerid_{a}\cdot \ptd{f}$.
This adjunction $\Sigma \dashv \Omega$ leads recursively to the adjunction (denoted $\adj{n}$) $\Sigma^n \dashv \Omega^n$.
We will use it to identify the sets $\pi_n(X) : \equiv ||\Omega^n X||_0$ and $||S^n \pto X||_0$.

\begin{de}
    For $A, B: U_*$, the \textit{wedge} $A \vee B : U_*$ of $A$ and $B$ is the higher inductive type generated by $\inl: A \to A \vee B, \inr: B \to A \vee B, \bs{\vee}: A \vee B, \gluel: \inl(\bs{A}) = \bs{\vee}, \gluer: \inr(\bs{B}) = \bs{\vee}$. It is pointed by $\bs{\vee}$, and $\inl$ and $\inr$ can be considered as pointed maps, respectively by $\gluel$ and $\gluer$.
\end{de}

\begin{de}
    Like in \cite{Flo}, 4.1.6, we define the \textit{smash-product} $A \wedge B$ of pointed types $A$ and $B$ as the higher inductive type generated by
    \begin{align*}
    & \auxl : A \wedge B; \\
    & \auxr : A \wedge B; \\
    & \langle \cdot, \cdot \rangle: A \times B \to A \wedge B; \\
    & \gluel: \depfun{a: A}\langle a, \bs{B}\rangle = \auxl; \\
    & \gluer: \depfun{b: B}\langle \bs{A},b \rangle = \auxr.
    \end{align*}
    In other words, $A \wedge B$ is the pushout of
    \begin{tikzcd}
    	{A + B} & \s \\
    	{A \times B.}
    	\arrow[from=1-1, to=1-2]
    	\arrow[from=1-1, to=2-1]
    \end{tikzcd}
\end{de}
It will be pointed by $\auxr$.
We will also denote the projection map $A \times B \to A \wedge B$ by $\ctos_{A,B}$, or by $\ctos$ when $A, B$ are clear from the context. (``$\ctos$'' meaning ``cross to smash''.)

The smash-product is functorial (\cite{Flo}, 4.3.8).

Let $i: A \vee B \to A \times B$ be the map defined by $\inl(a) \mapsto (a, \bs{B})$, $\inr(b) \mapsto (\bs{A}, b)$, $\bs{\vee} \mapsto (\bs{A}, \bs{B})$ and $\gluel, \gluer \mapsto \refl_{(\bs{A}, \bs{B})}$; it is pointed by reflexivity. Note that all paths and higher paths coming from the map $$A \vee B \xrightarrow{i} A \times B\to A \wedge B$$ are contractible. We will denote by $byContr$ the witness of any such contractibility.

\begin{lem}\label{le:c-intro}
    Given maps $f: A \pto B$, $g: C \pto D,$ the diagram
\[\begin{tikzcd}
	{A \wedge C} & {A \wedge D} \\
	{B \wedge C} & {B \wedge D}
	\arrow["{f \wedge C}"', from=1-1, to=2-1]
	\arrow["{A \wedge g}", from=1-1, to=1-2]
	\arrow["{f \wedge D}", from=1-2, to=2-2]
	\arrow["{B \wedge g}", from=2-1, to=2-2]
\end{tikzcd}\]
commutes.
\end{lem}

\begin{de}\label{de:c}
    The proof below will be denoted by $c_{f,g} : (f \wedge D) \circ (A \wedge g) =_* (B \wedge g) \circ (f \wedge C)$.
\end{de}
\begin{proof}
    By path induction, we can assume that $f,g$ are pointed by $\refl$.

    On all 0-dimensional generators of $A \wedge C$ set $c_{f,g}$ to be reflexivity. Suppose we move from $\langle \bs{A}, c \rangle$ to $\auxr$ along $\gluer_c$. Then we need a filler for
\[\begin{tikzcd}
	{\langle f(\bs{A}), g(c) \rangle} & {\langle f(\bs{A}), g(c) \rangle} \\
	\auxr & {\auxr.}
	\arrow["\refl", squiggly, from=1-1, to=1-2]
	\arrow["{\ap_{(f \wedge D) \circ (A \wedge g)} (\gluer_{c})}"', squiggly, from=1-1, to=2-1]
	\arrow["{\ap_{(B \wedge g) \circ (f \wedge C)}( \gluer_{c})}", squiggly, from=1-2, to=2-2]
	\arrow["\refl", squiggly, from=2-1, to=2-2]
\end{tikzcd}\]
The left path computes as $\ap_{(f \wedge D) \circ (A \wedge g)} (\gluer_{c}) = \ap_{f \wedge D} \ap_{A \wedge g} (\gluer_{c}) = \ap_{f \wedge D} \gluer_{gc} = \gluer_{gc}.$ The right path similarly computes to the same.

A similar filler is exhibited for $\gluel_a$.

It is easily seen that the constructed homotopy is based. Indeed, its value on $\auxr$ equals reflexivity, and so do the proofs of pointedness for all the functions involved.
\end{proof}

We will repeatedly use that the construction above is natural: the cube constructed from two squares as in the definition of $c$ commutes. See Lemma \ref{le:smashes-cube}.

\begin{de}
    Given two types $X, Y: U$, define their \textit{join} $X \join Y$ as the pushout of the span $X \xleftarrow{\pi_1} X \times Y \xrightarrow{\pi_2} Y$ (thus pointed by $\inl(\bs{X})$ when $X$ is pointed).
\end{de}

It is easy to show that the join is functorial.

\begin{re}
    Using the $3 \times 3$-lemma (\cite{Bru}, Lemma 1.8.3), one can show that $A \join B \simeq \Sigma (A \wedge B)$. We won't need this, using only the natural map $\jtos: A \join B \to \Sigma (A \wedge B)$ induced from the span map
    \[\begin{tikzcd}
    	{A \times B} & A \\
    	B & {A \wedge B} & \sgtn \\
    	& {\sgtn.}
    	\arrow[from=1-1, to=1-2]
    	\arrow[from=1-1, to=2-1]
    	\arrow[from=1-1, to=2-2]
    	\arrow[from=1-2, to=2-3]
    	\arrow[from=2-1, to=3-2]
    	\arrow[from=2-2, to=2-3]
    	\arrow[from=2-2, to=3-2]
    \end{tikzcd}\]
    (Note that the maps above are equivalently $A \times \sgtn \to A \wedge \sgtn$ and $\sgtn \times B \to \sgtn\wedge B$.)
\end{re}

\subsubsection{Paths in product types}\label{sec:paths-in-prod}

We introduce notation to abstract out some intuition behind technical definitions and proofs.

\begin{re}\label{re:balance}
We try to strike a balance between being precise and being understandable. Most notably, this manifests itself in how we treat higher coherences. For example, a square with two sides $p$ and two sides $\refl$ can be seen as many equalities, such as $p \cdot \refl = p \cdot \refl,$ $p^{} \cdot p^{-1} \cdot \refl = \refl$, $p \cdot p^{-1} = \refl \cdot \refl,$ etc. While we usually are explicit about the types of such cells, often we prefer to just speak about ``the square''. Indeed, many proof ideas can be easily expressed in such a language but become impenetrable when written as formulas; at the same time, knowing the geometric idea makes extracting a precise formula (when it is necessary) mechanical.

In many of our technical proofs, in the end the boundary of a cell gets contracted, in which case the 2-cell, as a path between two halves of its boundary, is only important up to a rotation. One has to remember, however, that these cells are oriented; the same geometric shape can represent both a (higher) path and its opposite. Thus, if the boundary of a cell is symmetric, focusing only on its shape can lead to a sign error. For example, one can keep in mind whether a cell is traversed clockwise or counterclockwise.
\end{re}

Recall that for any product type $X \times Y$ there are equivalences $(x=_X x') \times (y=_Y y') \to (x,y) =_{X \times Y} (x',y')$. When it is not ambiguous, we use $(p, y)$ to denote the image of $(p, \refl_y)$, and similarly for $(x, q)$. We have $(p, y) = \ap_{(-,y)}(p)$. The image of $(p, y)$ in $X \wedge Y$ will be denoted by $\lps{p}{y}$, and similarly for $\rps{x}{q}$. Thus, $\lps{p}{y} = \ap_{\langle -, y \rangle}(p)$.
\\
Next, given $X,Y:U_*, p : \Omega X, q: \Omega Y, x : X, y : Y$, define the versions of the same paths but
made into loops at the basepoint:
    $$[p, y\rangle : \Omega (X \wedge Y) :\equiv \gluer_y^{-1} \cdot \lps{p}y{} \cdot \gluer_y;$$
    $$\langle x, q] : \Omega (X \wedge Y) :\equiv \gluer_{\bs{Y}}^{-1} \cdot \gluel_{\bs{X}} \cdot \gluel_{x}^{-1} \cdot \rps{x}{q} \cdot (\gluer_{\bs{Y}}^{-1} \cdot \gluel_{\bs{X}} \cdot \gluel_{x}^{-1})^{-1}.$$

If $A, B$ are types and $p: a =_A a', q: b =_B b'$ are paths, there is an obvious square filler of type $\ap_{(-, b)}(p) \cdot \ap_{(a', -)}(q) =_{(a,b)=(a',b')} \ap_{(a, -)}(q)  \cdot \ap_{(-, b')}(p)$ in $A \times B$, defined by path induction. We will denote its image in $A \wedge B$ by $[p,q]$. Thus, $[p,q]$ is of type $\lps{p}{b} \cdot \rps{a'}{q} =_{\langle a,b\rangle=\langle  a',b'\rangle} \rps{a}{q} \cdot \lps{p}{b'}$.
If $p,q$ are paths in $A \vee B$ (i.e. $a\equiv a' \equiv \bs{A}, b\equiv b' \equiv \bs{B}$), then the boundary of $[p,q]$ is contractible, and we denote by $[p,q]^l$ the double loop $byContr \cdot [p,q] \cdot byContr^{-1}: \Omega \Omega (A \wedge B)$. For example, if $\loo$ generates $S^1$, then $[\loo, \loo]^l : \Omega \Omega (S^{1} \wedge S^{1})$ can be seen as representing the sphere as a square with identified boundary points.

\subsubsection{-in and -out functions}

Let $X, Y$ be types and $y:Y$. Let
$$\cone_y:\depfun{x,x':X}{\depfun{p:x=x'}{\ap_{\inl}(p) =_{\inl(x) =_{X \join Y} \inl(x')} \glue(x,y) \cdot \glue(x',y)^{-1}}}$$
be defined in the obvious way by path induction. (This dependently equals $\ap_{\glue(-,y)}(p^{-1}).$)
For $X : U_*, p: \Omega X$, denote $c(p,y) :\equiv cancel_* (\cone_y(\bs{X}, \bs{X}, p)) : \ap_{\inl}(p) =_{\Omega (X \join Y)} \refl_{\inl(\bs{X})}$.
Intuitively, $c(p,y)$ is a cone in $X \join Y$ connecting $p$ to $y$.

We now define several -in and -out maps that enable recursive definitions and proofs. The concrete details in the formulas are not as important as the fact that natural functions with such domains and codomains exist\footnote{Recall that an endofunctor $F$ on a monoidal category $(\CC, I, \otimes)$ is called \textit{strong} if there is a well-behaved natural transformation $F A  \otimes B \to F ( A \otimes B)$.}. For a detailed discussion of their properties, see Section \ref{se:-in-out}.
Note that we abuse notation by using the same symbol for different functions. Such polymorphism should not lead to confusion and frees formulas from clutter.

\begin{de}
    Define $\oout: \depfun{A, B: U_*}{(\Omega A \wedge B \to \Omega (A \wedge B)})$ by $\langle p, b \rangle \mapsto \lpsl{p}{b}$, $\auxr \mapsto \refl_{\bs{}}, \auxl \mapsto \refl_{\bs{}}$, $\gluel_p \mapsto byContr \cdot cancel, \gluer_b \mapsto byContr \cdot cancel$.
\end{de}
\begin{de}
    Define $\oout: \depfun{A: U_*,B: U}{(\Omega A \join B \to \Omega (A \join B)})$ by $\inl(p) \mapsto \ap_{\inl}(p); \inr(b) \mapsto \refl_{\bs{}}; \glue(p, b) \mapsto c(p,b)$.
\end{de}
\begin{de}
    Define $\siin: \depfun{A, B: U_*}{(\Sigma(A \wedge B) \pto \Sigma A \wedge B})$ by $N \mapsto \auxr; S \mapsto \auxr; \merid_{\langle a, b \rangle} \mapsto [\mmerid_a, b \rangle; \merid_{\auxr} \mapsto \refl_{\bs{}}; \merid_{\auxl} \mapsto \refl_{\bs{}}; \ap_{\merid}(\gluer_b) \mapsto (\gluer_b^{-1} \whisk cancel \whisk \gluer_b) \cdot cancel; \ap_{\merid}(\gluel_a) \mapsto byContr \cdot cancel$.
\end{de}
\begin{de}
    Define $\siin: \depfun{A, B: U}{(\Sigma (A \join B) \pto \Sigma A \join B})$ by $N \mapsto \inl(N); S \mapsto \inl(S); \merid(\inl(a)) \mapsto \ap_{\inl}(\merid(a)); \merid(\inr(b)) \mapsto \glue(N, b) \cdot \glue(S, b)^{-1};$ $\ap_{\merid}(\glue(a,b)) \mapsto \cone_b(N, S, \merid_a)$.
\end{de}

When regarded as pointed maps, the join $\oout$ is pointed by the evident path $\ap_{\inl}(\refl_{\bs{}}) = \refl_{\bs{}},$ and the other three maps are pointed by reflexivity.

\begin{de}
    We also define in a similar way functions $\siin_2$ (and $\oout_2$) that move the functor into (or out of) the second factor. For example, $\oout_2: \depfun{A: U, B:U_*}{(A \join \Omega B \to \Omega (A \join B))}$
    is defined on $A, B: U_*$ by  $\inr(p) \mapsto \glue(\bs{A}, \bs{B})\ap_{\inr}(p)\glue(\bs{A}, \bs{B})^{-1}; \inl(a) \mapsto \refl_{\bs{}}; \glue(a, p) \mapsto c_2(a,p)$, with the ``left-pointing'' cone $c_2(a,p): \refl_{\bs{}} = \glue(\bs{A}, \bs{B})\ap_{\inr}(p)\glue(\bs{A}, \bs{B})^{-1}$.
\end{de}

Note that applying $\Sigma$ is the same as smashing with $S^1$ (\cite{Bru}), so $\Sigma (A \wedge B) \simeq S^1 \wedge(A \wedge B) \simeq (S^1 \wedge A) \wedge B \simeq  \Sigma A \wedge B$, so unsurprisingly, $\siin$ is invertible (Lemma \ref{le:sout-w}). We may refer to its inverse as $\sout$.

Our definitions of the pairing will rely on the following construction.
\begin{de}
    Given a commutative square, its \textit{cogap map} is the induced map from the pushout to the bottom right object, as in the diagram below.
    \[\begin{tikzcd}
    	\bullet & \bullet \\
    	\bullet & \bullet \\
    	&& \bullet
    	\arrow[from=1-1, to=1-2]
    	\arrow[from=1-1, to=2-1]
    	\arrow[from=1-2, to=2-2]
    	\arrow[from=1-2, to=3-3]
    	\arrow[from=2-1, to=2-2]
    	\arrow[from=2-1, to=3-3]
    	\arrow["\lrcorner"{anchor=center, pos=0.125, rotate=180}, draw=none, from=2-2, to=1-1]
    	\arrow[dashed, from=2-2, to=3-3]
    \end{tikzcd}\]
\end{de}

\begin{de}
    Given a pair of maps $f: A \to B, g: C \to D$, define the \textit{pushout product} $P^\times(f,g)$ as the canonical map $(A \times D) \coprod^{(A \times C)} (B \times C) \to (B \times D)$. In other words, the pushout product of $f$ and $g$ is the following cogap map:
\[\begin{tikzcd}
	{A \times C} & {A \times D} \\
	{B \times C} & \bullet \\
	&& {B \times D.}
	\arrow["{A \times g}", from=1-1, to=1-2]
	\arrow["{f \times C}"', from=1-1, to=2-1]
	\arrow[from=1-2, to=2-2]
	\arrow["{f \times D}", from=1-2, to=3-3]
	\arrow[from=2-1, to=2-2]
	\arrow["{B \times g}"', from=2-1, to=3-3]
	\arrow["\lrcorner"{anchor=center, pos=0.125, rotate=180}, draw=none, from=2-2, to=1-1]
	\arrow[dashed, from=2-2, to=3-3]
\end{tikzcd}\]

\end{de}

\begin{de}\label{de:PO-smash}
    Given a pair of pointed maps $f: A \pto B, g: C \pto D$, define the \textit{pushout smash} $P^\wedge(f,g)$ as the canonical map $(A \wedge D) \coprod^{(A \wedge C)} (B \wedge C) \to (B \wedge D)$.
\end{de}
For pointed $f$ and $g$, the maps $P^\times(f,g)$ and $P^\wedge(f,g)$ are pointed in the obvious way.

The cogap maps are natural in the following sense. Assume given a filler for the outer parallelepiped in
\[\begin{tikzcd}
	&& \bullet & \bullet \\
	\bullet & \bullet & \bullet & \bullet \\
	\bullet & \bullet &&& \bullet \\
	&& {\bullet.}
	\arrow[from=1-3, to=1-4]
	\arrow[from=1-3, to=2-3]
	\arrow["\inl"', from=1-4, to=2-4]
	\arrow[from=1-4, to=3-5]
	\arrow[from=2-1, to=1-3]
	\arrow[from=2-1, to=2-2]
	\arrow[from=2-1, to=3-1]
	\arrow[from=2-2, to=1-4]
	\arrow[from=2-2, to=3-2]
	\arrow[from=2-2, to=4-3]
	\arrow["\inr", from=2-3, to=2-4]
	\arrow[from=2-3, to=3-5]
	\arrow["\lrcorner"{anchor=center, pos=0.125, rotate=180}, draw=none, from=2-4, to=1-3]
	\arrow[dashed, from=2-4, to=3-5]
	\arrow[from=3-1, to=2-3]
	\arrow[from=3-1, to=3-2]
	\arrow[from=3-1, to=4-3]
	\arrow["\lrcorner"{anchor=center, pos=0.125, rotate=180}, draw=none, from=3-2, to=2-1]
	\arrow[dashed, from=3-2, to=2-4]
	\arrow[dashed, from=3-2, to=4-3]
	\arrow[from=4-3, to=3-5]
\end{tikzcd}\]
Then we get the commutativity of the rectangle whose three sides are dashed maps. Its commutativity on the image of $\inl$ is the composition of the other 4 sides of the upper triangular prism; for $\inr$ we take the lower prism. We now need a dependent equality between these commutativities as we move along some $\glue$. By the construction of the dashed maps, it follows from the filler of the outer parallelepiped.

\begin{de}\label{de:cogap-nat}
    Denote this construction of a 2-cell by $\square$ (as a function of the outer parallelepiped).
\end{de}

For formal reasons, this construction respects juxtaposition:
\begin{lem}\label{le:cogap-juxt-commute}
    If $J$ denotes the operation of cube juxtaposition and $j$ denotes square juxtaposition, we have, for any juxtaposable cubes $a$ and $b$: $\square J(a,b) = j(\square a, \square b).$
\end{lem}

\subsection{Spheres}\label{sec:spheres}

Define $S^0 :\equiv \s$, pointed by $0_{\s}$, and $S^{1+n} :\equiv \Sigma S^{n}$, pointed by $N :\equiv \inl(*)$.
\\
Denote $S^{n,m} :\equiv S^n \wedge S^m$. We will tacitly transport along arithmetic equalities, identifying the corresponding spheres; for example we won't differentiate between $S^{k+t}$ and $S^{t+k}$. This is justified by the fact that all such equalities, and hence all transports along them, are equal.

We identify $S^{n,m}$ with $S^{n+m}$ via the recursively defined sequence of maps $j^{n,m}: S^{n+m} \pto S^{n,m}$ satisfying
\begin{align*}
    & j^{0,0}(0_{\s}) \equiv \langle 0_{\s}, 0_{\s}\rangle, j^{0,0}(1_{\s}) \equiv \langle 1_{\s}, 1_{\s} \rangle; \\
    & j^{0, m+1} \equiv \siin_2 \circ \Sigma j^{0, m}; \\
    & j^{n+1, m} \equiv \siin \circ \Sigma j^{n,m}.
\end{align*}\label{j}
The map $j^{n,m}$ is pointed by $\gluer_{0_{\s}}$ for $n=m=0$ and as a composition of pointed maps otherwise.

We identify $S^{n} \join S^{m}$ with $S^{n+m+1}$ via the recursively defined sequence of maps $e^{n,m}: S^{n+m+1} \pto S^{n} \join S^{m}$ satisfying
\begin{align*}
    & e^{0,0}(N) \equiv \inl (0_{\s}), e^{0,0}(S) \equiv \inl(1_{\s}),\\
    & \> \> \> \> \ap_{e^{0,0}}(\merid_{0_{\s}})=\glue(0_{\s},1_{\s})\cdot \glue(1_{\s},1_{\s})^{-1},\\
    & \> \> \> \> \ap_{e^{0,0}}(\merid_{1_{\s}})=\glue(0_{\s},0_{\s})\cdot \glue(1_{\s},0_{\s})^{-1}; \\
    & e^{0, m+1} \equiv - \siin_2 \circ \Sigma e^{0,m}; \\
    & e^{n+1,m} \equiv - \siin \circ \Sigma e^{n,m}.
\end{align*}
The map $e^{n,m}$ is pointed by $\refl_{\bs{}}$ for $n=m=0$ and as a composition of pointed maps otherwise.

The minus signs in this recursion denote precomposition with the \textit{flip map}. For $A : U_*$, let $\operatorname{flip}_A : \Sigma A \pto \Sigma A$ be defined by $N \mapsto S$, $S \mapsto N$ and $\ap_{\operatorname{flip}_A}(\merid_a) = \merid_a^{-1}$, pointed by $\merid_{\bs{A}}^{-1}$, and for a pointed map $f : \Sigma A \pto B$ set $-f :\equiv f \circ \operatorname{flip}_A$. Since $\ap_{\operatorname{flip}_A}(\mmerid_a) = \merid_{\bs{A}}^{-1} \cdot \mmerid_a^{-1} \cdot \merid_{\bs{A}}$, where the outer factors cancel against the pointing of $\operatorname{flip}_A$, we get $\adj{}(-f)(a) = \adj{}(f)(a)^{-1}$; that is, under the adjunction $\Sigma \dashv \Omega$, the map $-f$ corresponds to the pointwise inverse of $\adj{}(f)$. Accordingly, for a pointed map $g$ into a loop space we write $-g :\equiv \lambda a.\, g(a)^{-1}$. We also write $(-1)^k f$ for $f$ if $k$ is even and for $-f$ if $k$ is odd, and $-1 :\equiv -\id$. The flip map is natural: $\Sigma h \circ \operatorname{flip}_A =_* \operatorname{flip}_{A'} \circ \Sigma h$ for $h : A \pto A'$, so a minus sign can be moved past suspended maps; e.g. in the recursion above it does not matter whether it applies to $\siin$ or to the whole composite. Finally, under the identification $\pi_n(B) \simeq ||S^n \pto B||_0$ ($n \geq 1$), the operation $f \mapsto -f$ induces the inverse in the group $\pi_n(B)$.

\begin{lem}\label{le:sphere-identifications}
    The composition $S^{n+m+1} \xrightarrow{e^{n,m}} S^n \join S^m \xrightarrow{\jtos}
    \Sigma S^{n,m}$ is equal as a pointed map to $\Sigma j^{n,m}$.
\end{lem}
\begin{proof}
    We induct on $n$ and $m$.

    When $n=m=0,$ we equate the composition and $\Sigma j^{n,m}$ by $N \mapsto \refl_{N}; S \mapsto \merid_{\langle 1_{\s}, 1_{\s} \rangle}; \merid_{0_{\s}} \mapsto \ap_{\merid}(byContr) \cdot cancel; \merid_{1_{\s}} \mapsto  \ap_{\merid}(byContr) \cdot cancel$. This is trivially based.

    The inductive step follows from Lemma \ref{le:siin-compatible}.
\end{proof}

\begin{re}
This lemma explains the alternating signs for $e^{n,m}$: as $A \join B \simeq \Sigma (A \wedge B)$, taking $\Sigma$ inside from $\Sigma (A \join B)$ requires switching it with another $\Sigma$.
The recurrence for $j^{n,m}$ will be used in the proof of the Leibniz formula (Proposition \ref{4.1}). The recurrence for $e^{n,m}$ will give rise to the recursive definition of the pairing (Theorem \ref{th:defs-of-pairing}), which is the one we will use most often; this suggests alternating the signs in the definition of $j^{n,m}$.
However, for proper sign hygiene, we want $j^{n,m}$ to have a recurrence without alternating signs, as it preserves, intuitively, the order of coordinates. For example, otherwise Lemmas \ref{le:Sn-smash-merid} and \ref{le:merid-smash-Sm} would give minus the result one is used to seeing.
\end{re}

\begin{lem}
    $j^{n,m}$ and $e^{n,m}$ are equivalences for all $n,m$.
\end{lem}
\begin{proof}
    Follows inductively from the fact that $\siin$ ($\siin_2$) is an equivalence (Lemma \ref{le:sout-w}, Lemma \ref{le:sout-j}).
\end{proof}

\subsection{Relative homotopy groups}

\begin{de}
    Let $i_A: A \pto B$ be a pointed pair. For $n \geq 1$, define the set of ``relative spheroids'' $\pi_n(i_A)$ (or $\pi_n(B, A)$ when the map $i_A$ is clear from the context) to be the $0$-truncation of the type of diagrams in $U_*$ of the form

    \[\begin{tikzcd}
    	{S^{n-1}} & A \\
    	\sgtn & {B.}
    	\arrow["{\bar{x}}", from=1-1, to=1-2]
    	\arrow[from=1-1, to=2-1]
    	\arrow["\alpha"', shorten <=8pt, shorten >=8pt, Rightarrow, from=1-2, to=2-1]
    	\arrow["{i_A}", from=1-2, to=2-2]
    	\arrow["{\bs{B}}", from=2-1, to=2-2]
    \end{tikzcd}\]
\end{de}
Thus, $\pi_n(B, A)$ can be thought of as the set of elements of $\pi_{n-1}(A)$ with a trivialization (of their representatives) inside of $B$.

\begin{re} \label{addition}
    Let $F$ be the fiber of $i_A$ over the basepoint:

    \[\begin{tikzcd}
    	F & A \\
    	\sgtn & {B.}
    	\arrow[from=1-1, to=1-2]
    	\arrow[from=1-1, to=2-1]
    	\arrow["\lrcorner"{anchor=center, pos=0.125}, draw=none, from=1-1, to=2-2]
    	\arrow[shorten <=8pt, shorten >=8pt, Rightarrow, from=1-2, to=2-1]
    	\arrow["{i_A}", from=1-2, to=2-2]
    	\arrow["{\bs{B}}", from=2-1, to=2-2]
    \end{tikzcd}\]
    Then the universal property of the fiber implies that $\pi_n(B,A)$ is equivalently $\pi_{n-1}(F)$ (the induced map $S^{n-1} \to F$ is pointed because $\alpha$ is based). This leads to a group structure on $\pi_n(B,A)$ for $n > 1$.
    \\
    Translating it back to the level of diagrams,
    to add $(\bar{x}, \alpha)$ and $(\bar{y}, \beta),$ one has to form the diagram
    \[\begin{tikzcd}
    	{S^{n-1}} & {S^{n-1} \vee S^{n-1}} && A \\
    	\sgtn &&& B
    	\arrow["{\pin_{S^{n-2}}}", from=1-1, to=1-2]
    	\arrow[from=1-1, to=2-1]
    	\arrow["{(\bar{x}, \bar{y})}", from=1-2, to=1-4]
    	\arrow[between={0.3}{0.6}, Rightarrow, from=1-4, to=2-1]
    	\arrow[from=1-4, to=2-4]
    	\arrow[from=2-1, to=2-4]
    \end{tikzcd}\]
    with the obvious 2-cell. Here $\pin_{(-)}: \depfun{A:U_*}{(\Sigma A \pto  \Sigma A \vee  \Sigma A})$ is defined on $A: U_*$ by $N \mapsto \inr(S), S \mapsto \inl(S), \merid_a \mapsto \ap_{\inr}(\merid_a)^{-1} \cdot \gluer \cdot \gluel^{-1} \cdot \ap_{\inl}(\merid_a)$, pointed by $\ap_{\inr}(\merid_{\bs{A}})^{-1} \cdot \gluer$.
    We won't need this description.
\end{re}

\begin{re}
    $\pi_n(-)$ is natural: a map of pointed pairs $(B, A) \to (D, C)$ (that is, a commutative square in $U_*$) induces a map of relative spheroid sets $\pi_n(B,A) \to \pi_n(D,C)$, which is a homomorphism of groups for $n>1$.
\end{re}

\begin{re}
    To see a map $f: S^{n} \to A$ as a relative spheroid lifting an element of $\pi_{n}(A,\sgtn)$, we need to juxtapose with the pushout diagram, as in
    \[\begin{tikzcd}
    	{S^{n-1}} & \sgtn & \sgtn \\
    	\sgtn & {S^{n}} & {A.}
    	\arrow[from=1-1, to=1-2]
    	\arrow[from=1-1, to=2-1]
    	\arrow[from=1-2, to=1-3]
    	\arrow["\mmerid"{description}, between={0.3}{0.7}, Rightarrow, from=1-2, to=2-1]
    	\arrow[from=1-2, to=2-2]
    	\arrow["{\ptd{f}^{-1}}"{description}, between={0.3}{0.7}, Rightarrow, from=1-3, to=2-2]
    	\arrow[from=1-3, to=2-3]
    	\arrow[from=2-1, to=2-2]
    	\arrow["{\bs{A}}"', curve={height=12pt}, from=2-1, to=2-3]
    	\arrow["f", from=2-2, to=2-3]
    \end{tikzcd}\]
    The inverse equivalence $\pi_n(A, \sgtn) \to \pi_n(A)$ applies the universal property of the suspension to a diagram.
    \\
To see that this is a group homomorphism, it is enough to note that (in $U_*$) for $D$ the type of commutative diagrams
\begin{tikzcd}
	X & \sgtn \\
	\sgtn & A
	\arrow[from=1-1, to=1-2]
	\arrow[from=1-1, to=2-1]
	\arrow["{\bs{A}}"', from=1-2, to=2-2]
	\arrow["{\bs{A}}", from=2-1, to=2-2]
\end{tikzcd}, the diagram
\[\begin{tikzcd}
	& D \\
	{\Sigma X \to A} && {X \to \Omega A}
	\arrow["{\text{cogap}}"', from=1-2, to=2-1]
	\arrow["{\text{gap}}", from=1-2, to=2-3]
	\arrow["{\adj{}}", from=2-1, to=2-3]
\end{tikzcd}\]
commutes. Indeed, both maps take a 2-cell $\alpha$ to $\lambda x.\alpha(x)$.
\end{re}

\begin{re}\label{LES}
    Under this representation, the map $\pi_n(F) \to \pi_n(A)$ becomes the ``boundary homomorphism'' $\kappa: \pi _{n+1}(B,A) \to \pi_{n}(A),$ which picks out the top arrow $S^{n} \to A$ of a diagram. The long exact sequence

    $$\ldots \to \pi_n(F) \to \pi_n(A) \to \pi_n(B) \to \pi_{n-1}(F) \to \ldots$$
    takes the form
    $$\ldots \to \pi_{n+1}(B,A) \to \pi_n(A) \to \pi_n(B) \to \pi_n(B,A) \to \ldots,$$
    where the map $\pi_n(B) \to \pi_n(B,A)$ is the composition of $\pi_n(B) \xrightarrow{\simeq} \pi_n(B, \sgtn)$ from the above remark with the map induced by the map of pointed arrows $(B, \sgtn) \to (B, A)$. The fact that the square
    \begin{tikzcd}
    	{\pi_{n}(B, \sgtn)} & {\pi_{n}(B,A)} \\
    	{\pi_{n}(B)} & {\pi_{n-1}(F)}
    	\arrow[from=1-1, to=1-2]
    	\arrow[from=1-1, to=2-1]
    	\arrow[from=1-2, to=2-2]
    	\arrow[from=2-1, to=2-2]
    \end{tikzcd}
    commutes is not obvious (e.g. it could anticommute and both sequences would still be exact). See Lemma \ref{le:LES-isom} for its proof.
    As a warm-up we provide a direct proof of exactness, which is easier and self-contained:
\end{re}
\begin{prp}\label{prp:LESproof}
    Given a pointed pair $i: A \pto B$, the sequence of pointed sets
     $$\ldots \to \pi_n(A) \xrightarrow{\pi_n(i)} \pi_n(B) \xrightarrow{j_n} \pi_n(B,A) \xrightarrow{\kappa_n} \pi_{n-1}(A) \to \ldots$$
    is exact. Here the sequence ends at $n = 1$ (inclusive), the map $\kappa_n$ picks out the top arrow of a diagram, and the map $j_n$ is defined as $\pi_n(B) \xrightarrow{\sim} \pi_n(B,\sgtn) \to \pi_n(B,A)$, where the first arrow assigns to a pointed map out of $\Sigma S^{n-1}$ the corresponding cocone, and the second arrow is induced by a map of pointed pairs.
\end{prp}
\begin{proof}
    We prove exactness at every term of the sequence.
    \begin{enumerate}
        \item
        \begin{itemize}
            \item The composition $\pi_n(B) \to \pi_n(B, A) \to \pi_{n-1}(A)$ is obviously trivial.
            \item Assume given an element $|x|: \pi_n(B, A)$ such that its top map $\bar{x}: S^{n-1} \to A$ is based-homotopic to $\lambda \_.\bs{A}$. Then compose $x$ with the given homotopy, as in
            \[\begin{tikzcd}
            	{S^{n-1}} & A \\
            	\sgtn & {B.}
            	\arrow["{\bs{A}}", curve={height=-12pt}, from=1-1, to=1-2]
            	\arrow["{\bar{x}}"', from=1-1, to=1-2]
            	\arrow[from=1-1, to=2-1]
            	\arrow[from=1-2, to=2-2]
            	\arrow[from=2-1, to=2-2]
            \end{tikzcd}\]
            The outer square then exhibits a lift for $|x|$ to an element of $\pi_n(B, \sgtn) \simeq \pi_n(B)$.
        \end{itemize}
        \item
        \begin{itemize}
            \item The composition $\pi_n(B,A) \to \pi_{n-1}(A) \to \pi_{n-1}(B)$ is trivial by assigning to $|x|$ the 2-cell of $x$.
            \item Assume given an element $|f|$ of $\pi_{n-1}(A)$ such that $\pi_{n-1}(i) (|f|) = 0 : \pi_{n-1}(B)$. This means that $i \circ f$ is trivializable by a certain based homotopy $\alpha$ as a pointed map $S^{n-1} \pto B$. But then the class of
            \[\begin{tikzcd}
            	{S^{n-1}} & A \\
            	\sgtn & B
            	\arrow["f", from=1-1, to=1-2]
            	\arrow[from=1-1, to=2-1]
            	\arrow["\alpha"{description}, shorten <=4pt, shorten >=2pt, Rightarrow, from=1-2, to=2-1]
            	\arrow["i", from=1-2, to=2-2]
            	\arrow[from=2-1, to=2-2]
            \end{tikzcd}\]
            lifts $|f|$ to an element of $\pi_n(B,A)$.
        \end{itemize}
        \item
        \begin{itemize}
            \item The composition $\pi_n(A) \to \pi_n(B) \to \pi_n(B,A)$ factors through $\pi_n(A,A) = \pi_n(\id_A)$, so for its triviality it is enough to show that $\pi_n(A,A)$ is trivial. This follows from the contractibility of the type of commutative (in $U_*$) diagrams of the form
            \[\begin{tikzcd}
            	{S^{n-1}} & A \\
            	{\sgtn.}
            	\arrow[from=1-1, to=1-2]
            	\arrow[from=1-1, to=2-1]
            	\arrow[from=2-1, to=1-2]
            \end{tikzcd}\]
            To prove it, apply path induction twice:
            $\tot{f: S^{n-1} \to A}{\tot{p: f(\bs{}) = \bs{}}{\tot{h: f = \lambda \_.\bs{}}{(p^{-1} \cdot h(\bs{}) = \refl_{\bs{}})}}} \simeq \tot{p: \bs{} = \bs{}}{(p^{-1} = \refl_{\bs{}})} \simeq \sgtn$. (From a more abstract point of view, one instance of path induction is sufficient as the type in question is $\tot{f: \operatorname{Hom}_{U_*}(S^{n-1}, A)}{f= \lambda \_.\bs{A}}$.)
            \item Assume given a map $f: S^n \pto B$ such that the outer square in
            \[\begin{tikzcd}
            	{S^{n-1}} & \sgtn & A \\
            	\sgtn & {S^n} & B
            	\arrow[from=1-1, to=1-2]
            	\arrow[from=1-1, to=2-1]
            	\arrow[from=1-2, to=1-3]
            	\arrow["\merid"{description}, Rightarrow, from=1-2, to=2-1]
            	\arrow[from=1-2, to=2-2]
            	\arrow[from=1-3, to=2-3]
            	\arrow[from=2-1, to=2-2]
            	\arrow["f", from=2-2, to=2-3]
            \end{tikzcd}\]
            is equal to a trivial one.

            This means having the diagram
            \[\begin{tikzcd}
            	{S^{n-1}} & \sgtn & A \\
            	\sgtn & {S^n} & {B,}
            	\arrow[from=1-1, to=1-2]
            	\arrow[curve={height=-12pt}, from=1-1, to=1-3]
            	\arrow[from=1-1, to=2-1]
            	\arrow[from=1-2, to=1-3]
            	\arrow["\merid"{description}, Rightarrow, from=1-2, to=2-1]
            	\arrow[from=1-2, to=2-2]
            	\arrow["i", from=1-3, to=2-3]
            	\arrow[from=2-1, to=2-2]
            	\arrow[dashed, from=2-2, to=1-3]
            	\arrow["f", from=2-2, to=2-3]
            \end{tikzcd}\]
            where the top cell is some self-equality $\alpha$ of the trivial map such that the cells given by $\ap_i(\alpha)$ and $\ap_f(\merid)$ are equal. But this equality means that $f$ lifts to the map $S^{n} \pto A$ which is induced from $\alpha$ by the universal property of the suspension.
        \end{itemize}
    \end{enumerate}
\end{proof}

\begin{prp}\label{pr:NaturalityLES}
    The assignment $(A \to B) \mapsto (\ldots \to \pi_{n+1}(B,A) \to \pi_n(A) \to \pi_n(B) \to \pi_n(B,A) \to \ldots)$ is natural: if
    \[\begin{tikzcd}
    	A & {A'} \\
    	B & {B'}
    	\arrow["f", from=1-1, to=1-2]
    	\arrow[from=1-1, to=2-1]
    	\arrow[from=1-2, to=2-2]
    	\arrow["g"', from=2-1, to=2-2]
    \end{tikzcd}\]
    is a map of pointed pairs, then there is an induced map of long exact sequences.

\end{prp}
\begin{proof}
    The commutativity of the arising squares follows from the naturality of $\pi_*$ for all maps except the boundary maps $\pi_n(B,A) \to \pi_{n-1}(A)$. For the latter, the commutativity says that juxtaposing two squares and taking the top map is the same as taking the top maps of the squares and then composing them.
\end{proof}

\section{Pairings}\label{se:pairings}

We now define a product operation on relative homotopy groups. The first definition is inspired by the corresponding construction in classical algebraic topology, see \cite{Dug1}, Section 4.

\subsection{On homotopy groups}

\subsubsection{Pairing via cogap maps}

Let $n,m \geq 1$ and assume $f: A \pto B$, $g: C \pto D$ are pointed pairs. Let $x :  \pi_n(B,A) \equiv \pi_n(f)$ and $y :  \pi_m(D,C) \equiv \pi_m(g).$ We will define an element of a set, so we can assume we have diagrams
\[\begin{tikzcd}
	{S^{n-1}} & A & {S^{m-1}} & C \\
	\sgtn & B & \sgtn & {D.}
	\arrow["{\bar{x}}", from=1-1, to=1-2]
	\arrow[from=1-1, to=2-1]
	\arrow["{\alpha_x}"', shorten <=6pt, shorten >=6pt, Rightarrow, from=1-2, to=2-1]
	\arrow[from=1-2, to=2-2]
	\arrow["{\bar{y}}", from=1-3, to=1-4]
	\arrow[from=1-3, to=2-3]
	\arrow["{\alpha_y}"', shorten <=6pt, shorten >=6pt, Rightarrow, from=1-4, to=2-3]
	\arrow[from=1-4, to=2-4]
	\arrow["{\bs{B}}", from=2-1, to=2-2]
	\arrow["{\bs{D}}", from=2-3, to=2-4]
\end{tikzcd}\]
Then let $P : \equiv {(A \wedge D) \coprod^{(A \wedge C)} (B \wedge C)}$ and let $P^\wedge(f,g) : P \pto B \wedge D$ be the pushout smash of $A \xrightarrow{f} B$ and $C \xrightarrow{g} D$ (Definition \ref{de:PO-smash}).

\begin{de}\label{de:circ}
    Define $x \circ y : \pi_{n+m}(B \wedge D, P)$ as the element represented by
\[\begin{tikzcd}
	{S^{n+m-1}} & {\Sigma (S^{n-1} \wedge S^{m-1}) \simeq (S^{n-1} \wedge \sgtn) \coprod^{(S^{n-1} \wedge S^{m-1})} (\sgtn \wedge S^{m-1})} & P \\
	\sgtn & \sgtn & {B \wedge D.}
	\arrow["{\Sigma j^{n-1,m-1}}", from=1-1, to=1-2]
	\arrow[from=1-1, to=2-1]
	\arrow["{(x \wedge \bs{}) \sqcup (\bs{} \wedge y)}", from=1-2, to=1-3]
	\arrow[from=1-2, to=2-2]
	\arrow["{\alpha_{x \circ y}}"{description}, shorten <=20pt, shorten >=20pt, Rightarrow, from=1-3, to=2-2]
	\arrow[from=1-3, to=2-3]
	\arrow[from=2-1, to=2-2]
	\arrow[from=2-2, to=2-3]
\end{tikzcd}\]
Here $j^{n-1, m-1}$ is from Section \ref{sec:spheres}, the equivalence comes from the fact that for any $X$ the types $X \wedge \sgtn$ and $\sgtn \wedge X$ are contractible (we take $\auxr$ as the center), and $x \wedge \bs{} \sqcup \bs{} \wedge y$ is induced from a map of diagrams (using Lemma \ref{le:smashes-cube}, see below). Finally, the 2-cell is constructed via the naturality of the pushout smash as a cogap map (see the construction of Definition \ref{de:cogap-nat}). For this we need a filler of the outer parallelepiped in the following diagram:

\[\begin{tikzcd}
	&& {A \wedge C} & {A \wedge D} \\
	{S^{n-1,m-1}} & {S^{n-1} \wedge \sgtn} & {B \wedge C} & \bullet \\
	{\sgtn \wedge S^{m-1}} & \bullet &&& {B \wedge D} \\
	&& {\sgtn \wedge \sgtn.}
	\arrow[from=1-3, to=1-4]
	\arrow[from=1-3, to=2-3]
	\arrow["\inl"', from=1-4, to=2-4]
	\arrow[from=1-4, to=3-5]
	\arrow[from=2-1, to=1-3]
	\arrow[from=2-1, to=2-2]
	\arrow[from=2-1, to=3-1]
	\arrow[from=2-2, to=1-4]
	\arrow[from=2-2, to=3-2]
	\arrow[from=2-2, to=4-3]
	\arrow["\inr", from=2-3, to=2-4]
	\arrow[from=2-3, to=3-5]
	\arrow["\lrcorner"{anchor=center, pos=0.125, rotate=180}, draw=none, from=2-4, to=1-3]
	\arrow[dashed, from=2-4, to=3-5]
	\arrow["{\alpha_{x \circ y}}", shorten <=12pt, shorten >=12pt, Rightarrow, from=2-4, to=4-3]
	\arrow[from=3-1, to=2-3]
	\arrow[from=3-1, to=3-2]
	\arrow[from=3-1, to=4-3]
	\arrow["\lrcorner"{anchor=center, pos=0.125, rotate=180}, draw=none, from=3-2, to=2-1]
	\arrow[dashed, from=3-2, to=2-4]
	\arrow[dashed, from=3-2, to=4-3]
	\arrow[from=4-3, to=3-5]
\end{tikzcd}\]
It exists by Lemma \ref{le:smashes-cube}. The constructed 2-cell is based: indeed, its composition with $\inl$ is based as a juxtaposition of based cells (see Corollary \ref{cor:comp_of_based}), which is enough because $\inl$ is pointed by $\refl$.
\end{de}

\begin{de}
    Define the $\times$-pairing $\odot: \pi_n(f) \times \pi_m(g) \to \pi_{n+m}(P^\times(f,g))$ in the same way but using $\times$ instead of $\wedge$ and the map $e^{n-1,m-1}$
    instead of $\Sigma j^{n-1, m-1}$. That is, $x \odot y$ is represented by the diagram

\[\begin{tikzcd}
	{S^{n+m-1}} & {S^{n-1} \join S^{m-1} \simeq (S^{n-1} \times \sgtn) \coprod^{(S^{n-1} \times S^{m-1})} (\sgtn \times S^{m-1})} & P \\
	\sgtn & \sgtn & {B \times D,}
	\arrow["{e^{n-1,m-1}}", from=1-1, to=1-2]
	\arrow[from=1-1, to=2-1]
	\arrow["{(x \times \bs{}) \sqcup (\bs{} \times y)}", from=1-2, to=1-3]
	\arrow[from=1-2, to=2-2]
	\arrow["{\alpha_{x \odot y}}"{description}, shorten <=19pt, shorten >=19pt, Rightarrow, from=1-3, to=2-2]
	\arrow[from=1-3, to=2-3]
	\arrow[from=2-1, to=2-2]
	\arrow[from=2-2, to=2-3]
\end{tikzcd}\]
where $P :\equiv A \times D \coprod^{A \times C} B \times C$ and $\alpha_{x \odot y}$ is constructed via $\square$ like above.

\end{de}

\subsubsection{On the definitions of the pairing}

First of all, we show that the pairing $\circ$ factors through the cross-product pairing $\odot$ via the projection $P^\times \to P^\wedge$.

Note that there is a map of pointed pairs $\ctos: (B \times D, {(A \times D) \coprod^{A \times C} (B \times C)}) \pto (B \wedge D, {(A \wedge D) \coprod^{A \wedge C} (B \wedge C)})$. To define it, first form the parallelepiped with 4 sides being the maps $\ctos$ (directly or using Lemma \ref{le:nat-tr-cube}) and then take $\square$ of it.
\begin{prp}\label{pr:cross-factors-wedge}
    Let $(B, A)$, $(D, C)$ be pointed pairs.
    Then, for any $n,m \ge 0$ and $x : \pi_{n+1}(B,A), y : \pi_{m+1}(D,C)$, the following relation holds:
    $$x \circ y = \ctos_* (x \odot y).$$
\end{prp}
\begin{proof}
    Suppose given elements $|x| : \pi_{n+1}(B,A), |y| : \pi_{m+1}(D,C)$. Upon pre-pasting the diagrams defining $|x| \circ |y|$ and $\ctos_* (|x| \odot |y|)$ with $(e^{n,m}, \id_{\sgtn})^{-1}$ and applying Lemma \ref{le:sphere-identifications}, they both become represented by $\square$ of two parallelepipeds with the same left corners ($S^n \xleftarrow{} S^n \times S^m \xrightarrow{} S^m$) and right corners ($A \wedge D \xleftarrow{} A \wedge C \xrightarrow{} B \wedge C$). Both are juxtapositions of two smaller parallelepipeds: one is ``first $\ctos$ on spheres and then map to spaces'', and the other is ``first map to spaces and then $\ctos$ on them''. We need to prove that these parallelepipeds are equal.

    They both start at the product $S^n \times S^m$. Thus, even when dealing with the fillers from Lemma \ref{le:smashes-cube}, we only consider their values on elements of the form $\langle x, y \rangle$, which compute to reflexivity. This makes the check purely formal.
\end{proof}

Now note that if $F_f$ and $F_g$ are the fibers of $f$ and $g$, then the fiber of $P^{\times}(f,g)$ is equivalent to $F_f \join F_g$ (Lemma \ref{le:fiber-of-pp}). Thus, since relative homotopy groups are absolute homotopy groups of fibers, one can ask if the pairing has a corresponding description.

\begin{thm}\label{th:defs-of-pairing}

Let $n,m > 0$ and let $f: A \pto B$ and $g: C \pto D$ be pointed pairs. The $\times$-pairing $\pi_{n}(f) \times \pi_{m}(g) \xrightarrow{\odot} \pi_{n+m}(P^{\times}(f,g))$ is equivalently the operation $\pi_{n-1}(F_f) \times \pi_{m-1}(F_g) \to \pi_{n+m-1}(F_f \join F_g)$ obtained by truncating from $I_{n-1,m-1}^{F_f, F_g}$, where $I_{k,t}^{X,Y}: \Omega^k X \times \Omega^t Y \to \Omega^{k+t+1}(X \join Y)$ is the following inductively defined function:
\begin{itemize}
    \item
for $k=t=0$, $I^{X,Y}_{k,t}: X \times Y \to \Omega (X \join Y)$ equals $\lambda(x,y).\glue(\bs{X},\bs{Y})\glue(x,\bs{Y})^{-1}\glue(x,y)\glue(\bs{X},y)^{-1}$;
    \item
for $k>0$, the map $I_{k,t}^{X,Y}$ is $\Omega^k X \times \Omega^t Y \xrightarrow{-I^{\Omega X,Y}_{k-1,t}} \Omega^{k+t} (\Omega X \join Y) \xrightarrow{\Omega^{k+t}\oout} \Omega^{k+t+1}(X \join Y)$;
    \item
for $k=0, t>0$, use $\oout$ for the second factor: $X \times \Omega^t Y \xrightarrow{-I^{X,\Omega Y}_{0,t-1}} \Omega^{t} (X \join \Omega Y) \xrightarrow{\Omega^{t}\oout_2} \Omega^{t+1}(X \join Y)$.
\end{itemize}
\end{thm}

\begin{proof}

Let $F_f$, $F_g$ be the fibers of $f$ and $g$. There are maps of pairs $(\sgtn, F_f) \to (B, A)$ and $(\sgtn, F_g) \to (D, C)$. Note that the equivalence of Lemma \ref{le:fiber-of-pp} for maps $f, g$ factors as $F_f \join F_g \xrightarrow{\sim} \fib P^\times(F_f \to \sgtn, F_g \to \sgtn) \to \fib P^\times(f,g)$. Moreover, the $\times$-pairing $\odot$ is natural (Proposition {\ref{odot-prod-nat}}). This allows us to assume that the codomains of $f$ and $g$ are $\sgtn$. For $\pi_-(\sgtn, -)$ the operation $\odot$ works via ``joining'' the representatives.
Thus, we need to show that the inductively defined operation $I^{X, Y}_{n,m}: \Omega^n X \times \Omega^m Y \to \Omega^{n+m+1}(X \join Y)$ takes $(x, y): \Omega^n X \times \Omega^m Y$ adjunct to $\bar{x}: S^n \pto X, \bar{y}: S^m \pto Y $ to the element adjunct to $S^{n+m+1} \xrightarrow{e^{n,m}} S^n \join S^m \xrightarrow{\bar{x} \join \bar{y}} X \join Y$.
\begin{enumerate}
    \item The base case $n = m = 0$ follows from the definitions of $I^{X, Y}_{0,0}$ and $e^{0,0}$: for $\loo : \equiv \mmerid_{1_\s}$, we have $\ap_{((\bar{x} \join \bar{y}) \circ e^{0,0})}(\loo) = \glue(\bar{x}(0_{\s}), \bar{y}(0_{\s}))\glue(\bar{x}(1_{\s}), \bar{y}(0_{\s}))^{-1}(\glue(\bar{x}(0_{\s}),\bar{y}(1_{\s}))\glue(\bar{x}(1_{\s}),\bar{y}(1_{\s}))^{-1})^{-1} =
    \glue(\bs{X},\bs{Y})\glue(x ,\bs{Y})^{-1}(\glue(\bs{X},y)\glue(x, y)^{-1})^{-1} = I^{X, Y}_{0,0}(x,y).$

    \item Let us do the inductive step, for example, on $n$: let $n>0$ and $m$ arbitrary, and consider
    $$I_{n,m}^{X,Y} \equiv (\Omega^n X \times \Omega^m Y \xrightarrow{-I^{\Omega X, Y}_{n-1, m}} \Omega^{n+m}(\Omega X \join Y) \xrightarrow{\Omega^{n+m} (\Omega\operatorname{-out})} \Omega^{n+m+1}(X \join Y)).$$
    By induction, $\adj{n+m}(I^{\Omega X, Y}_{n-1, m} (x,y)) =_* S^{n+m} \xrightarrow{e^{n-1,m}} S^{n-1} \join S^m \xrightarrow{\adj{n-1}x \join \adj{m} y} \Omega X \join Y$, and so $\adj{n+m}$ will take the composite above to
    $$S^{n+m} \xrightarrow{-e^{n-1,m}} S^{n-1} \join S^m \xrightarrow{\adj{n-1}x \join \adj{m}y} \Omega X \join Y \xrightarrow{\Omega\operatorname{-out}} \Omega^{}(X \join Y).$$
    Now, $\adj{}$ of this has to be shown to equal $S^{n+m+1} \xrightarrow{e^{n,m}} S^{n} \join S^m \xrightarrow{\adj{n}x \join \adj{m} y} X \join Y$.

    Denote $S^{n-1} $ by $A$, $S^m$ by $B$, $\adj{n-1}x$ by $\bar{x}$ and $\adj{m}y$ by $\bar{y}$. Then we have to prove that the outer diagram commutes in
    \[\begin{tikzcd}
    	& {\Sigma A \join B} & {\Sigma \Omega X \join Y} & {X \join Y} \\
    	{S^{m+n+1}} & {\Sigma (A \join B)} & {\Sigma (\Omega X \join Y)} & {\Sigma \Omega (X \join Y).}
    	\arrow["{\Sigma \bar{x} \join \bar{y}}", from=1-2, to=1-3]
    	\arrow["{\operatorname{adj}\bar{x} \join \bar{y}}", curve={height=-12pt}, from=1-2, to=1-4]
    	\arrow["{\epsilon_X \join Y}", from=1-3, to=1-4]
    	\arrow["{\Sigma e^{n-1,m}}", from=2-1, to=2-2]
    	\arrow["{\siin }", from=2-2, to=1-2]
    	\arrow["{\Sigma (\bar{x} \join \bar{y})}"', from=2-2, to=2-3]
    	\arrow["{\siin }", from=2-3, to=1-3]
    	\arrow["{\Sigma \oout}"', from=2-3, to=2-4]
    	\arrow["{\epsilon_{X \join Y}}"', from=2-4, to=1-4]
    \end{tikzcd}\]
    (Here $\epsilon_{-}$ is the counit of the $\Sigma - \Omega$ adjunction.) The left square commutes by the naturality of $\siin$ (Lemma \ref{le:siin-nat-j}). For the commutativity of the right square, see Lemma \ref{le:siin-oout-adj}.
\end{enumerate}
\end{proof}

\begin{re}
    The minus sign in the definition of $I$ is the minus sign from the definition of $e$, as seen in the proof. We can think of it as expressing that the $\Omega$ appearing from $I_{0,0}$ has to be pushed inside after taking out all other loop space functors for the ``correct'' order. This adds the sign $(-1)^{n+m}$.
\end{re}

\begin{cor}\label{cor:distr1}
    The operation $\odot$ is distributive over addition in $\pi_*$: for $n \geq 2$, $x, x' : \pi_n(B,A), y, y' : \pi_m(D,C),$ we have $(x+x') \odot y = x \odot y + x' \odot y$ and for $m \geq 2,$ $x \odot (y + y') = x \odot y + x \odot y'$.
\end{cor}
\begin{proof}
    Distributivity holds on the space level already (before 0-truncation).

    Let us first show it over the first factor for the case $n=2, m=1$. We need
    a term of the type
    $$\depfun{X,Y:U_*}{\depfun{p,q:\Omega X}{I_{1,0}^{X,Y}}(p \cdot q, -) = I_{1,0}^{X,Y}(p, -) \cdot I_{1,0}^{X,Y}(q, -)}.$$
    This follows from the explicit formula for $I_{1,0}^{X,Y}$. Concretely, denote the wedge-shaped path $\glue (x,y) \cdot \glue (x',y)^{-1}$ by $w(x,x',y): \inl(x) = \inl(x')$. Let $p:x=x', q: x' = x''$. In the diagram
    \[\begin{tikzcd}
    	& {\refl_{\bs{}}} \\
    	{\ap_{\inl}(p) w(x',x'',\bs{Y})} & {w(x,x',\bs{Y}) w(x',x'',\bs{Y})} & {w(x,x'',\bs{Y})} \\
    	{w(x,x',y) \ap_{\inl}(q)} & {w(x,x',y) w(x',x'',\bs{Y})} & {\ap_{\inl}(p \cdot q)} \\
    	& {w(x,x',y)w(x',x'',y)} & {w(x,x'',y)} \\
    	& {\refl_{\bs{}}}
    	\arrow[dashed, from=1-2, to=2-2]
    	\arrow["{\text{if }x \equiv x''}", dashed, from=1-2, to=2-3]
    	\arrow[squiggly, from=2-1, to=3-2]
    	\arrow[squiggly, from=2-2, to=2-1]
    	\arrow[squiggly, from=2-2, to=2-3]
    	\arrow[squiggly, from=2-3, to=3-3]
    	\arrow[squiggly, from=3-1, to=4-2]
    	\arrow[squiggly, from=3-2, to=3-1]
    	\arrow[squiggly, from=3-3, to=4-3]
    	\arrow[dashed, from=4-2, to=5-2]
    	\arrow[squiggly, from=4-3, to=4-2]
    	\arrow[dashed, from=4-3, to=5-2]
    \end{tikzcd}\]
    the left and the right paths $w(x,x',\bs{Y}) w(x',x'',\bs{Y}) \rightsquigarrow w(x,x',y)w(x',x'',y)$ are equal by path induction. When $x \equiv x' \equiv x'' \equiv \bs{X}$, the right loop on $\refl$ equals $I_{1,0}^{X,Y}(pq, y)$. Moving from $\refl$ to the left, the first three arrows composed with $w(x,x',y) \rightsquigarrow \refl$ would be $I_{1,0}^{X,Y}(p, y)$ if they didn't involve $w(x',x'',\bs{Y})$. The second factor in the composite $I_{1,0}^{X,Y}(p, y) \cdot I_{1,0}^{X,Y}(q, y)$ then proceeds to de-contract $\refl$ to $w(x',x'',\bs{Y})$, and move along the next two arrows as if they were not whiskered with $w(x,x',y)$. Elementary path algebra equalities witnessing that it does not matter in which order to cancel the inverses and remove $\refl$ then imply that the left loop on $\refl$ equals $I_{1,0}^{X,Y}(p, y) \cdot I_{1,0}^{X,Y}(q, y)$.

    Similarly, we get $I_{0,1}^{X,Y}(-,p \cdot q) = I_{0,1}^{X,Y}(-,p) \cdot I_{0,1}^{X,Y}(-,q)$ for any $p,q: \Omega Y$.
    Now, the distributivity over the second factor follows by induction.
    Let $m>0, n>0$ and $p,q: \Omega^m Y, r: \Omega^n X$. By induction, $I^{\Omega X,Y}_{n-1,m}(r, p \cdot q) = I^{\Omega X,Y}_{n-1,m}(r, p) \cdot I^{\Omega X,Y}_{n-1,m}(r, q)$. The map $-\Omega^{n+m} \oout$ then takes it to an element equal to $ I^{X,Y}_{n,m}(r,p) \cdot  I^{X,Y}_{n,m}(r, q)$, as desired. The case $n=0, m>1$ follows similarly.

    Distributivity over the first factor now follows from graded-commutativity (Proposition \ref{pr:graded-comm}).
\end{proof}

\begin{cor}\label{cor:distr}
    The operation $\circ$ is distributive over addition in $\pi_*$: for $n \geq 2$, $x, x' : \pi_n(B,A), y, y' : \pi_m(D,C),$ we have $(x+x') \circ y = x \circ y + x' \circ y$ and for $m \geq 2,$ $x \circ (y + y') = x \circ y + x \circ y'$.
\end{cor}

\begin{prp}\label{odot-prod-nat}
    $\odot$ is natural: suppose given a pointed pair $g: C \pto D$ and a map between pointed pairs $f: A \pto B$ and $f': A' \pto B'$. This gives rise to maps $\pi_n(B,A) \to \pi_n(B',A')$ and
    $\pi_{n+m}(P^{\times}(f,g)) \to \pi_{n+m}(P^{\times}(f',g))$. If the former map sends $x$ to $x',$ then the latter, for any $y : \pi_m(D,C)$, sends $x \odot y$ to $x' \odot y$.
\end{prp}
\begin{proof}
    We are in a situation of having the diagram
    \[\begin{tikzcd}
    	&&&& \bullet & \bullet \\
    	&& \bullet & \bullet & \bullet & \bullet \\
    	\bullet & \bullet & \bullet & \bullet &&& \bullet \\
    	\bullet & \bullet &&& \bullet \\
    	&& {\bullet.}
    	\arrow[from=1-5, to=1-6]
    	\arrow[from=1-5, to=2-5]
    	\arrow[from=1-6, to=2-6]
    	\arrow[from=1-6, to=3-7]
    	\arrow[from=2-3, to=1-5]
    	\arrow[from=2-3, to=2-4]
    	\arrow[from=2-3, to=3-3]
    	\arrow[from=2-4, to=1-6]
    	\arrow["\inl"', from=2-4, to=3-4]
    	\arrow[from=2-4, to=4-5]
    	\arrow[from=2-5, to=2-6]
    	\arrow[from=2-5, to=3-7]
    	\arrow["\lrcorner"{anchor=center, pos=0.125, rotate=180}, draw=none, from=2-6, to=1-5]
    	\arrow[dashed, from=2-6, to=3-7]
    	\arrow[from=3-1, to=2-3]
    	\arrow[from=3-1, to=3-2]
    	\arrow[from=3-1, to=4-1]
    	\arrow[from=3-2, to=2-4]
    	\arrow[from=3-2, to=4-2]
    	\arrow[from=3-2, to=5-3]
    	\arrow[from=3-3, to=2-5]
    	\arrow["\inr", from=3-3, to=3-4]
    	\arrow[from=3-3, to=4-5]
    	\arrow["\lrcorner"{anchor=center, pos=0.125, rotate=180}, draw=none, from=3-4, to=2-3]
    	\arrow[dashed, from=3-4, to=2-6]
    	\arrow[dashed, from=3-4, to=4-5]
    	\arrow[from=4-1, to=3-3]
    	\arrow[from=4-1, to=4-2]
    	\arrow[from=4-1, to=5-3]
    	\arrow["\lrcorner"{anchor=center, pos=0.125, rotate=180}, draw=none, from=4-2, to=3-1]
    	\arrow[dashed, from=4-2, to=3-4]
    	\arrow[dashed, from=4-2, to=5-3]
    	\arrow[from=4-5, to=3-7]
    	\arrow[from=5-3, to=4-5]
    \end{tikzcd}\]

The claim follows immediately from the fact that juxtaposition and $\square$ commute (Lemma \ref{le:cogap-juxt-commute}). (Proving that the juxtaposition of the smaller parallelepipeds equals the outer parallelepiped is straightforward.\footnote{This is one of the main reasons for introducing the $\times$-pairing $\odot$. Such a proof would be very technical if the parallelepipeds were those with $\wedge$.})
\end{proof}

\begin{lem}\label{le:cross-to-wedge-natural}
    $P^\times \to P^\wedge$ is natural: suppose given an arrow $g$ and a map of arrows $f \to f'$. Then the diagram of arrow maps below commutes (meaning the two cells are equal in $U_*$).
    \[\begin{tikzcd}
    	{P^\times(f,g)} & {P^\wedge(f,g)} \\
    	{P^\times(f',g)} & {P^\wedge(f',g)}
    	\arrow[from=1-1, to=1-2]
    	\arrow[from=1-1, to=2-1]
    	\arrow[from=1-2, to=2-2]
    	\arrow[from=2-1, to=2-2]
    \end{tikzcd}\]
\end{lem}
\begin{proof}
    Straightforward.
\end{proof}

\begin{prp}\label{circ-prod-nat}
    $\circ$ is natural: suppose given a pointed pair $g: C \pto D$ and a map between pointed pairs $f: A \pto B$ and $f': A' \pto B'$. This gives rise to maps $\pi_n(B,A) \to \pi_n(B',A')$ and
    $\pi_{n+m}(P^{\wedge}(f,g)) \to \pi_{n+m}(P^{\wedge}(f',g))$. If the former map sends $x$ to $x',$ then the latter, for any $y : \pi_m(D,C)$, sends $x \circ y$ to $x' \circ y$.
\end{prp}
\begin{proof}
Immediate from the naturality of the map $P^\times \to P^\wedge$ (Lemma \ref{le:cross-to-wedge-natural}) and the naturality of $\odot$ (Proposition \ref{odot-prod-nat}).
\end{proof}

\begin{prp}\label{pr:graded-comm}
    Let $x: \pi_m(f), y: \pi_n(g)$. Consider the obvious map $\mir: P^\wedge(f, g) \to P^\wedge(g,f)$ (we can define it as $\square$ of the cube from Lemma \ref{le:nat-tr-cube} applied to the natural transformation $\swap$ that interchanges the left and right factors of a smash-product). We have $\mir_*(x \circ y) = (-1)^{nm+n+m} y \circ x$.
\end{prp}
\begin{proof}
    By naturality, it is enough to check the universal case, when $x$ and $y$ correspond to identity maps. By Lemma \ref{le:siin-anticom}, $\siin$ and $\siin_2$ anticommute, from which it follows inductively that $\swap \circ j^{p,q} = (-1)^{pq} j^{q,p}$ for any $p,q$ (see Corollary \ref{cor:oout-anticom-w} and Lemma \ref{le:oout-mirrors-oout2-w} for the exact same argument with $\oout$). Then we get $\mir \circ \Sigma j^{n-1,m-1} = (-1)^{(n-1)(m-1)+1}\Sigma j^{m-1,n-1}$ and the claim follows.
\end{proof}

\begin{re}
    An alternative way to prove the graded-commutativity would be to use the absolute homotopy group description from Theorem \ref{th:defs-of-pairing}. Indeed, the case $n=m=1$ corresponds to $I_{0,0}^{F_f,F_g}$, which differs from $I_{0,0}^{F_g,F_f}$ by a minus sign because of the inversion in the relation $\ap_{\mir}\glue(x,y)=_{\inr(x) = \inl(y)}\glue({y,x)^{-1}}$. The cases of higher dimension then go like above: swapping $\ooutk{m-1}$ and $\ooutk{n-1}$ adds $(n-1)(m-1)$ to the exponent of the sign. However, the proof of the base case for the join version of Corollary \ref{cor:oout-anticom-w} seems to be much less straightforward than its smash counterpart (see Lemma \ref{le:oout-anticom-j}).
\end{re}

\begin{de}\label{de:ooutk}
    For $X, Y: U_*$, define inductively $\ooutk{k}: \Omega^k X \wedge Y \pto \Omega^k (X \wedge Y)$ via
    \begin{align*}
        & \ooutk{0} \equiv \id; \\
        & \ooutk{k+1} \equiv \Omega^k \oout \circ \ooutk{k}.
    \end{align*}
\end{de}
\noindent
Define $\ooutk{k}_2$ similarly.

Note that $P^\wedge(\sgtn \to X, \sgtn \to Y) \simeq (X \wedge Y, \sgtn)$, and so the product $\circ$ defines a map $\pi_n(X,\sgtn) \times \pi_m(Y,\sgtn) \to \pi_{n+m}(X \wedge Y,\sgtn)$.
\begin{prp}\label{pr:smashing}
    For $n,m \geq 1$, the map $\pi_n(X,\sgtn) \times \pi_m(Y,\sgtn) \to \pi_{n+m}(X \wedge Y,\sgtn)$
    under the isomorphism $\pi_\bullet(-, \sgtn) \simeq \pi_{\bullet}(-)$
    equals the truncation of
    $$\Omega^n X \times \Omega^m Y \to \Omega^n X \wedge \Omega^m Y \xrightarrow{\ooutk{m}_2} \Omega^m (\Omega^n X \wedge Y) \xrightarrow{\Omega^m \ooutk{n}} \Omega^{n+m} (X \wedge Y)$$
    up to sign $(-1)^n$.
\end{prp}

\begin{proof}
    The proof is technical and is deferred to the last section, see Proposition \ref{pr:smashing-proof}.
\end{proof}

So far, the pairing $\circ$ has been defined on relative homotopy groups of dimension $\geq 1$, but the operation defined above makes sense for $0$-th homotopy groups as well. We can thus extend the definition.
\begin{de}\label{de:circ-absolute}
    Define the \textit{absolute homotopy group pairing}
    $$\circ: \pi_n(X) \times \pi_m(Y) \to \pi_{n+m}(X \wedge Y)$$
    to be the pairing from Proposition \ref{pr:smashing}. That is, $\circ$ is the set-truncation of the recursively defined map
    $$\Omega^n X \times \Omega^m Y \to \Omega^n X \wedge \Omega^m Y \xrightarrow{\oout_2} \Omega (\Omega^{n} X \wedge \Omega^{m-1} Y) \to \ldots \xrightarrow{\Omega^{n+m-1}\oout} \Omega^{n+m}(X \wedge Y).$$
\end{de}
\noindent
This will make some diagrams from Section \ref{se:SSS} slightly cuter, but otherwise we won't make essential use of this definition.

In conclusion, we have one pairing to show the distributivity, one for the naturality and one (as we'll show now) for the Leibniz formula. The first two are equivalent and the last one is their image under the projection $\ctos: P^\times \to P^\wedge$.

\subsubsection{The Leibniz formula}

We now investigate the Leibniz formula for $\circ$. Let $(B,A)$ and $(D,C)$ be pointed pairs and set $P :\equiv (A \wedge D) \coprod^{(A \wedge C)} (B \wedge C)$. Given $x : \pi_{n+1}(B,A), y : \pi_{m+1}(D,C)$ (we shift the degrees by one, as this simplifies the indexing of the spheres below), one can take the boundary $\kappa x : \pi_{n}(A) \simeq \pi_{n}(A,\sgtn)$ and then multiplying by $y$ gives an element $(\kappa x) \circ y : \pi_{n+m+1}(A \wedge D, (\sgtn \wedge D) \coprod^{\sgtn \wedge C} (A \wedge C) ).$
Denote by $l$ the obvious pointed pair map $(A \wedge D, (\sgtn \wedge D) \coprod^{\sgtn \wedge C} (A \wedge C) ) \xrightarrow{\sim} (A \wedge D, A \wedge C) \to (P, A \wedge C)$ where the pair on the right is $A \wedge C \to A \wedge D \xrightarrow{\inl} P$.
Similarly, denote by $r$ the map $(B \wedge C,  (A \wedge C) \coprod^{A \wedge \sgtn} (B \wedge \sgtn)) \xrightarrow{\sim} (B \wedge C, A \wedge C) \to (P, A \wedge C)$. Note that the square $(A \wedge D, A \wedge C) \to (P, A \wedge C)$ commutes by $\refl$ while $(B \wedge C, A \wedge C) \to (P, A \wedge C)$ commutes by $\glue$. There is also the map $j: (P, \sgtn) \to (P, A \wedge C)$.

\begin{prp}\label{4.1}
    With the above notation,
    $j_* \kappa (x \circ y) = -l_* ((\kappa x) \circ y) + (-1)^{n+1} r_* (x \circ \kappa y ): \pi_{n+m+1}(P, A \wedge C)$ for $n,m>0.$
\end{prp}

\begin{re}
    Recall that $x$ and $y$ have degrees $n+1$ and $m+1$. The formula \textit{a priori} can't look like the standard $d(xy) = (dx)y+(-1)^{n+1}xdy$ as long as the skew-commutativity coefficient is $(-1)^{nm+1}$ (Proposition \ref{pr:graded-comm}), or we'll get $d(xy) = -d(xy)$. We can, however, make it look like $d(xy) = (dx)y+(-1)^{n}xdy$ if we twist the definition of the product by $(-1)^{n+m}$. For example, this would be the case if we chose the recursive formula for $j^{n,m}$ instead of $e^{n,m}$ to have sign alternation.
\end{re}

Before proving the formula, let us record a lemma.

\begin{lem}\label{le:pi(S,S)}
    For $n > 1,$ there is an isomorphism $\pi_{n}(S^{n}, S^{n-1}) \simeq \ZZ\times \ZZ$, where $S^{n-1} \to S^n$ is the constant map. To a diagram
    \[\begin{tikzcd}
    	{S^{n-1}} & {S^{n-1}} \\
    	\sgtn & {S^{n}}
    	\arrow["f", from=1-1, to=1-2]
    	\arrow[from=1-1, to=2-1]
    	\arrow["\alpha"', Rightarrow, nfold, from=1-2, to=2-1]
    	\arrow[from=1-2, to=2-2]
    	\arrow[from=2-1, to=2-2]
    \end{tikzcd}\]
    it assigns the pair consisting of the class $[f] : \pi_{n-1}(S^{n-1}) \simeq \ZZ$ and the class of the map $S^{n} \equiv \Sigma S^{n-1} \to S^n$ induced from $\alpha$ by the universal property of the suspension.
\end{lem}
\begin{proof}
    Since the map $S^{n-1} \to S^n$ is constant, it induces the trivial map on homotopy groups, so the long exact sequence of the pair (Proposition \ref{prp:LESproof})
    $$\pi_n(S^{n-1}) \xrightarrow{0} \pi_n(S^n) \xrightarrow{j_n} \pi_n(S^n, S^{n-1}) \xrightarrow{\kappa_n} \pi_{n-1}(S^{n-1}) \xrightarrow{0} \pi_{n-1}(S^n)$$
    yields a short exact sequence of groups
    $$0 \to \pi_n(S^n) \xrightarrow{j_n} \pi_n(S^n, S^{n-1}) \xrightarrow{\kappa_n} \pi_{n-1}(S^{n-1}) \to 0.$$
    
    It is easy to see that this splits via a section of $\kappa_n$ that is the truncation of the map that completes a map $f: S^{n-1} \to S^{n-1}$ into a square with a trivial $\alpha$. Moreover, by definition $j_n$ factors as $\pi_n(S^n) \to \pi_n(S^n, \sgtn) \to \pi_n(S^n, S^{n-1})$, where the first map is exactly the one using the universal property of the suspension. It follows that the splitting has the form described by the statement of the lemma.
\end{proof}

\begin{proof}[Proof of the Leibniz formula (Proposition \ref{4.1})]
    As always, we can assume dealing with elements $|x|_0, |y|_0$ for diagrams $x, y$.

    By naturality (Proposition \ref{circ-prod-nat}), it is enough to check the case $(B, A) \equiv (\sgtn, S^{n}), (D,C) \equiv (\sgtn, S^{m}),$ with $x, y$ corresponding to identity maps of arrows. Then $(P, A \wedge C) \simeq (\Sigma S^{n,m}, S^{n, m})$. We will compare the images of the right- and left-hand sides of the formula in $\pi_{n+m+1}(S^{n+m+1}, S^{n+m}) \simeq \ZZ \times \ZZ$ under the inverse of $(\Sigma j^{n, m}, j^{n, m})$. For $x$ the generator of $\pi_{n+1}(\sgtn, S^{n}),$ $\kappa x$ is the generator of $\pi_{n}(S^{n}),$ which can be represented by the following diagram:

    \[\begin{tikzcd}
    	{S^{n-1}} & \sgtn \\
    	\sgtn & {S^{n} \equiv \Sigma S^{n-1}.}
    	\arrow[from=1-1, to=1-2]
    	\arrow[from=1-1, to=2-1]
    	\arrow["\mmerid"', between={0.2}{0.8}, Rightarrow, from=1-2, to=2-1]
    	\arrow["N", from=1-2, to=2-2]
    	\arrow["N"', from=2-1, to=2-2]
    \end{tikzcd}\]
    Then for the first term of the right-hand side we consider the image of $\kappa x \circ y$ in $\pi_{n+m+1}(\sgtn, S^{n+m})$. For that we analyze the diagram

    \[\begin{tikzcd}
    	{S^{n+m}} &[1.8em] {\sgtn \coprod^{S^{n-1,m}} \sgtn} & {\sgtn \wedge \sgtn  \coprod^{\sgtn \wedge S^{m}} S^{n,m}} &[4.8em] {S^{n+m}} \\
    	\sgtn & {\sgtn \wedge \sgtn} & {S^{n} \wedge \sgtn} & {\sgtn.}
    	\arrow["{\Sigma j^{n-1,m}}", from=1-1, to=1-2]
    	\arrow[from=1-1, to=2-1]
    	\arrow[from=1-2, to=1-3]
    	\arrow[from=1-2, to=2-2]
    	\arrow["{(\bs{}, (j^{n,m})^{-1})}", from=1-3, to=1-4]
    	\arrow[from=1-3, to=2-3]
    	\arrow[from=1-4, to=2-4]
    	\arrow[from=2-1, to=2-2]
    	\arrow[from=2-2, to=2-3]
    	\arrow[from=2-3, to=2-4]
    \end{tikzcd}\]
    Let us compute the top morphism. The middle map is from the cube
\[\begin{tikzcd}
	{S^{n-1,m}} && \sgtn \\
	& \sgtn && \sgtn \\
	\sgtn && {\sgtn \coprod^{S^{n-1,m}} \sgtn} \\
	& {S^{n,m}} && {\sgtn \coprod^{\sgtn \wedge S^{m}} S^{n,m}.}
	\arrow[""{name=0, anchor=center, inner sep=0}, line width=1pt, -{cm bold to}, from=1-1, to=1-3]
	\arrow[from=1-1, to=2-2]
	\arrow[line width=1pt, -{cm bold to}, from=1-1, to=3-1]
	\arrow["\refl"{description}, shorten <=9pt, shorten >=9pt, Rightarrow, from=1-3, to=2-2]
	\arrow[from=1-3, to=2-4]
	\arrow[from=1-3, to=3-3]
	\arrow[shift right=4, line width=1pt, -{cm bold to}, curve={height=-12pt}, shorten <=9pt, shorten >=9pt, from=1-3, to=4-4]
	\arrow[""{name=1, anchor=center, inner sep=0}, from=2-2, to=2-4]
	\arrow["\alpha"', shorten <=5pt, shorten >=5pt, Rightarrow, from=2-2, to=3-1]
	\arrow[from=2-2, to=4-2]
	\arrow["\refl", shorten <=40pt, shorten >=40pt, Rightarrow, from=2-4, to=4-2]
	\arrow[from=2-4, to=4-4]
	\arrow[from=3-1, to=3-3]
	\arrow[from=3-1, to=4-2]
	\arrow[line width=1pt, -{cm bold to}, curve={height=6pt}, from=3-1, to=4-4]
	\arrow[dashed, from=3-3, to=4-4]
	\arrow[from=4-2, to=4-4]
	\arrow["\lrcorner"{anchor=center, pos=0.125, rotate=180}, draw=none, from=3-3, to=0]
	\arrow["\lrcorner"{anchor=center, pos=0.125, rotate=180}, draw=none, from=4-4, to=1]
\end{tikzcd}\]
    A map from the pushout of $\sgtn \xleftarrow{} S^{n-1,m} \xrightarrow{} \sgtn$ is specified by
    a witness of commutativity of the square drawn with bold arrows in the diagram. This square commutes by the composition of the three 2-cells shown in the diagram: the two labelled $\refl$ are trivial, and the third one, $\alpha$, is $\merid_{S^{n-1}}$ smashed with $S^{m}$. By Lemma \ref{le:merid-smash-Sm}
    we conclude that our map $S^{n+m} \to S^{n+m}$ is pointed homotopic to the identity.

    Similarly, for the second term of the right-hand side we consider the diagram
    \[\begin{tikzcd}
    	{S^{n+m}} &[1.8em] {\sgtn \coprod^{S^{n,m-1}} \sgtn} & {S^{n,m} \coprod^{S^{n} \wedge \sgtn} \sgtn \wedge \sgtn} & {S^{n+m}} \\
    	\sgtn & {\sgtn \wedge \sgtn} & {\sgtn \wedge S^{m}} & {\sgtn.}
    	\arrow["{\Sigma j^{n,m-1}}", from=1-1, to=1-2]
    	\arrow[from=1-1, to=2-1]
    	\arrow[from=1-2, to=1-3]
    	\arrow[from=1-2, to=2-2]
    	\arrow[from=1-3, to=1-4]
    	\arrow[from=1-3, to=2-3]
    	\arrow[from=1-4, to=2-4]
    	\arrow[from=2-1, to=2-2]
    	\arrow[from=2-2, to=2-3]
    	\arrow[from=2-3, to=2-4]
    \end{tikzcd}\]
    This time the nontrivial cell is the top of the cube, which is \textit{minus} the meridian smashed with $S^{n}$ on the left. This produces the map $(-1)^{n+1}$, per Lemma \ref{le:Sn-smash-merid}. (Intuitively, the meridian goes along the first coordinate, while $S^{n} \wedge \merid$ goes along the $(n+1)$-st coordinate, so $n$ coordinate transpositions are required to connect the two.)

    We are interested in the images of these under the map
    $(\sgtn, S^{n+m}) \to (S^{n+m+1}, S^{n+m})$.

The first summand, $(\Sigma j^{n, m}, j^{n, m})^{-1}_*l_* ((\kappa x) \circ y)$, is represented by
\[\begin{tikzcd}[column sep=1.4em]
	{S^{n+m}} &[2.6em] {\sgtn \coprod^{S^{n-1,m}} \sgtn} & {\sgtn \wedge \sgtn  \coprod^{\sgtn \wedge S^{m}} S^{n,m}} &[-0.7em] &[4.2em] {S^{n+m}} \\
	\sgtn & {\sgtn \wedge \sgtn} & {S^{n} \wedge \sgtn} & {\sgtn \coprod^{S^{n,m}} \sgtn} & {S^{n+m+1}.}
	\arrow["{\Sigma j^{n-1,m}}", from=1-1, to=1-2]
	\arrow[from=1-1, to=2-1]
	\arrow[from=1-2, to=1-3]
	\arrow[from=1-2, to=2-2]
	\arrow[from=1-3, to=1-5]
	\arrow[from=1-3, to=2-3]
	\arrow[curve={height=-12pt}, from=1-3, to=2-4]
	\arrow[from=1-5, to=2-5]
	\arrow[from=2-1, to=2-2]
	\arrow[from=2-2, to=2-3]
	\arrow[from=2-3, to=2-4]
	\arrow["{(\Sigma j^{n, m})^{-1}}", from=2-4, to=2-5]
\end{tikzcd}\]
The witness of the triangle's commutativity is $\refl,$ so the resulting diagram has the standard generator on top and a trivial 2-cell.

For $r_* (x \circ \kappa y )$, we get
\[\begin{tikzcd}[column sep=1.6em]
	{S^{n+m}} &[2.6em] {\sgtn \coprod^{S^{n,m-1}} \sgtn} & {S^{n,m} \coprod^{S^{n} \wedge \sgtn}\sgtn \wedge \sgtn} && {S^{n+m}} \\
	\sgtn & {\sgtn \wedge \sgtn} & {\sgtn \wedge S^{m}} & {\sgtn \coprod^{S^{n,m}} \sgtn} & {S^{n+m+1}.}
	\arrow["{\Sigma j^{n,m-1}}", from=1-1, to=1-2]
	\arrow[from=1-1, to=2-1]
	\arrow[from=1-2, to=1-3]
	\arrow[from=1-2, to=2-2]
	\arrow[from=1-3, to=1-5]
	\arrow[from=1-3, to=2-3]
	\arrow[from=1-3, to=2-4]
	\arrow[from=1-5, to=2-5]
	\arrow[from=2-1, to=2-2]
	\arrow[from=2-2, to=2-3]
	\arrow["{\inr\circ !}", from=2-3, to=2-4]
	\arrow[from=2-4, to=2-5]
\end{tikzcd}\]
This time the commutativity of the triangle is witnessed by $\inl(s) \mapsto \glue(s)$. Denote the top middle map by $t$. Then we know from above that $(j^{n,m})^{-1} \circ t \circ \Sigma j^{n, m-1} =_* (-1)^{n+1}$. Then our 2-cell computes as
\begin{align*}
(\Sigma j^{n,m})^{-1} \whisk \glue \whisk (t \circ \Sigma j^{n, m-1})
&=_* (\Sigma j^{n,m})^{-1} \whisk \glue \whisk j^{n,m} \whisk ((j^{n,m})^{-1} \circ t \circ \Sigma j^{n, m-1}) \\
&=_* (-1)^{n+1} (\Sigma j^{n,m})^{-1} \whisk \glue \whisk  j^{n,m} \\
&=_* (-1)^{n+1} \glue \equiv (-1)^{n+1} \merid.
\end{align*}
Hence, under the isomorphism $\pi_{n+m+1}(S^{n+m+1}, S^{n+m}) \simeq \mathbb{Z} \times \mathbb{Z}$ from Lemma \ref{le:pi(S,S)}, we get
$(-1)^{n+1}(1, 1)$.

Finally, for $j_* \kappa (x \circ y)$ we get
\[\begin{tikzcd}
	{S^{n+m}} & \sgtn &[1.8em] \sgtn & \sgtn & {S^{n+m}} \\
	\sgtn & {S^{n+m+1}} & {\Sigma S^{n,m}} & {\Sigma S^{n,m}} & {S^{n+m+1}.}
	\arrow[from=1-1, to=1-2]
	\arrow[from=1-1, to=2-1]
	\arrow[from=1-2, to=1-3]
	\arrow["\mmerid"', shorten <=5pt, shorten >=8pt, Rightarrow, from=1-2, to=2-1]
	\arrow[from=1-2, to=2-2]
	\arrow[from=1-3, to=1-4]
	\arrow[from=1-3, to=2-3]
	\arrow[from=1-4, to=1-5]
	\arrow[from=1-4, to=2-4]
	\arrow[from=1-5, to=2-5]
	\arrow[from=2-1, to=2-2]
	\arrow["{\Sigma j^{n,m}}", from=2-2, to=2-3]
	\arrow["\id"', from=2-3, to=2-4]
	\arrow[from=2-4, to=2-5]
\end{tikzcd}\]
(Here the middle two arrows at the bottom form the top arrow of the product of generators; we then post-whisker the meridian with it to get an element of a relative homotopy group, and regard the image of this in $\pi_{m+n+1}(S^{m+n+1}, S^{m+n})$.)

Our proposition thus becomes $(0, 1) = -(1, 0) + (-1)^{n+1} \cdot ((-1)^{n+1}, (-1)^{n+1})$.
\end{proof}

\begin{re}
    When $n=0$ or $m=0$, the proposition's statement does not make sense: addition is not defined in relative $\pi_1$, and we'd also have to pass through the undefined $\pi_0(-, \sgtn)$. Moreover, relative $\pi_2$ is not generally abelian. The plus on the right-hand side of the Leibniz formula is still justified since, as the proof shows, the two terms commute.
\end{re}

\begin{re}
    The universal case $(B, A) = (\sgtn, S^n), (D, C) = (\sgtn, S^{m})$ with $x,y$ identity diagrams can also be described in the language of the homotopy groups of fibers. (We omit the derivation as it is straightforward but technical. We won't use description further.) The element $\kappa x \circ y$ is then the image of $(\adj{n}\id, \adj{m}\id)$ under the composition
    \begin{align*}
        \Omega^{n}S^{n} \times \Omega^{m}S^{m} \xrightarrow{I_{n-1,m}^{\Omega S^n, S^m}} \Omega^{n+m}(\Omega S^{n} \join S^{m}) \to \\
        \Omega^{n+m} \Sigma (\Omega S^{n} \wedge S^{m})
        \to \Omega^{n+m} \Sigma \Omega (S^{n} \wedge S^{m}) \to \\
        \Omega^{n+m} S^{n, m}.
    \end{align*}
    This then gets mapped into $\Omega^{n+m} S^{n, m} \times \Omega^{n+m+1} \Sigma S^{n, m} = \Omega^{n+m} \fib (S^{n,m} \to \Sigma S^{n, m})$ by the function $\Omega^{n+m}\lambda x.(x, 0)$.
    \\
    A similar composite computes $x \circ \kappa y$, which then gets sent to an element of
    $\Omega^{n+m} S^{n, m} \times \Omega^{n+m+1} \Sigma S^{n, m}$ by the map $\Omega^{n+m}\lambda x.(x, \merid_x)$ (since this time the map of pairs from which it is induced commutes by $\merid$).

    For example, the intuitive reason why the first coordinates of $\kappa x \circ y$ and $x \circ \kappa y$ differ by the sign $(-1)^{n-1}$ is clear if one keeps track of the order of the functors $\Omega$ and $\Sigma$. The former applies to the canonical element the composition
    \begin{align*}
    &(\onum{n} \ldots \onum{1} \snum{1} \ldots \snum{n} S^0) \times (\onum{n+m} \ldots \onum{n+1}  \snum{n+1} \ldots \snum{n+m} S^0) \to \\
    &\onum{0}\onum{n+m} \ldots \onum{2} ((\onum{1} \snum{1}\ldots \snum{n} S^0) \join \snum{n+1} \ldots \snum{n+m} S^0) \to \\
    &\onum{0}\onum{n+m} \ldots \onum{2} \snum{0} \onum{1} ((\snum{1}\ldots \snum{n} S^0) \wedge \snum{n+1} \ldots \snum{n+m} S^0) \to \\
    &\onum{0}\onum{n+m} \ldots \onum{2}  ((\snum{1}\ldots \snum{n} S^0) \wedge \snum{n+1} \ldots \snum{n+m} S^0) \to \\
    & \onum{0}\onum{n+m} \ldots \onum{2}  \snum{1} \ldots \snum{n+m} S^0,
    \end{align*}
    where we label each $\Sigma$ and $\Omega$ to keep track of their order.
    Thus the first $\Omega$ annihilates with the join coordinate, and the order of the others is preserved.
    For $x \circ \kappa y$, we take the map
    $$(\onum{n} \ldots \onum{1} \snum{1} \ldots \snum{n} S^0) \times (\onum{n+m} \ldots \onum{n+1}  \snum{n+1} \ldots \snum{n+m} S^0) \to \onum{0}\onum{n+m} \ldots \onum{n+2} \onum{n} \ldots \onum{1}  \snum{1}\ldots \snum{n+1} \ldots \snum{n+m} S^0.$$
    To connect the two, we can push $\onum{n+1}$ to the front of $\snum{0}$ (which requires $n-1$ transpositions) and then use the triangle identity to appeal to the diagram of natural transformations
    \[\begin{tikzcd}
    	{\onum{n+1} \snum{0}  \onum{1}} & {\onum{n+1}} & {\onum{1} \snum{0}  \onum{n+1} } \\
    	&& {\onum{1}.}
    	\arrow[from=1-1, to=1-2]
    	\arrow["{\eta_{\onum{n+1}}}", from=1-2, to=1-3]
    	\arrow["\simeq", from=1-2, to=2-3]
    	\arrow["{\onum{1}\epsilon}", from=1-3, to=2-3]
    \end{tikzcd}\]
    This approach may be easier to formalize.
\end{re}

\subsection{On towers and spectral sequences}\label{pairing_on_towers}

\subsubsection{Setup}

\begin{de}
    A \textit{tower} $X_\bullet$ is a $\mathbb{Z}$-indexed sequence of pointed maps
    $$\ldots \xrightarrow{t_{q-2}} X_q \xrightarrow{t_{q-1}} X_{q-1} \xrightarrow{t_q} \ldots,$$
    called the \textit{structure maps}.
\end{de}

For simplicity, assume that the structure maps are pointed by $\refl$. This does not lose generality if only a finite number of structure maps is nontrivial, which will be the case in our applications.

Recall from \cite{Flo} that a tower $W_\bullet$ gives rise to a spectral sequence $E^{\bullet,\bullet}_\bullet(W_\bullet)$, where $E_1^{p,q} = \pi_{p+q}(\fib(W_{q} \to W_{q-1})) = \pi_{p+q+1}(W_{q-1}, W_{q})$. Under boundedness assumptions, it converges to $\pi_{\bullet}(\lim W_{\bullet})$.

\begin{exa}
    Given a truncated type family $B: A \to U$, setting $W_q :\equiv \depfun{a:A}{||B(a)||_q}$\footnote{Note that $W_q$ is not, generally, the $q$-th truncation of $\depfun{a:A}{B(a)}$. If this is the case, the spectral sequence turns out to be degenerate.} yields the \textit{Federer spectral sequence} (\cite{Fed}).
    $$H^{-p}(A, \lambda a.\pi_q(B(a))) \Rightarrow{} \pi_{p+q}(\depfun{a:A}{B(a)}).$$
    Given an action of a group $G$ on a space $X$, this specializes to
    $$H^{-p}(G, \pi_q(X)) \Rightarrow{} \pi_{p+q}(X^G).$$
    Here $\pi_q(X)$ is a $G$-module in the obvious way and $X^G$ denotes the homotopy fixed points of $X$. This is classically known as the \textit{homotopy fixed points spectral sequence}.

\end{exa}

We define $E_1^{p, q}(W_\bullet) : \equiv \pi_{p+q}(W_q, W_{q+1})$ when $p+q \geq 1$ and $0$ otherwise. Note that this is a $(2, -1)$-shift of \cite{Flo}'s pages. This indexing is the obvious choice for the pairing to make sense, but it is also consistent with the convergence statement: see Section \ref{conv}.

\begin{de}\label{de:good}
    Call a tower \textit{good} if in the associated spectral sequence all terms $E_1^{\bullet, \bullet}$ are abelian groups.
\end{de}
This is the case, for example, when each $\fib(W_{q} \to W_{q-1})$ is an H-space, as will happen in the case of the Serre spectral sequence. To avoid dealing with corner cases, from now on \textit{we consider only good towers}.

We give a geometric description of the elements of all pages of the spectral sequence, which allows for explicit manipulations with geometric notions such as the differential and the multiplication.

Define, for any $r \geq 1$, $(\ZZ\times\ZZ)$-graded groups $Z_r, B_r$ by
$$Z^{p,q}_r : \equiv \im (\iota_r: \pi_{p+q}(W_q, W_{q+r}) \to \pi_{p+q}(W_q, W_{q+1})),$$
$$B^{p,q}_r : \equiv \im (\epsilon_r: \pi_{p+q+1}(W_{q-r+1}, W_q) \to \pi_{p+q}(W_q) \to \pi_{p+q}(W_q, W_{q+1})).$$ Clearly, $B_r^{p,q} \subset Z_r^{p,q}$ as elements in the image of $\pi_{p+q}(W_q, \sgtn) \to \pi_{p+q}(W_q, W_{q+1})$ are \textit{a fortiori} in the image of $\pi_{p+q}(W_q, W_{q+r}) \to \pi_{p+q}(W_q, W_{q+1})$. Thus, we get a chain of groups $0 = B^{p,q}_1 \subset \ldots \subset B^{p,q}_r \subset B^{p,q}_{r+1} \subset \ldots \subset Z^{p,q}_{r+1} \subset Z^{p,q}_{r} \subset \ldots \subset Z^{p,q}_1 = \pi_{p+q
}(W_q, W_{q+1}).$

Let $E_r^{p,q} : \equiv Z_r^{p,q} / B_r^{p,q}.$

We define $d_{r}: E^{p,q}_{r} \to E^{p-r-1,q+r}_{r}.$ For an element of the form $y = [\iota_r(x)]$, we'd like its image to be represented by the image of $x$ under $\pi_{p+q}({W_{q}, W_{q+r}})\xrightarrow{\kappa} \pi_{p+q-1}({W_{q+r}}) \xrightarrow{\rho} \pi_{p+q-1}({W_{q+r}, W_{q+r+1}})$. Clearly, this element is in $Z^{\bullet,\bullet}_r$ because it is in the image of $\pi_{p+q-1}(W_{q+r}, \sgtn) \to \pi_{p+q-1}(W_{q+r}, W_{q+2r}) \to \pi_{p+q-1}(W_{q+r}, W_{q+r+1}).$
Next, we need to check that if $\iota_r (x) = 0 $ then $[\rho \kappa x] = 0$, that is, $\rho \kappa (x) \in B_r^{{p-r-1, q+r}}.$ This is seen from the diagram below, in which the middle row is exact and the square commutes:

\begin{equation}\label{eq:d-well-defined}
\begin{tikzcd}
	& {\pi_{p+q}(W_q, W_{q+r})} & {\pi_{p+q}(W_q, W_{q+1})} \\
	{\pi_{p+q}(W_{q+1}, W_{q+r})} & {\pi_{p+q-1}(W_{q+r})} & {\pi_{p+q-1}(W_{q+1})} \\
	& {\pi_{p+q-1}(W_{q+r}, W_{q+r+1}).}
	\arrow[from=1-2, to=1-3]
	\arrow[from=1-2, to=2-2]
	\arrow[from=1-3, to=2-3]
	\arrow[from=2-2, to=2-3]
	\arrow[from=2-1, to=2-2]
	\arrow[from=2-2, to=3-2]
\end{tikzcd}
\end{equation}

Finally, if $\iota_r(x) : B_r^{p,q}$ then $\rho \kappa x : B_r^{p-r-1,q+r};$ indeed, it is actually $0$.
Thus, we can set $d_r([\iota_r(x)]) :\equiv [\rho \kappa x].$ Note that $\im d_r = B_{r+1}/B_{r}$ and $\ker d_r = Z_{r+1}/B_r$, so $E_{r+1} = \ker d_r / \im d_r.$

Recall that in \cite{Flo} the passage from a tower to a spectral sequence is mediated by an \textit{exact couple}. That is, one first forms graded groups $\bar{E}_1 : \equiv \bigoplus_{n,s}\pi_n(\fib(W_s \to W_{s-1}))$ and $\bar{D}_1 : \equiv \bigoplus_{n,s}\pi_n(W_s)$; the structure maps then determine a natural map $\bar{i}_1: \bar{D}_1 \to \bar{D}_1$, fiber inclusions give rise to a map $\bar{k}_1: \bar{E}_1 \to \bar{D}_1$, and the differentials from the long exact sequences of pairs $(W_{s-1}, W_{s})$ give a map $\bar{j}_1: \bar{D}_1 \to \bar{E}_1$. These data determine a so-called exact couple, from which one iteratively constructs \textit{derived} exact couples $(\bar{E}_r, \bar{D}_r, \bar{i}_r, \bar{j}_r, \bar{k}_r)$. Concretely, $\bar{D}_r^{p+q, q}$ is defined as $\im\bar{i}_{r-1}^{p+q,q} = \im (\iota^{abs}_r: \pi_{p+q}(W_{q+r-1}) \to \pi_{p+q}(W_{q}))$, and $\bar{E}_r$ is defined as the homology of $\bar{E}_{r-1}$ with respect to the differential $\bar{d}_{r-1} : \equiv \bar{j}_{r-1} \circ \bar{k}_{r-1}$. See \cite{Flo}, 5.2--5.3 for details.
\begin{prp}\label{prp:r-cycles}
$d_r : E_r \to E_r$ coincides with the differential from the exact couple construction.
\end{prp}

\begin{proof}
    Denote the exact couple graded groups and maps as above.
    Set $k_r: E_r \to \bar{D}_r$ by $k_r[\iota_r x] :\equiv \kappa \iota_r x,$ which is clearly well-defined ($\kappa \iota_r x \in \im \iota_r^{abs}$ and elements of $B_r$ are sent to $0$ by $\kappa$).
    Define $j_r: \bar{D}_r \to E_r$ by $j_r (\iota^{abs}_r(x)) :\equiv [\rho (x)].$ The lower two rows of diagram \eqref{eq:d-well-defined} show that this is well-defined: if $\iota^{abs}_r(x) = 0$ for some $x : \pi_{p+q-1}(W_{q+r})$, then, by the exactness of the middle row, $x = \kappa z$ for some $z : \pi_{p+q}(W_{q+1}, W_{q+r})$, and hence $\rho(x) = \rho \kappa z$ lies in $B_r$. Clearly, $d_r = j_r \circ k_r.$

    For $r=1$, it is clear that $(E_r, j_r, k_r)$ equals $(\bar{E}_r, \bar{j}_r, \bar{k}_r)$. Thus, we only need to show inductively that $(j_r, k_r) = (\bar{j}_r, \bar{k}_r)$.

    Recall that the derived couple $(\bar{E}_{r+1}, \bar{D}_{r+1}, \bar{i}_{r+1}, \bar{j}_{r+1}, \bar{k}_{r+1})$ is defined as follows: $\bar{D}_{r+1} :\equiv \im \bar{i}_r$, $\bar{E}_{r+1} :\equiv \ker \bar{d}_r / \im \bar{d}_r$, the map $\bar{i}_{r+1}$ is the restriction of $\bar{i}_r$, and
    $$\bar{k}_{r+1}([e]) :\equiv \bar{k}_r(e) \text{ for } e : \ker \bar{d}_r, \qquad \bar{j}_{r+1}(\bar{i}_r(a)) :\equiv [\bar{j}_r(a)] \text{ for } a : \bar{D}_r.$$
    Now assume that $(E_r, j_r, k_r) = (\bar{E}_r, \bar{j}_r, \bar{k}_r)$. Then $\bar{d}_r = \bar{j}_r \circ \bar{k}_r = j_r \circ k_r = d_r$, so $\bar{E}_{r+1} = \ker d_r / \im d_r = (Z_{r+1}/B_r)/(B_{r+1}/B_r)$, which we identify with $E_{r+1} \equiv Z_{r+1}/B_{r+1}$. Under this identification, the class $[\iota_{r+1} x] : E_{r+1}$ of an element $x : \pi_{p+q}(W_q, W_{q+r+1})$ corresponds to the class $[[\iota_{r+1} x]]$ of $[\iota_{r+1} x] : \ker d_r$. For such $x$ we have
    $$\bar{k}_{r+1} [[\iota_{r+1}  x]] \equiv \bar{k}_r ([\iota_{r+1}  x]) = k_r ([\iota_{r+1} x]) \equiv \kappa \iota_{r+1} x \equiv k_{r+1}[\iota_{r+1} x].$$
    Similarly, every element of $\bar{D}_{r+1}$ is of the form $\iota^{abs}_{r+1} x = \bar{i}_r (\iota^{abs}_{r} x)$, and
    \[\bar{j}_{r+1} (\iota^{abs}_{r+1} x) \equiv \bar{j}_{r+1}  (\bar{i}_r (\iota^{abs}_{r} x)) \equiv [\bar{j}_{r} (\iota^{abs}_{r} x)] = [{j}_{r} (\iota_{r}^{abs} x)] \equiv [[\rho (x)]] = {j}_{r+1} (\iota^{abs}_{r+1} x). \qedhere\]
\end{proof}

\subsubsection{The main theorem}

\begin{de}\label{coherent-pairing}
Suppose given towers $X_\bullet, W_\bullet, Y_\bullet$. A \textit{coherent pairing} is a function of type $\depfun{n,m: \mathbb{Z}}{(W_n \wedge X_m \pto Y_{n+m})}$ such that all squares (denoted $cohL,$ $cohR$)
    \[\begin{tikzcd}
    	{W_n \wedge X_m} & {Y_{n+m}} & {Y_{n+m}} & {W_n \wedge X_{m}} \\
    	{W_{n-1} \wedge X_m} & {Y_{n+m-1}} & {Y_{n+m-1}} & {W_n \wedge X_{m-1}}
    	\arrow[from=1-1, to=1-2]
    	\arrow[from=1-1, to=2-1]
    	\arrow[from=1-2, to=2-2]
    	\arrow[from=1-3, to=2-3]
    	\arrow[from=1-4, to=1-3]
    	\arrow[from=1-4, to=2-4]
    	\arrow[from=2-1, to=2-2]
    	\arrow[from=2-4, to=2-3]
    \end{tikzcd}\]
    commute in $U_*$, and all cubes
\[\begin{tikzcd}
	& {W_q \wedge X_t} \\
	{W_{q-1} \wedge X_t} & {Y_{q+t}} & {W_{q} \wedge X_{t-1}} \\
	{Y_{q+t - 1}} & {W_{q-1} \wedge X_{t-1}} & {Y_{q+t - 1}} \\
	& {Y_{q+t - 2}}
	\arrow[from=2-1, to=3-2]
	\arrow[from=1-2, to=2-3]
	\arrow[from=2-3, to=3-2]
	\arrow[from=1-2, to=2-2]
	\arrow[from=2-2, to=3-1]
	\arrow[from=2-2, to=3-3]
	\arrow[from=3-3, to=4-2]
	\arrow[from=3-1, to=4-2]
	\arrow[from=3-2, to=4-2]
	\arrow[from=2-3, to=3-3]
	\arrow[from=1-2, to=2-1]
	\arrow[from=2-1, to=3-1]
\end{tikzcd}\]
    have fillers (in $U$).
\end{de}

\begin{thm}\label{th:product-SS}
    Suppose given towers $X_\bullet, W_\bullet, Y_\bullet$ and a coherent pairing $t: \depfun{n,m: \mathbb{Z}}{(W_n \wedge X_m \pto Y_{n+m})}$.
    Then there is an induced pairing of the spectral sequences $E_\bullet^{p,q}(W) \times E_\bullet^{s,t}(X) \to E_\bullet^{p+s,q+t}(Y)$ that enjoys distributivity and the Leibniz rule: $ -(d_\bullet^{p,q} x) \cdot y + (-1)^{p+q} x \cdot d_\bullet^{s,t}y = d_\bullet^{p+s, q+t} (x \cdot y)$.
    Moreover, if $E^{p,q}_{r+1}$ is seen as a subquotient of $E^{p,q}_r$, then we have $[x] \cdot [y] = [x \cdot y]$ for $x, y : Z^{\bullet,\bullet}_r$.
\end{thm}

\begin{proof}
We want to define the pairing $E_r^{p,q}(W) \otimes E_r^{s,t}(X) \to E_r^{p+s,q+t}(Y) = Z_r^{p+s,q+t}/B_r^{p+s,q+t}$ by its action on representatives:
$$\pi_n({W_{q}, W_{q+1}}) \otimes \pi_m(X_{t}, X_{t+1}) \xrightarrow{\circ} \pi_{n+m}\left(W_{q} \wedge X_{t},  (W_{q+1} \wedge X_{t}) \coprod^{W_{q+1} \wedge X_{t+1}} (W_{q} \wedge X_{t+1})\right) \to \pi_{n+m}(Y_{q+t}, Y_{q+t+1}).$$
Here $n = p+q, m = s+t.$ Denote this pairing by $\cdot$, so that $x \cdot y = \rho_*(x \circ y),$ where $\rho$ is the obvious pointed map $(W_{q} \wedge X_{t},  W_{q+1} \wedge X_{t} \coprod^{W_{q+1} \wedge X_{t+1}} W_{q} \wedge X_{t+1}) \to (Y_{q+t}, Y_{q+t+1})$ (we need the second coherence above (the cube fillers) for $\rho$ to exist).
We need to check several things:

1) We need to prove
$\iota_r(x) = 0 \to \iota_r(x) \cdot y = 0$ (and $\iota_r(y) = 0 \to x \cdot \iota_r(y) = 0$), which is clear.

2) We need to show that if $x,y : Z_r $ then $x \cdot y \in Z_r$. By assumption, the boundaries of the spheroids $x,y$ lift $r-1$ steps ``up'', giving spheroids $\bar{x} : \pi_n(W_q, W_{q+r})$ and $\bar{y}: \pi_m(X_t, X_{t+r})$. We had no access to $\bar{x}, \bar{y}$ when constructing the element $x \cdot y,$ but now we are proving a proposition, so we can assume given explicit lifts $\bar{x}, \bar{y}$.
Then we'd like $x \cdot y \in \im \iota_r$ to follow from the diagram

\[\begin{tikzcd}
	{\pi_n(W_q, W_{q+r}) \times \pi_m(X_t, X_{t+r})} & {\pi_{n+m}(W_q \wedge X_t, (W_{q+r} \wedge X_t) \coprod^{\ldots} (W_q \wedge X_{t+r}))} & {\pi_{n+m}(Y_{q+t}, Y_{q+t+r})} \\
	{\pi_n(W_q, W_{q+1}) \times \pi_m(X_t, X_{t+1})} & {\pi_{n+m}(W_q \wedge X_t, (W_{q+1} \wedge X_t) \coprod^{\ldots} (W_q \wedge X_{t+1}))} & {\pi_{n+m}(Y_{q+t}, Y_{q+t+1}).}
	\arrow["\circ", from=1-1, to=1-2]
	\arrow[from=1-1, to=2-1]
	\arrow["\blacksquare", from=1-2, to=1-3]
	\arrow[from=1-2, to=2-2]
	\arrow[from=1-3, to=2-3]
	\arrow["\circ", from=2-1, to=2-2]
	\arrow["{\rho_*}", from=2-2, to=2-3]
\end{tikzcd}\]
Here the down-right-right composite takes $(\bar{x},\bar{y})$ to $x \cdot y$. The rightmost map is $\iota_r$, so to prove $x \cdot y \in \im \iota_r$ it's enough to show that such a commutative diagram exists.
The left square here is the naturality of the product $\circ$.
The map labeled by $\blacksquare$ is defined so as to make the right-hand square commute; namely, it can be defined as induced from the diagram below, which we'd like to define as the composite of the two squares:
\[\begin{tikzcd}
	{W_{q+r} \wedge X_t \coprod^{W_{q+r} \wedge X_{t+r}} W_q \wedge X_{t+r}} & {Y_{q+t+r}} \\
	{W_{q+1} \wedge X_t \coprod^{W_{q+1} \wedge X_{t+1}} W_{q} \wedge X_{t+1}} & {Y_{q+t+1}} \\
	{W_q \wedge X_t} & {Y_{q+t}.}
	\arrow[from=2-2, to=3-2]
	\arrow[from=1-2, to=2-2]
	\arrow[from=1-1, to=1-2]
	\arrow[from=2-1, to=2-2]
	\arrow[from=3-1, to=3-2]
	\arrow[from=2-1, to=3-1]
	\arrow[from=1-1, to=2-1]
\end{tikzcd}\]
For that we need to prove that the top square actually commutes.

Its commutativity amounts to a filler of the following cube:
    \[\begin{tikzcd}[row sep=3.8em, column sep=1.6em]
    	& {W_{q+r} \wedge X_t} \\
    	{W_{q+r} \wedge X_{t+r}} & {Y_{q+t+2r}} &[1.6em] {Y_{q+t+r}} && {W_{q+1} \wedge X_t} \\
    	& {W_{q} \wedge X_{t+r}} && {W_{q+1} \wedge X_{t+1}} & {Y_{q+t+2}} &[1.6em] {Y_{q+t+1}} \\
    	&&&& {W_{q} \wedge X_{t+1}.}
    	\arrow[from=2-1, to=1-2]
    	\arrow[from=1-2, to=2-3]
    	\arrow[from=2-1, to=3-2]
    	\arrow[from=3-2, to=2-3]
    	\arrow[dashed, from=2-1, to=2-2]
    	\arrow[dashed, from=2-2, to=2-3]
    	\arrow[from=3-4, to=2-5]
    	\arrow[from=2-5, to=3-6]
    	\arrow[from=3-4, to=4-5]
    	\arrow[from=4-5, to=3-6]
    	\arrow[dashed, from=3-4, to=3-5]
    	\arrow[dashed, from=3-5, to=3-6]
    	\arrow[color=black!50, from=1-2, to=2-5]
    	\arrow[color=black!50, from=2-1, to=3-4]
    	\arrow[color=black!50, from=2-2, to=3-5]
    	\arrow[color=black!50, from=2-3, to=3-6]
    	\arrow[color=black!50, from=3-2, to=4-5]
    \end{tikzcd}\]
Note that the commutativity of the left and the right (diamond-shaped) sides of the cube is exhibited by connecting their pairs of edges with the dashed compositions via 2-dimensional coherences of the pairing. But then the cube is seen to be filled by the fillers of its upper and lower halves, which are filled by the 3-dimensional coherences of the pairing.

Having established the commutativity of the top square, we obtain the desired commutativity of the right square of the diagram of homotopy groups.

3) To factor through $Z_r \to E_r,$ it is enough to show the distributivity and that multiplying by a boundary (an element of $B_r$) gives a boundary. The former follows from the distributivity of $\circ$ (Corollary \ref{cor:distr}), and the latter immediately follows from the Leibniz formula, which we prove next.

4) We prove the Leibniz formula.

We need to show that $d_r^{p+s, q+t} (x \cdot y) =  -(d_r^{p,q} x) \cdot y + (-1)^{n} x \cdot d_r^{s,t}y$. Consider the diagram
    \[\begin{tikzcd}
    	{\pi_n({W_{q}, W_{q+r}}) \otimes \pi_m(X_{t}, X_{t+r})} & {\pi_{n+m}(W_{q} \wedge X_{t},  W_{q+r} \wedge X_{t} \coprod^{\ldots} W_{q} \wedge X_{t+r})} \\
    	{\pi_{n+m-1}(W_{q+r} \wedge X_{t} \coprod^{\ldots} W_{q} \wedge X_{t+r})} & {\pi_{n+m}(Y_{q+t}, Y_{q+t+r})} \\
    	{\pi_{n+m-1}(W_{q+r} \wedge X_{t} \coprod^{\ldots} W_{q} \wedge X_{t+r}, W_{q+r} \wedge X_{t+r})} & {\pi_{n+m-1}(Y_{q+t+r})} \\
    	& {\pi_{n+m-1}(Y_{q+t+r}, Y_{q+t+r+1}).}
    	\arrow["\circ", from=1-1, to=1-2]
    	\arrow["\kappa", from=1-2, to=2-1]
    	\arrow["{\rho_*}", from=1-2, to=2-2]
    	\arrow["{j_*}", from=2-1, to=3-1]
    	\arrow["{\rho_*}"{description}, from=2-1, to=3-2]
    	\arrow["\kappa", from=2-2, to=3-2]
    	\arrow["d"{pos=0.6}, shift left=2, curve={height=-30pt}, from=2-2, to=4-2]
    	\arrow["{s_*}"{description}, from=3-1, to=4-2]
    	\arrow[from=3-2, to=4-2]
    \end{tikzcd}\]
Here the top arrow of the square $s$ is $W_{q+r} \wedge X_{t+r} \to Y_{q+r+t+r} \to Y_{q+t+r+1}$.
If $x = \iota_r(\bar{x}), y = \iota_r(\bar{y})$, then the top-right composite takes $\bar{x}, \bar{y}$ to $d_r^{p+s, q+t} (x \cdot y)$. Both the triangle and the parallelogram commute by the naturality of the long exact sequence of homotopy groups (Proposition \ref{pr:NaturalityLES}).

By Proposition \ref{4.1} (applied with $n-1$ and $m-1$ in place of $n$ and $m$), we know that $j_* \kappa (\bar{x} \circ \bar{y}) = -l_* ((\kappa \bar{x}) \circ \bar{y}) + (-1)^{n}  r_* (\bar{x} \circ \kappa \bar{y})$ for the corresponding $l$ and $r$. So it only remains to show that $s_* l_* ((\kappa \bar{x}) \circ \bar{y}) = (d_r^{p,q} {x}) \cdot {y}$ and $s_* r_* (\bar{x} \circ \kappa \bar{y}) = {x} \cdot d_r^{s,t}{y}.$
This follows from naturality: the first equality comes from the commutativity in

\[\begin{tikzcd}[column sep=small, cells={inner xsep=0.4ex}]
	{\pi_n(W_q, W_{q+r}) \otimes \pi_m(X_t, X_{t+r})} \\
	{\pi_{n-1}(W_{q+r}, \sgtn) \otimes \pi_m(X_t, X_{t+r})} & {\pi_{n-1}(W_{q+r}, W_{q+r+1}) \otimes \pi_m(X_t, X_{t+r})} \\
	{\pi_{m+n-1}(W_{q+r} \wedge X_t, \sgtn \wedge X_t \coprod^{\ldots} W_{q+r} \wedge X_{t+r})} & {\pi_{m+n-1}(W_{q+r} \wedge X_t, W_{q+r+1} \wedge X_t \coprod^{\ldots} W_{q+r} \wedge X_{t+r})} \\
	{\pi_{m+n-1}(W_{q+r} \wedge X_{t} \coprod^{\ldots} W_{q} \wedge X_{t+r}, W_{q+r} \wedge X_{t+r})} & {\pi_{m+n-1}(Y_{q+t+r}, Y_{q+t+r+1}).}
	\arrow["{\kappa \otimes \id}"', from=1-1, to=2-1]
	\arrow["{d \otimes \id}", from=1-1, to=2-2]
	\arrow["{i_* \otimes \id}", from=2-1, to=2-2]
	\arrow["\circ"', from=2-1, to=3-1]
	\arrow["\circ", from=2-2, to=3-2]
	\arrow["{i_*}", from=3-1, to=3-2]
	\arrow["{l_*}"', from=3-1, to=4-1]
	\arrow["{\rho_*}", from=3-2, to=4-2]
	\arrow["{s_*}"', from=4-1, to=4-2]
\end{tikzcd}\]

The left column of this diagram followed by $s_*$ takes $(\bar{x}, \bar{y})$ to $s_* l_*((\kappa \bar{x}) \circ \bar{y})$, while $d \otimes \id$ followed by the right column takes it to $(d_r^{p,q} x) \cdot y$; here $\rho_*$ is induced by the analogue of the map $\rho$ from the definition of the product. The upper triangle commutes by the definition of $d$, and the middle square by the naturality of the pairing (Proposition \ref{circ-prod-nat}). It remains to show that the lower square commutes, i.e. that the two maps $s_* l_*$ and $\rho_* i_*$ from $\pi_{m+n-1}(W_{q+r} \wedge X_t, W_{q+r} \wedge X_{t+r})$ to $\pi_{m+n-1}(Y_{q+t+r}, Y_{q+t+r+1})$ are equal. It is enough to see this equality on the space pair level, i.e. to prove that the following cube has a filler in $U_*$:

\[\begin{tikzcd}
	& { W_{q+r+1} \wedge X_t \coprod^{W_{q+r+1} \wedge X_{t+r}} W_{q+r} \wedge X_{t+r}} && {Y_{q+t+r+1}} \\
	{W_{q+r} \wedge X_{t+r}} && {W_{q+r} \wedge X_{t+r}} \\
	& {W_{q+r} \wedge X_t} && {Y_{q+t+r}} \\
	{W_{q+r} \wedge X_t} && {W_{q+r} \wedge X_{t} \coprod^{\ldots} W_{q} \wedge X_{t+r}.}
	\arrow[from=1-2, to=1-4]
	\arrow[from=1-2, to=3-2]
	\arrow[from=1-4, to=3-4]
	\arrow[from=2-1, to=1-2]
	\arrow[from=2-1, to=2-3]
	\arrow[from=2-1, to=4-1]
	\arrow[from=2-3, to=1-4]
	\arrow["{\operatorname{via } \inl }", from=2-3, to=4-3]
	\arrow[from=3-2, to=3-4]
	\arrow[from=4-1, to=3-2]
	\arrow[from=4-1, to=4-3]
	\arrow[from=4-3, to=3-4]
\end{tikzcd}\]
This is again a formal consequence of pairing coherences.

For $s_* r_* (\bar{x} \circ \kappa \bar{y}) = \bar{x} \cdot d_r^{s,t}\bar{y},$ the proof is analogous, though slightly more involved as the ``front-right'' composite from the cube above involves $\glue$ as a commutativity witness.

This concludes the proof that $d_r$ is a derivation with respect to the multiplication just constructed and thus the proof of the theorem.
\end{proof}

\begin{re}
    One could expect that we'd need the cubes in the second coherence condition of Definition \ref{coherent-pairing} to have based fillers. However, we don't use this anywhere. (At the same time, for applications we only need coherences for the Whitehead tower, which we anyway show to hold in $U_*$ in Lemma \ref{le:Whitehead-coh}.)
    As seen in the proof above, the cubes are used to define certain 2-cells, for which to commute \textit{in} $U$ the 3-fillers are necessary. At the same time, these cells are based by virtue of 2-dimensional coherences: a cell from a pushout is based if and only if its whisker with $\inl$ is based, so the way we glue doesn't matter. This is similar to how we didn't need based 3-fillers in the definitions of the pairing (e.g. Definition \ref{de:circ}).  One could also expect $U_*$-fillers to be needed for parallelepipeds at the end of the proof above, but they are formal.
\end{re}

\begin{re}\label{re:graded-comm-SS}
    Assume that $X \equiv W$ and the diagrams
    \[\begin{tikzcd}
    	{\pi_k (X_p \wedge X_q)} & {\pi_k Y_{p+q}} \\
    	{\pi_k (X_q \wedge X_p)}
    	\arrow[from=1-1, to=1-2]
    	\arrow[from=1-1, to=2-1]
    	\arrow[from=2-1, to=1-2]
    \end{tikzcd}\]
    commute up to sign $(-1)^{s(p,q)}$, and such triangles are compatible with structure maps (in a certain natural sense that we won't specify). Then from Proposition \ref{pr:graded-comm} it follows that the product $E^{p,q}(X_\bullet) \otimes E^{s,t}(X_\bullet) \to E^{p+s,q+t}(Y_\bullet)$ commutes up to sign $(-1)^{s(q,t) + (p+q-1)(t+s-1)+1}$.
\end{re}

\subsubsection{Convergence}\label{conv}

We now discuss the convergence of the product.

    Recall that \cite{Flo} gives a convergence statement for spectral sequences associated to ``lim-towers'' -- that is, towers such that $Y_n = \sgtn$ for all $n$ sufficiently small and $Y_n = Y$ for all $n$ sufficiently large. (Actually, for a slightly broader class of towers: we only need the corresponding equalities on homotopy groups.)

    To state the convergence of the product we first switch to ``colim-towers'' -- towers such that $Y_n = Y$ for all $n$ sufficiently small and $Y_n = \sgtn$ for all $n$ sufficiently large. We will default to using upper indexing for such towers: $\ldots \to Y^k \to Y^{k-1} \to \ldots $.

    Convergence for the associated spectral sequence looks as follows. Let $D_r^{p,q} :\equiv \im (\pi_{p+q}(Y^{q+r-1}) \to \pi_{p+q}(Y^{q}))$. From the exact couple construction we get, for $p+q \geq 1,$ the exact sequence
    $$E_r^{p+r+2, q-r} \xrightarrow{k} D_r^{p+r, q-r+1} \xrightarrow{i} D_r^{p+r+1, q-r} \xrightarrow{j} E_r^{p+2,q-1} \xrightarrow{k} D_r^{p,q}.$$
    \\The bounding terms $E_r^{p+r+2, q-r}$ and $D_r^{p,q}$ are eventually (as $r$ grows, for given $p$ and $q$) zero. Hence, instead of a cofiltration (as was the case for lim-towers) we get a filtration $F_{\bullet, \bullet}$ on $\pi_\bullet(Y)$ and convergence of the form $0 \to F \to F \to E_\infty \to 0$, which is the standard setup to state the convergence of the product.

    To formulate the convergence statement precisely, assume $Y^0 = Y$ and $\pi_{p+q}(Y^i) = 0$ for all $i > T$. Define the types $F_{p, q} :\equiv F_{p,q}(Y^\bullet) :\equiv \lim\limits_{r \to \infty}D_r^{p+r-1, q-r+1} \equiv \lim\limits_{r \to \infty}\im(\pi_{p+q}(Y^{q}) \to \pi_{p+q}(Y^{q-r+1})) = \im(\pi_{p+q}(Y^{q}) \to \pi_{p+q}(Y))$.
    Then we get sequences
    \begin{align*}
        0 \to F_{p+q+1-T, T} \xrightarrow{\sim} E_\infty^{p+q+1-T,T} \to 0 \\
        \ldots \\
        0 \to F_{p, q+1} \to F_{p+1, q} \to E_\infty^{p+1,q} \to 0 \\
        \ldots \\
        0 \to F_{p+q, 1} \to \pi_{p+q+1}(Y) \to E_\infty^{p+q+1,0} \to 0.
    \end{align*}
Thus, we get
\begin{prp}
    Assume given a colim-tower $Y^\bullet$ with colimit $Y$ such that for any $n$, the group $\pi_{n}(Y^r)$ is trivial for all sufficiently large $r$. Then there is an associated converging spectral sequence
    $$E_{\bullet}^{p,q}(Y^\bullet) \implies \pi_{p+q}(Y)$$
    with $E_{1}^{p,q}(Y^\bullet) \equiv \pi_{p+q}(Y^q, Y^{q+1})$, where the convergence means the existence of a filtration $\ldots \subset F_{p+q-t, t } \subset \ldots \subset  F_{p+q-1, 1} \subset F_{p+q, 0} =\pi_{p+q}(Y)$ fitting in a chain of short exact sequences as above.
\end{prp}
\begin{thm}\label{th:prod-convergence}
    Assume given towers $X^\bullet, W^\bullet, Y^\bullet$ with colimits $X, W, Y$ as in the above proposition, and a coherent pairing $W^\bullet \wedge X^\bullet \to Y^\bullet$. Then one gets a multiplication $F_{p,q}(W^\bullet) \otimes F_{p',q'}(X^\bullet) \to F_{p+p', q+q'}(Y^\bullet)$.
    Moreover, it descends to a multiplication on $F_{\bullet,\bullet}/F_{\bullet-1,\bullet+1}$, which coincides with the multiplication on $E_\infty^{\bullet,\bullet} \simeq F_{\bullet,\bullet}/F_{\bullet-1,\bullet+1}$.
\end{thm}

\begin{proof}
    This follows from the description of $E_r$ as $Z_r/B_r$. Namely, there is a multiplication on every $F_{\bullet, \bullet}$ induced from the pairing
    $$\pi_{p+q}(W) \times \pi_{p'+q'}(X^{}) \xrightarrow{\sim} \pi_{p+q}(W, \sgtn) \times \pi_{p'+q'}(X^{}, \sgtn) \to \pi_{p+p'+q+q'}(Y^{}, \sgtn) \xrightarrow{\sim} \pi_{p+p'+q+q'}(Y^{}).$$
    (The check here repeats the one from the construction of Theorem \ref{th:product-SS}.) By the naturality of the product, the maps of the exact sequence $0 \to  F_{\bullet-1, \bullet+1} \to  F_{\bullet, \bullet} \to E^{\bullet, \bullet}_{\infty} \to 0$ preserve it.
\end{proof}

\subsubsection{Pairings on Whitehead towers}

The standard example of a bounded tower is the Postnikov tower of a truncated pointed type

$$Y = Y_q \to Y_{q-1} \to \ldots \to Y_{-1} = \sgtn$$
\noindent
with $Y_i :\equiv ||Y||_i$.
Note that the type of interest is the limit of this sequence; Postnikov towers are lim-towers, while our convergence statement is for {colim}-towers. On the other hand,
given a pairing $W \wedge X \pto Y$, there is no obvious way to define a pairing on Postnikov towers $W_\bullet, X_\bullet, Y_\bullet.$ This motivates switching to \textit{Whitehead towers}, which are obtained by substituting connected covers for truncations:
$$\sgtn = Y^{q+1} \to Y^{q} \to Y^{q-1} \to \ldots \to Y^{0} = Y.$$
Here $Y^i$ is defined as the fiber of $Y \to Y_{i-1}$.

Now in the square (in $U_*$)
\[\begin{tikzcd}
	{W^q \wedge X^t} & {Y^{q+t}} \\
	{W \wedge X } & Y
	\arrow[from=2-1, to=2-2]
	\arrow[from=1-1, to=2-1]
	\arrow[from=1-2, to=2-2]
	\arrow[dashed, from=1-1, to=1-2]
\end{tikzcd}\]
the dashed lift exists because the domain is $(q+t-1)$-connected (\cite{Bru}, Proposition 4.3.1, \cite{Hurewicz}, Lemma 2.15).

Note that such a passage from a lim-tower to a colim-tower is always possible:
\begin{de}
    Given a lim-tower $Y_{\bullet}$, its \textit{associated colim-tower} $Y^\bullet$ is defined by $Y^n :\equiv \fib((\lim Y_{\bullet}) \to Y_{n-1})$.
\end{de}
Moreover, the resulting spectral sequences are the same, as we now show.

\begin{lem}\label{le:fiber-grid}
    Assume that in a commutative diagram in $U_*$
\[\begin{tikzcd}
	{\Omega F} & {F'} & \sgtn & F \\
	{\Omega A} & {A'} & X & A \\
	{\Omega B} & {B'} & X & B
	\arrow["i", from=1-1, to=1-2]
	\arrow[from=1-1, to=2-1]
	\arrow[from=1-2, to=1-3]
	\arrow[from=1-2, to=2-2]
	\arrow[from=1-3, to=1-4]
	\arrow[from=1-3, to=2-3]
	\arrow[from=1-4, to=2-4]
	\arrow[from=2-1, to=2-2]
	\arrow[from=2-1, to=3-1]
	\arrow[from=2-2, to=2-3]
	\arrow[from=2-2, to=3-2]
	\arrow[from=2-3, to=2-4]
	\arrow["\equiv", from=2-3, to=3-3]
	\arrow[from=2-4, to=3-4]
	\arrow[from=3-1, to=3-2]
	\arrow[from=3-2, to=3-3]
	\arrow[from=3-3, to=3-4]
\end{tikzcd}\]
    the columns and the lower two rows are exact sequences of pointed types, and $i$ is the induced fiber map for the lower left square. Then $i$ is an equivalence.
\end{lem}
\begin{proof}
    Follows directly from Lemma \ref{le:octahedral}.
\end{proof}

\begin{de}
    \textit{A morphism (resp.\ isomorphism) of spectral sequences} $E^{p,q}_r \to E'^{p,q}_r$ is a collection of morphisms (resp.\ isomorphisms) $E^{p,q}_r \to E'^{p,q}_r$ for all $p,q,r$ that commute with differentials.
\end{de}

\begin{de}
    Given two towers $\ldots \to Y_q \to Y_{q-1} \to \ldots$ and $\ldots \to Y'_q \to Y'_{q-1} \to \ldots,$ a \textit{map of towers} is a sequence of maps $Y_q \pto Y'_q$ that commute with structure maps.
\end{de}

\begin{prp}
    A map of towers $Y_\bullet \to Y'_\bullet$ induces a morphism of spectral sequences.
\end{prp}

\begin{proof}
    Indeed, a map of towers gives rise to a map of the corresponding exact couples, because the terms and maps there are defined as homotopy groups and maps from long exact sequences, which are natural with respect to pointed maps of (pairs of) types.

    It is an exercise in algebra to show that a map of exact couples induces a map of spectral sequences.
\end{proof}

Now we notice that there is a map of towers $\Omega Y_{\bullet-1} \to Y^\bullet$ and get the following result:
\begin{thm}
    Given a lim-tower $Y_\bullet$ with limit $Y$, the associated colim-tower $Y^\bullet$ has colimit $Y$ and there is a map of towers $\Omega Y_{\bullet-1} \to Y^\bullet$ inducing an isomorphism of spectral sequences $E(Y_\bullet)^{p,q}_r$ and $E(Y^\bullet)^{p-2,q+1}_r$.
\end{thm}
\begin{proof}
    Follows from Lemma \ref{le:fiber-grid} with $X : \equiv Y$, $A :\equiv Y_{\bullet}$, $B :\equiv Y_{\bullet - k}$ (and thus $A' \equiv Y^{\bullet+1}$, $B' \equiv Y^{\bullet-k+1}$) for different $k$, together with the description of $E^{\bullet,\bullet}_\bullet$ from Section \ref{pairing_on_towers}.
\end{proof}
\begin{re}
    Recall that \cite{Flo}'s indexing is $E^{p,q}_1 = \pi_{p+q}(\fib(Y_q \to Y_{q-1}))$. Denote by $\bar{E}^{p,q}_r(Y_\bullet)$ the corresponding spectral sequence. Then we get an isomorphism $\bar{E}(Y_\bullet)^{p,q}_r \to E(Y^\bullet)^{p,q}_r$.
\end{re}

Whitehead towers are good enough for our goal of constructing multiplicative spectral sequences:
\begin{lem}\label{le:Whitehead-coh}
    Given a pairing $W \wedge X \pto Y$ of pointed types,
    the induced pairing on Whitehead towers $\depfun{p,q: \mathbb{N}}{(W^{p} \wedge X^{q} \pto Y^{p+q})}$
    satisfies the coherence conditions from Theorem \ref{th:product-SS}.
\end{lem}
\begin{proof}
    Denote $||Y||_{p+q-1}$ by $Y_1$, $||Y||_{p+q-2 }$ by $Y_2$, $W^{p} \wedge X^q$ by $A$ and $W^{p-1} \wedge X^{q}$ by $B$. Then for $cohL$ we need the dashed cell in the diagram
    \[\begin{tikzcd}
    	A && Y & {Y_1} \\
    	B & {F_1} & Y & {Y_2.} \\
    	& {F_2}
    	\arrow[from=1-1, to=1-3]
    	\arrow[from=1-1, to=2-1]
    	\arrow[dashed, from=1-1, to=2-2]
    	\arrow[from=1-3, to=1-4]
    	\arrow["\equiv", from=1-3, to=2-3]
    	\arrow[from=1-4, to=2-4]
    	\arrow[from=2-1, to=2-3]
    	\arrow[dashed, from=2-1, to=3-2]
    	\arrow[from=2-2, to=1-3]
    	\arrow[from=2-2, to=3-2]
    	\arrow[from=2-3, to=2-4]
    	\arrow[from=3-2, to=2-3]
    \end{tikzcd}\]
    Here $F_i \to Y \to Y_i$ are two fiber sequences. We equate the two maps $A \pto F_2$ by appealing to the uniqueness of the induced pointed map into a fiber. For this we observe that 1) their composites with the map $F_2 \pto Y$ are equal; 2) the trivializations of composites $A \pto Y_2$ are equal. Indeed, the first holds because all involved sides of the prism commute in $U_*$. The second holds because the type $A \pto Y_2$ is contractible (for connectedness reasons), so all proofs that a pointed map from $A$ to $Y_2$ is trivial are equal.

    For the 3-cell, denote $W^q \wedge X^p$ by $D$. We have two elements $x,y: D \to Y^{q+p-2}$ and we have to equate certain equalities $r,s$ between them. For pointed types $A, B$, denote by $A^B$ the type $B \pto A$. Note that the homotopies $r$ and $s$ happen to be based; it is enough to prove $r \cdot s^{-1} =_{\Omega (Y^{q+p-2})^{D} } \refl $\footnote{Note that we shift to seeking pointed fillers now, to be able to appeal to a fiber sequence; accordingly, in what follows we use that we can juxtapose 3-fillers in $U_*$, which is possible in any wild category with a pentagon -- see Remark \ref{re:pentagon}.}. There is a fiber sequence
    $$(\Omega^2 ||Y||_{q+p-3})^D \to \Omega (Y^{q+p-2})^{D} \to \Omega Y^D.$$
    Note that $D$ is ($q+p-1$)-connected, while $\Omega^2 ||Y||_{q+p-3}$ is ($q+p-5$)-truncated. Hence, the fiber in the above sequence is trivial, and the desired equality follows from the triviality of the image of $r \cdot s^{-1}$ in $\Omega Y^D$, which we now proceed to show.

    Consider the square cone

\[\begin{tikzcd}
	&& Y \\
	{Y^{p+q}} & {Y^{p+q-1}} \\
	{Y^{p+q-1}} & {Y^{q+p-2}.}
	\arrow[from=2-1, to=1-3]
	\arrow[from=2-1, to=2-2]
	\arrow[from=2-1, to=3-1]
	\arrow[from=2-2, to=1-3]
	\arrow[from=2-2, to=3-2]
	\arrow[from=3-1, to=1-3]
	\arrow[from=3-1, to=3-2]
	\arrow[from=3-2, to=1-3]
\end{tikzcd}\]
It has an obvious filler. Add it to the cube we are trying to fill, to get the diagram

\[\begin{tikzcd}
	{W^{q} \wedge X^{p}} &&&& {W^{q} \wedge X^{p-1}} \\
	& {Y^{p+q}} && {Y^{p+q-1}} \\
	&& Y \\
	& {Y^{p+q-1}} && {Y^{q+p-2}} \\
	{W^{q-1} \wedge X^{p}} &&&& {W^{q-1} \wedge X^{p-1}.}
	\arrow[from=1-1, to=1-5]
	\arrow[from=1-1, to=2-2]
	\arrow[from=1-1, to=5-1]
	\arrow[from=1-5, to=2-4]
	\arrow[from=1-5, to=5-5]
	\arrow[from=2-2, to=2-4]
	\arrow[from=2-2, to=3-3]
	\arrow[from=2-2, to=4-2]
	\arrow[from=2-4, to=3-3]
	\arrow[from=2-4, to=4-4]
	\arrow[from=4-2, to=3-3]
	\arrow[from=4-2, to=4-4]
	\arrow[from=4-4, to=3-3]
	\arrow[from=5-1, to=4-2]
	\arrow[from=5-1, to=5-5]
	\arrow[from=5-5, to=4-4]
\end{tikzcd}\]
\\
     The filler for the back square cone we've just added,
     together with the above argument, reduces our task to finding a filler for the outer square cone, with tip $Y$ and vertices $W^\bullet \wedge X^\bullet$.

     We can also add the tip $W \wedge X$ to the front (big) side of the square, that is, append the cone

\[\begin{tikzcd}
	{W^{q} \wedge X^{p}} && {W^{q} \wedge X^{p-1}} \\
	& {W \wedge X} \\
	{W^{q-1} \wedge X^{p}} && {W^{q-1} \wedge X^{p-1}.}
	\arrow[from=1-1, to=1-3]
	\arrow[from=1-1, to=2-2]
	\arrow[from=1-1, to=3-1]
	\arrow[from=1-3, to=2-2]
	\arrow[from=1-3, to=3-3]
	\arrow[from=3-1, to=2-2]
	\arrow[from=3-1, to=3-3]
	\arrow[from=3-3, to=2-2]
\end{tikzcd}\]
    It has a filler by Lemma \ref{le:smashes-cube}.

    Then our goal decomposes into four prisms such as

\[\begin{tikzcd}
	{W^q \wedge X^{p}} && {W^q \wedge X^{p-1}} \\
	& {W \wedge X} \\
	{Y^{p+q}} && {Y^{p+q-1}} \\
	& {Y.}
	\arrow[from=1-1, to=1-3]
	\arrow[from=1-1, to=2-2]
	\arrow[from=1-1, to=3-1]
	\arrow[from=1-3, to=2-2]
	\arrow[from=1-3, to=3-3]
	\arrow[from=2-2, to=4-2]
	\arrow[from=3-1, to=3-3]
	\arrow[from=3-1, to=4-2]
	\arrow[from=3-3, to=4-2]
\end{tikzcd}\]

The existence of a filler for this prism follows from the definition of the back square. It was defined by appealing to the uniqueness of the induced fiber map; such uniqueness is proven by eliminating in a family of equality types. The filler we seek then witnesses the $\beta$-reduction for this map.
\end{proof}

\section{Multiplication on the Serre spectral sequence}\label{se:SSS}

Recall that in homotopy type theory one defines, for a ring $G$, the cohomology ring of a type $X$ by $H^n(X, G) :\equiv ||X \to K_nG||_0$, where $K_nG$ is the $n$-th Eilenberg--Mac Lane space. The set $H^{\bullet}(X, G)$ inherits additive and multiplicative structures from the corresponding operations on $K_{\bullet}G$.
Namely, for any $n,m$ there is addition $K_n G \times K_n G \xrightarrow{+} K_n G$ and multiplication $K_n G \times K_m G \xrightarrow{\cup} K_{n+m} G$, so we can define for example the cohomology cup-product by sending $(|a|_0, |b|_0)$ to the truncation of $X \xrightarrow{\Delta} X \times X \xrightarrow{a \times b} K_{n} G \times K_{m} G \xrightarrow{\cup} K_{n+m} G$.

\cite{Flo} (\cite{NLab}) shows that given a fibration (i.e. a type family) $F: X \to U$, one can construct a Federer spectral sequence that turns out to be the Serre spectral sequence $H^p(X , H^q(F_x, G)) \implies H^{p+q}(\tot{x:X}{F_x}, G)$. Note that the $E_1$ page as well as the target have multiplicative structures, which are known to be connected in classical algebraic topology. One could then hope that our pairing can be specialized to these products.

When constructing the Serre spectral sequence, \cite{Flo} uses spectra in order to have the delooping functor: given a type $X$, one can encapsulate its $G$-cohomology in the spectrum $\Sigma^{\infty} X \to HG$, taking the ($-n$)-th homotopy group of which produces the $n$-th cohomology group of $X$. Our framework has been developed for types only. However, we can encapsulate the $G$-cohomology of a type $X$ up to dimension $N$ in the type $C :\equiv (X \to K_N G)$. Indeed, due to the equivalences $\sigma_n: K_nG \to \Omega K_{n+1}G,$ we have $\pi_k(C) = H^{N-k}(X, G)$. Moreover, this keeps track of the addition (\cite{LjuMor}, Lemma 20) as well as the product in the cohomology ring:

\begin{lem}\label{le:29}
The diagram
    \[\begin{tikzcd}
    	{K_n G \times K_m G} & {K_{n+m} G} \\
    	{\Omega K_{n+1}G \times  K_{m}G} & {\Omega K_{m+n+1}G}
    	\arrow[from=1-1, to=1-2]
    	\arrow["{\sigma_n \times K_m G}"', from=1-1, to=2-1]
    	\arrow["{\sigma_{n+m}}", from=1-2, to=2-2]
    	\arrow[from=2-1, to=2-2]
    \end{tikzcd}\]
    commutes.
    In other words, the pointwise looping of the adjoint of the product map $K_m G \to K_{n+1} G \to K_{n+m+1}G$ coincides with the adjoint of the ``previous'' product map $K_m G \to K_{n} G \to K_{n+m}G$. (Note that it makes sense to speak of pointwise looping because the space-level cup-product factors through the smash-product by  Proposition 25 from \cite{LjuMor}.)
\end{lem}
\begin{proof}
    This is Lemma 29 from \cite{LjuMor}.
\end{proof}
\begin{re}
    Note that since the cup-product is graded-commutative, a similar diagram with $\sigma$ applied to the second coordinate commutes only up to $(-1)^n$. Intuitively, this is because one has to permute $\Omega$ and the $n$-dimensional $K_n G$.
\end{re}

Thus, our setup is as follows. Assume given a type $F$.
Let $Y :\equiv K_N G, Z :\equiv K_{2N}G$. Then there is a pairing $Y \wedge Y \to Z$. This gives rise to a map $\pi_n(F \to Y) \otimes \pi_m(F \to Y) \to \pi_{n+m}(F \to Z)$, which, as we show below, is equivalently the cup-product $H^{N-n}(F, G) \otimes H^{N-m}(F, G) \to H^{2N-n-m}(F, G)$. Now let this construction be parametrized by a type $X$:
assume given a fibration $F: X \to U$ and for $A \equiv Y,Z$ take $A :\equiv \depfun{x:X}{A_x}$ where $Y_x :\equiv F_x \to K_N G$ and $Z_x :\equiv F_x \to K_{2N} G$. Consider the Whitehead tower $Y_x^\bullet$ for $Y_x$ (so that $Y_x^q$ is the ($q-1$)-connected cover of $Y_x$). This gives a colim-tower for $Y$ by $Y^q :\equiv \depfun{x:X}{Y^q_x}$. Similarly, we get the tower $Z^\bullet$. There is an obvious pairing $Y^p \wedge Y^q \to Z^{p+q}$, which fulfills the coherence conditions from Theorem \ref{th:product-SS}.
\begin{prp}
    The pairing $Y^\bullet \wedge Y^\bullet \to Z^\bullet$ induced from $K_NG \wedge K_NG \to K_{2N}G$ is coherent (see Definition \ref{coherent-pairing}).
\end{prp}
\begin{proof}
    Follows from Lemma \ref{le:Whitehead-coh}: if $A_x \wedge B_x \to C_x$ is coherent for every $x: X$ then $(\depfun{x:X}{A_x}) \wedge (\depfun{x:X}{B_x}) \to \depfun{x:X}{C_x}$ is coherent by $\funext$ for the 2-cell and naturality of $\funext$ for the 3-cell. (In our case, $A_x \equiv B_x \equiv Y_x$ and $C_x \equiv Z_x$.)
\end{proof}

Note that $\fib({Y^{q+1}_x \to Y^{q}_x}) \simeq K_{q-1}\pi_qY_x,$ and since $\depfun{x:X}{}$ commutes with taking fibers, $\fib({Y^{q+1} \to Y^q})$ equals $\depfun{x:X}{K_{q-1}\pi_qY_x}$.
    Denote $p' :\equiv -p, q' :\equiv N - q$. Then there is a chain of group isomorphisms
    \begin{align*}
    E_1^{p',q'}(Y^{\bullet}) &\equiv \pi_{p'+q'}(Y^{q'}, Y^{q'+1}) \simeq \pi_{p'+q'-1}(\depfun{x:X}{K_{q'-1}\pi_{q'}(F_x \to K_N G)}) \\
    &\simeq ||\depfun{x:X}{\Omega^{p'+q'-1}K_{q'-1}\pi_{q'}(F_x \to K_N G)}||_0 \simeq ||\depfun{x:X}{K_p\pi_{q'}(F_x \to K_N G)}||_0 \\
    &\simeq ||\depfun{x:X}{K_p\pi_{0}(F_x \to K_{N-q'}G)}||_0 \equiv H^p(X, \lambda x. H^q(F_x, G)).
    \end{align*}
    Similarly, for $p'' :\equiv -p, q'' :\equiv 2N - q$, $E^{p'',q''}_1(Z^\bullet) \simeq H^p(X, \lambda x. H^q(F_x, G))$.

Likewise, for the target we have
$$\pi_{p+q}(Y) = ||\depfun{x:X}{\Omega^{p+q}(F_x\to K_NG}) ||_0 \simeq ||\depfun{x:X}{(F_x\to K_{N-p-q}G}) ||_0 \simeq ||(\tot{x:X}{{F_x}) \to K_{N-p-q}G}||_0 = H^{p'+q'}(\tot{x:X}{F_x}, G)$$
and $\pi_{p+q}(Z) \simeq H^{p''+q''}(\tot{x:X}{F_x}, G).$

Finally, recall that a product structure on cohomology with local coefficients can be defined as follows. Suppose given a type $X$ and a family $\lambda x.R_x: X \to \operatorname{Rings}$; by definition, $H^\bullet(X, \lambda x.R_x)$ is $||\depfun{x:X}{K_\bullet R_x}||_0$. Then define the product in the obvious way by $\lambda(|f|, |g|).|\mu^{p,q}_x \circ (f, g) \circ \Delta|$, where $f: \depfun{x:X}{K_p R_x}$, $g: \depfun{x:X}{K_q R_x}$, $\Delta: X \to X \times X$ is the diagonal, and $\mu^{p,q}_x: K_p R_x \times K_qR_x \to K_{p+q}
R_x$ is induced from $R_x \otimes R_x \to R_x$. For a constant type family $R_x \equiv R$ this is the definition of the standard cup-product. Now, since $\lambda x. H^\bullet(F_x, R)$ is a family of (graded) rings, this yields a cup-product structure on $H^\bullet(X, \lambda x. H^\bullet(F_x, R))$.

Thus, there are two products $E_1^{p,q}(Y^\bullet) \times E_1^{s,t}(Y^\bullet) \to E_1^{p+s,q+t}(Z^\bullet)$: one coming from the pairing of towers, as in the general theory, and the other from its identification with $H^{-p}(X , \lambda x.H^{N-q}(F_x, G)) \times H^{-s}(X , \lambda x.H^{N-t}(F_x, G)) \to H^{-p-s}(X , \lambda x.H^{2N-q-t}(F_x, G))$.

\begin{thm}\label{th:products-coincide}
    These two products are equivalent up to a sign.
\end{thm}

It will make sense to speak of the precise sign once we have fixed isomorphisms $E_1^{p',q'}(Y^\bullet) \simeq H^p(X, \lambda x. H^q(F_x, G))$ and $E^{p'',q''}_1(Z^\bullet) \simeq H^p(X, \lambda x. H^q(F_x, G))$. Namely, we need to fix, for a type $W$, the equivalence $\fib({W^{q+1} \to W^{q}}) \to K_{q-1}\pi_qW$, where $W^n$ is the $(n-1)$-connected cover of $W$.
Let us do it now.

For any tower $X^\bullet$, denote by $F^X_q$ the fiber of the structure map $t_q^X: X^{q+1} \to X^q$. When $X^\bullet$ is the Whitehead tower of a pointed type $X$, we want to fix a pointed equivalence $f^X_q$ between $F^X_q$ and $K_{q-1}\pi_qX$. These are ($q-2$)-connected and ($q-1$)-truncated types, so by \cite{Warn}, Corollary 12, it is determined by its action on loop spaces. We define $f^X_q$ as the (essentially unique) map whose $(q-1)$-th looping postcomposed with $\sigma_0^{-1} \circ \Omega \sigma_1^{-1} \circ \ldots \circ \Omega^{q-2} \sigma_{q-2}^{-1}: \Omega^{q-1}K_{q-1}\pi_{q} X \to \pi_{q} X$ equals the equivalence $\Omega^{q-1}F_{q}^{X} \xrightarrow{||-||_0} \pi_{q-1}F_{q}^{X} \to \pi_{q}X^{q} \xrightarrow{\sim} \pi_{q}X$ where the second map is the inverse to the isomorphism from the long exact sequence of $F_q^X \to X^{q+1} \to X^q$.

For a Whitehead tower $X^{\bullet}$, fix also the equivalence $e_q^X: \Omega F_q^X \to F_{q-1}^{\Omega X}$ defined as follows. For any $k$, there is an equivalence $\Omega X^k \to (\Omega X)^{k-1}$ defined by the universal property as the (essentially unique) map whose composition with the $(k-2)$-connected covering projection $(\Omega X)^{k-1} \to \Omega X$ is $\Omega$ of the $(k-1)$-connected covering projection $X^k \to X$. This produces an equivalence of pointed maps $\Omega t_q^X \to t_{q-1}^{\Omega X}$, whose induced fiber map we set to be $e^X_q$.

It is of course necessary to relate these maps:
\begin{lem}\label{le:twoangles}
    For any $q: \mathbb{N}$ and $Y:U$,
    we have
    $\Omega f^Y_q = - (\sigma_{q-2} \circ f^{\Omega Y}_{q-1} \circ e^Y_q).$
\end{lem}
\begin{proof}
    The relation is seen to hold upon applying $\Omega^{q-2}$ and whiskering with $\sigma_0^{-1} \circ \Omega \sigma_1^{-1} \circ \ldots \circ \Omega^{q-2} \sigma_{q-2}^{-1}: \Omega^{q-1}K_{q-1}\pi_qY \to \pi_q Y$.
    Indeed, $f^{\Omega Y}_{q-1}$ is defined as the map whose $(q-2)$-th looping postcomposed with $\sigma_0^{-1} \circ \Omega \sigma_1^{-1} \circ \ldots \circ \Omega^{q-3} \sigma_{q-3}^{-1}: \Omega^{q-2}K_{q-2}\pi_{q-1} \Omega Y \to \pi_{q-1} \Omega Y$ equals the map $\Omega^{q-2}F_{q-1}^{\Omega Y} \xrightarrow{\sim} \pi_{q-2}F_{q-1}^{\Omega Y} \to \pi_{q-1}(\Omega Y)^{q-1} \xrightarrow{\sim} \pi_{q-1}\Omega Y$. We need to show that, precomposed with $\Omega^{q-2} e^Y_q$, this coincides with the description of $- \sigma_0^{-1} \circ \Omega \sigma_1^{-1} \circ \ldots \circ \Omega^{q-2} \sigma_{q-2}^{-1} \circ \Omega^{q-1} f^Y_q$, which is $\Omega^{q-2} \Omega F^Y_q \xrightarrow{\sim} \pi_{q-2}\Omega F^Y_q \to \pi_{q-1} \Omega Y^{q} \to \pi_{q-1}(\Omega Y)^{q-1} \to \pi_{q-1}\Omega Y \equiv \pi_{q}Y$ (see the comment on signs below).
    Extend the map of fiber sequences defining $e^Y_q$ to a map of long fiber sequences:
    \[\begin{tikzcd}
    	{\Omega^{q} Y^{q}} & {\Omega^{q-1} F^{Y}_q} & \ldots & {\Omega F^{Y}_q} & {\Omega Y^{q+1}} & {\Omega Y^{q}} \\
    	{\Omega^{q-1}(\Omega Y)^{q-1}} & {\Omega^{q-2} F^{\Omega Y}_{q-1}} & \ldots & {F^{\Omega Y}_{q-1}} & {(\Omega Y)^q} & {(\Omega Y)^{q-1}} \\
    	{\Omega^{q} Y.}
    	\arrow["{\delta_{q-2}^{\Omega t_q^Y}}", from=1-1, to=1-2]
    	\arrow[from=1-1, to=2-1]
    	\arrow[from=1-2, to=1-3]
    	\arrow["{\Omega^{q-2} e^Y_{q}}"', from=1-2, to=2-2]
    	\arrow["{\delta_0^{\Omega t_q^Y}}", from=1-3, to=1-4]
    	\arrow[from=1-4, to=1-5]
    	\arrow["{e_q^Y}", dashed, from=1-4, to=2-4]
    	\arrow["{\Omega t_q^Y}", from=1-5, to=1-6]
    	\arrow[from=1-5, to=2-5]
    	\arrow[from=1-6, to=2-6]
    	\arrow["{\delta_{q-2}^{t_{q-1}^{\Omega Y}}}", from=2-1, to=2-2]
    	\arrow[from=2-1, to=3-1]
    	\arrow[from=2-2, to=2-3]
    	\arrow["{\delta_{0}^{t_{q-1}^{\Omega Y}}}", from=2-3, to=2-4]
    	\arrow[from=2-4, to=2-5]
    	\arrow["{t_{q-1}^{\Omega Y}}", from=2-5, to=2-6]
    \end{tikzcd}\]
    Here we see boundary homomorphisms $\delta$ as depending on the function we started with and the number of ``cycles'' we made.
    \\
    The leftmost square commutes. The two paths starting from its upper right corner and ending in $\Omega^q Y$ give us what we want.
    (Note that $\delta_{q-2}^{\Omega t_q^{Y}} = -\delta_{q-1}^{t_{q}^{Y}}$. Thus the top map of the square gives after truncation \textit{minus}
    the connecting homomorphism used in the definition of $f^Y_q$. This minus won't play a role in the future.)
\end{proof}

First we need to show that the pairing $Y \wedge Y \to Z$ indeed encapsulates the product in the inner cohomology groups. This means the following. Consider a type $F$ and $N, M, p, q \geq 0$ with $N \geq q, M \geq p$. There is a map
$$\pi_p (F \to K_M G) \otimes \pi_q (F \to K_N G) \xrightarrow{\circ} \pi_{p+q} ((F \to K_M G) \wedge (F \to K_N G)) \to \pi_{p+q}(F \to K_{M+N} G).$$
Here $\circ$ is the pairing from Definition \ref{de:circ-absolute} and the second map is induced from the (pre-truncation) cup-product map.
\begin{lem}\label{le:pairing-gives-inner-prod}
    The diagram
    \[\begin{tikzcd}
    	{\pi_p (F \to K_M G) \otimes \pi_q (F \to K_N G)} & {\pi_{p+q}(F \to K_{M+N}G)} \\
    	{H^{M-p}(F,G) \otimes H^{N-q}(F,G)} & {H^{N+M-p-q}(F,G)}
    	\arrow[from=1-1, to=1-2]
    	\arrow["\sim", from=1-1, to=2-1]
    	\arrow["\sim"', from=1-2, to=2-2]
    	\arrow["\cup", from=2-1, to=2-2]
    \end{tikzcd}\]
    commutes up to sign $(-1)^{(M-p)q}$.
\end{lem}
\begin{proof}
    We induct on $p, q$. The claim is obvious when $p=q=0$.
    Let the equivalence $\sigma \equiv \sigma(T,k)$ be the composite $(F \to K_{T-k}G) \to (F \to \Omega^k  K_TG) \to \Omega^k (F \to K_TG)$ used for the vertical isomorphisms.
    From Lemma \ref{le:29} it follows that the diagram
    \[\begin{tikzcd}
    	{(F \to K_{M-1}G) \times (F \to K_{N}G)} & {(F \to K_{M+N-1}G)} \\
    	{\Omega(F \to K_{M}G) \times (F \to K_{N}G)} & {(F \to \Omega K_{M+N}G)}
    	\arrow["\cup", from=1-1, to=1-2]
    	\arrow["{\sigma \times \id}"', from=1-1, to=2-1]
    	\arrow["\sigma"', from=1-2, to=2-2]
    	\arrow[from=2-1, to=2-2]
    \end{tikzcd}\]
    commutes.
    Note that precomposing with $\oout$ corresponds to coordinatewise looping as in this diagram (Remark \ref{re:oout-desription}).
    Then the case $p = 1, q = 0$ follows from these observations.

    The case $p = 0, q=1$ follows similarly, but this time we get sign $(-1)^M$ due to the graded-commutativity of the cup-product.

    Denote $F \to K_{M-i}G$ by $K_{-i}$, $F \to K_{N-i} G$ by $K'_{-i}$ and $F \to K_{M+N-i} G$ by $K''_{-i}$ for $i : \NN$.

    Let us induct on $q$ when $p = 0$.
        Consider the diagram
    \[\begin{tikzcd}
    	{K \wedge \Omega^q \Omega K'} & {K_{} \wedge \Omega^qK'_{-1}} & {K_{} \wedge K'_{-q-1}} \\
    	{\Omega^{q} (K \wedge \Omega K')} & {\Omega^{q} (K_{} \wedge K'_{-1})} & {K''_{-q-1}} \\
    	{\Omega^{q+1} (K \wedge K')} & {\Omega^{q} (K''_{-1})} \\
    	& {\Omega^{q+1} (K'').}
    	\arrow[from=1-1, to=1-2]
    	\arrow["{\ooutk{q}_2}"', from=1-1, to=2-1]
    	\arrow["{K \wedge \sigma^{-1}}"', from=1-2, to=1-3]
    	\arrow["{\ooutk{q}_2}"', from=1-2, to=2-2]
    	\arrow[from=1-3, to=2-3]
    	\arrow["{\Omega^q (K \wedge \sigma^{-1})}", from=2-1, to=2-2]
    	\arrow[from=2-1, to=3-1]
    	\arrow[from=2-2, to=3-2]
    	\arrow[from=2-3, to=3-2]
    	\arrow[from=3-1, to=4-2]
    	\arrow[from=3-2, to=4-2]
    \end{tikzcd}\]
    We need to show that the outer square commutes up to sign $(-1)^{M(q+1)}$. The small square commutes by the naturality of $\ooutk{-}_2$ (Lemma \ref{le:ooutk-nat}). The lower trapezoid commutes up to sign $(-1)^{M}$ by the case $(0,1)$. The right trapezoid commutes up to $(-1)^{Mq}$ by the case $(0,q)$, giving us what we want.

    Induction on $p$ when $q = 0$ is the same but without signs. For a general inductive step, assume $p,q>0$ and consider the diagram
    \[\begin{tikzcd}
    	{\Omega^{p}  K \wedge \Omega^q  K'} & {K_{-p} \wedge \Omega^qK'_{}} & {K_{-p} \wedge K'_{-q}} \\
    	{\Omega^{q} (\Omega^p K \wedge K')} & {\Omega^{q} (K_{-p} \wedge K'_{})} & {K''_{-p-q}} \\
    	{\Omega^{q+p} (K \wedge K')} & {\Omega^{q} (K''_{-p})} \\
    	& {\Omega^{q+p} (K'').}
    	\arrow[from=1-1, to=1-2]
    	\arrow["{\ooutk{q}_2}"', from=1-1, to=2-1]
    	\arrow["{K_{-p}\wedge \sigma^{-1}}"', from=1-2, to=1-3]
    	\arrow["{\ooutk{q}_2}"', from=1-2, to=2-2]
    	\arrow[from=1-3, to=2-3]
    	\arrow["{\Omega^q ( \sigma^{-1} \wedge K')}", from=2-1, to=2-2]
    	\arrow[from=2-1, to=3-1]
    	\arrow[from=2-2, to=3-2]
    	\arrow[from=2-3, to=3-2]
    	\arrow[from=3-1, to=4-2]
    	\arrow[from=3-2, to=4-2]
    \end{tikzcd}\]
    By induction, the lower trapezoid commutes, and the right trapezoid commutes up to sign $(-1)^{(M-p)q}$ by the case $(p,q, M) = (0,q, M-p)$.
\end{proof}

Thus, our pairing already knows about the ``inner'' product -- the one in $H^\bullet(F, G)$. We now show the essential part of the theorem: the induced fiber pairing recovers the ``outer'' product -- the one in the cohomology of $X$.
\begin{proof}[Proof of Theorem \ref{th:products-coincide}]
    We will prove the equivalence up to sign $(-1)^{(N-q)s+(t+s)q+(p+q)(t+s+1)}$.
    Denote $\bar{Y}_x \equiv \bar{W}_x :\equiv ({F_x \to K_N G}), \bar{Z}_x :\equiv (F_x \to K_{2N}G)$.

    We first reduce to considering each $x:X$ separately.
    \\
    For families $X \xrightarrow{A,B} U $ and a family of maps $\lambda x.(f_x:A_x \to B_x),$ denote by $\depfun{x:X}{f_x}: \depfun{x:X}{A_x} \to \depfun{x:X}{B_x}$ the function $\lambda g.\lambda x.f_x(g(x))$.
    \\
    For any $x:X$, there is a pairing $p_x: \bar{Y}_x \wedge \bar{W}_x \to \bar{Z}_x$.
    It gives rise to a pairing for colim-towers $Y^q \wedge W^{s} \to Z^{q+s}$, where for $A = W, Y, Z$, we set $A^q :\equiv \depfun{x:X}{\bar{A}^q_x}$ and $\bar{A}_x^q \to \bar{A}_x \to ||\bar{A}_x||_{q-1}$ is a fiber sequence.
    For different $T$, let $t^T_q$ denote the $q$-th structure map $T^{q+1} \to T^{q}$, so that $t^A_q = \depfun{x:X}{{t^{A_x}_q}}$ for $A = Y, W, Z$.
    There is a map of pairs $\mu: P^\times(t_q^Y, t_s^W) \to t_{q+s}^Z$ which on fibers corresponds to a map $\bar{p}: F^Y_q \join F^W_s \to F^{Z}_{q+s}$. By the construction of Theorem \ref{th:defs-of-pairing}, this gives rise to maps $\Omega^{n-1} F^Y_q \times \Omega^{m-1} F^W_s \to \Omega^{m+n-1} F^{Z}_{q+s}$.
    Note that $P^\times(t_q^Y, t_s^W) \to t_{q+s}^Z = \depfun{x:X}{{t_{q+s}^{Z_x}}}$ factors through $\depfun{x:X}{\mu_x}$ with $\mu_x : P^\times({t_q^{Y_x}}, {t_s^{W_x}}) \to {t_{q+s}^{Z_x}}$; hence $F^Y_q \join F^W_s \to F^{Z}_{q+s}$ factors through $\depfun{x:X}{\bar{p}_x}$ with $\bar{p}_x : \depfun{x:X}{({F^{Y_x}_q} \join {F^{W_x}_s} \to {F^{Z_x}_{q+s}}})$ where ${F^{A_x}_q} \equiv \fib({t^{A_x}_q})$ for $A = Y, W , Z$.
    Thus, the map $\Omega^{n-1} F^Y_q \times \Omega^{m-1} F^W_s \to \Omega^{m+n-1} F^{Z}_{q+s}$ also factors through $\depfun{x:X}{(\Omega^{n+m-1}\bar{p}_x \circ I^{Y_x,W_x}_{n-1,m-1}}):
    \depfun{x:X}{\Omega^{n-1} {F^{Y_x}_q} \times \Omega^{m-1} {F^{W_x}_s}} \to \depfun{x:X}{\Omega^{m+n-1} {F^{Z_x}_{q+s}}}$.
    In other words, the pairing between fibers happens pointwise with respect to $x:X$. Similarly, the cohomology multiplication pairing is done pointwise, as well as the identification of fibers of our towers with
    dependent products of Eilenberg--Mac Lane spaces.

    We can now fix $x:X$ and ease the notation by redefining $Y :\equiv \bar{Y}_x$ etc.

    We prove that the diagram
    \[
    \begin{tikzcd}
    	{\Omega^{p+q-1} F^Y_q \times \Omega^{t+s-1} F^W_s} & {\Omega^{p+q+t+s-1} (F_q^Y \join F^W_s)} & {\Omega^{p+q+t+s-1} F^Z_{q+s}} \\
    	{\Omega^{p+q-1}(K_{q-1}\pi_q Y) \times \Omega^{{t+s-1}}(K_{s-1}\pi_s W)} && {\Omega^{p+q+t+s-1}K_{q+s-1}\pi_{q+s} Z} \\
    	{K_{-p} \pi_q Y \times K_{-t} \pi_s W} & {K_{-p-t} (\pi_q Y \otimes \pi_s W)} & {K_{-p-t}\pi_{q+s} Z}
    	\arrow[from=1-1, to=1-2]
    	\arrow[from=1-1, to=2-1]
    	\arrow[from=1-2, to=1-3]
    	\arrow[from=1-3, to=2-3]
    	\arrow[from=2-1, to=3-1]
    	\arrow[from=2-3, to=3-3]
    	\arrow[from=3-1, to=3-2]
    	\arrow[from=3-2, to=3-3]
    \end{tikzcd}
    \]
    commutes up to sign $(-1)^{(t+s)q+(p+q)(t+s+1)}$. Here $W \equiv Y \equiv (F \to K_NG),$ $Z \equiv (F  \to K_{2N}G)$, and the lower right horizontal map is induced from $\pi_q Y \otimes \pi_s W  \to \pi_{q+s}(Y \wedge W)\to \pi_{q+s}Z$ where the first map is the absolute homotopy group pairing from Definition \ref{de:circ-absolute}. Since this pairing is
    $(-1)^{(N-q)s}$
    times the cohomology product pairing, this leads to the sign announced at the beginning of the proof.

    We denote $n:\equiv p+q, m: \equiv t+s$ and induct on $n, m, q, s$ (assuming the statement for all pointed types $Y, W, Z$ with a pairing $Y \wedge W \to Z$).

    Consider the case $n=m=1$. We need to prove that
    \[\begin{tikzcd}
    	{F^Y_q \times F^W_s} & {\Omega(F_q^Y \join F_s^W)} & {\Omega^{} F^Z_{q+s}} \\
    	{K_{q-1}\pi_{q}Y \times K_{s-1} \pi_s W} & {K_{q+s-2}\pi_{q+s}Z} & {\Omega K_{q+s-1}\pi_{q+s}Z}
    	\arrow["{I_{0,0}^{F_q^Y, F_s^W}}", from=1-1, to=1-2]
    	\arrow["{f^Y_q \times f^W_s}", from=1-1, to=2-1]
    	\arrow[from=1-2, to=1-3]
    	\arrow["{\Omega f^Z_{q+s}}", from=1-3, to=2-3]
    	\arrow[from=2-1, to=2-2]
    	\arrow["{\sigma_{q+s-2}}", from=2-2, to=2-3]
    \end{tikzcd}\]
    commutes up to sign $(-1)^{q}$. Now we induct on $q$ and $s$. For the base case, consider the diagram
    \[\begin{tikzcd}
    	{F_1^Y \times F_1^W} & {\pi_0F_1^Y \times \pi_0F_1^W} & {\pi_1(F^Y_1 \join F^W_1)} & {\pi_1 F^Z_2} & {\Omega F_2^Z} \\
    	{\pi_1 (Y^1,Y^2) \times \pi_1 (W^1, W^2)} & {\pi_2(P^\times(t^Y_1, t^W_1))} & {\pi_2(P^\wedge(t^Y_1, t^W_1))} & {\pi_2(Z^2, Z^3)} & {\Omega K_1 \pi_2 Z} \\
    	{\pi_1 Y^1 \times \pi_1 W^1} & {\pi_2(Y^1 \wedge W^1)} & {\pi_2Z^2} && {K_0 \pi_2 Z} \\
    	{\pi_1 Y \times \pi_1 W} & {\pi_2(Y \wedge W)} & {\pi_2Z.}
    	\arrow["\sim", from=1-1, to=1-2]
    	\arrow[from=1-2, to=1-3]
    	\arrow[from=1-3, to=1-4]
    	\arrow["\simeq", from=1-4, to=1-5]
    	\arrow["{\Omega f^Z_2}", from=1-5, to=2-5]
    	\arrow["\sim", from=2-1, to=1-1]
    	\arrow[from=2-1, to=2-2]
    	\arrow[from=2-2, to=1-3]
    	\arrow[from=2-2, to=2-3]
    	\arrow[from=2-3, to=2-4]
    	\arrow[from=2-4, to=1-4]
    	\arrow["{\sigma_0^{-1}}", from=2-5, to=3-5]
    	\arrow["\delta"{description}, curve={height=-12pt}, from=3-1, to=1-2]
    	\arrow["\sim", from=3-1, to=2-1]
    	\arrow[from=3-1, to=3-2]
    	\arrow[from=3-1, to=4-1]
    	\arrow[from=3-2, to=2-3]
    	\arrow[from=3-2, to=3-3]
    	\arrow[from=3-2, to=4-2]
    	\arrow["\delta"{description}, curve={height=-12pt}, from=3-3, to=1-4]
    	\arrow["\sim"', from=3-3, to=2-4]
    	\arrow[from=3-3, to=4-3]
    	\arrow["\equiv"{description}, from=3-5, to=4-3]
    	\arrow[from=4-1, to=4-2]
    	\arrow[from=4-2, to=4-3]
    \end{tikzcd}\]
    The maps labeled $\delta$ are boundary homomorphisms by Lemma \ref{le:LES-isom}.
    The rightmost cell of the diagram is the definition of $f_2^Z$. Thus, we need to show that the outer diagram anticommutes. In the top row, the first trapezoid commutes by Theorem \ref{th:defs-of-pairing} and the second one by naturality and the definition of the map $F^Y_1 \join F^W_1 \to F^Z_2$.
    In the second row, (anti)commutativities are by naturality and the case $(1,1)$ of Proposition \ref{pr:smashing} (which gives $-1$). In the last row, the first square commutes by the naturality of $\circ$ and the second one by the definition of the map $Y^1 \wedge W^1 \to Z^2$ and the naturality of $\circ$.

    Now let us induct on $q$ assuming $q>1$.
    We need to show that the two maps $F^Y_q \pto F^W_s \pto \Omega K_{q+s-1}\pi_{q+s}Z$ are equal, for which it is enough to show that their loopings are equal (\cite{Warn}, Corollary 11).

    For any pointed types $A, B, C$, there is a chain of equivalences
    \begin{align*}
    \Omega A \pto \Omega (B \pto C) \simeq  \\
    \Omega A \pto (B \pto \Omega C) \simeq  \\
    B \pto \Omega A \pto \Omega C,
    \end{align*}
    which sends $\Omega f$ to $\lambda b.\Omega(f(-)(b))$. With this in mind, we fix $w: F_s^W$ and prove the equality of two functions from loop spaces of the form $\lambda p.(\Omega (f(-)(w)))(p)$.

    Key to our inductive argument is the pairing $\Omega Y \wedge W \xrightarrow{\oout} \Omega (Y \wedge W) \to \Omega Z$. By induction, the two maps
    $F^{\Omega Y}_{q-1} \times F_s^W \pto \Omega K_{q+s-2}\pi_{q+s-1}\Omega Z$
    are equal. Note that the induced pairing $\pi_{q-1}\Omega Y \otimes \pi_s W \to \pi_{q+s-1}\Omega Z$ coincides with $\pi_{q} Y \otimes \pi_s W \to \pi_{q+s} Z$.
    \\
    Consider the diagram
\[\begin{tikzcd}
	{\Omega F^Y_q} & {\Omega \Omega(F_q^Y \star F_s^W)} & {\Omega\Omega^{} F^Z_{q+s}} \\
	{ F^{\Omega Y}_{q-1}} & { \Omega(F_{q-1}^{\Omega Y} \star F_s^W)} & {\Omega^{} F^{\Omega Z}_{q+s-1}} \\
	{K_{q-2} \pi_{q-1} \Omega Y} & {K_{q+s-3}\pi_{q+s-1}\Omega Z} & {\Omega K_{q+s-2}\pi_{q+s}Z} \\
	{\Omega K_{q-1} \pi_{q}Y} & {\Omega K_{q+s-2}(\pi_qY \otimes \pi_{s}W)} & {\Omega \Omega K_{q+s-1}\pi_{q+s}Z.}
	\arrow["{\Omega I_{0,0}^{F_q^Y, F_s^W}(-,w)}", from=1-1, to=1-2]
	\arrow["{e^Y_{q}}", from=1-1, to=2-1]
	\arrow["{-\Omega f^{Y}_q}"', curve={height=18pt}, from=1-1, to=4-1]
	\arrow[from=1-2, to=1-3]
	\arrow["{\Omega e^Z_{q+s}}", from=1-3, to=2-3]
	\arrow["{-\Omega \Omega f^Z_{q+s}}", curve={height=-18pt}, from=1-3, to=4-3]
	\arrow["{I_{0,0}^{F_{q-1}^{\Omega Y}, F_s^W}(-,w)}", from=2-1, to=2-2]
	\arrow["{f^{\Omega Y}_{q-1}}", from=2-1, to=3-1]
	\arrow["{\Omega \oout \circ \Omega ({e^Y_{q}}^{-1} \star F_s^W)  }"', from=2-2, to=1-2]
	\arrow[from=2-2, to=2-3]
	\arrow["{\Omega f^{\Omega Z}_{q+s-1}}"', from=2-3, to=3-3]
	\arrow["{- \cup f^W_sw}", from=3-1, to=3-2]
	\arrow["{ \sigma_{q-2}}"', from=3-1, to=4-1]
	\arrow["{\sigma_{q+s-3}}", from=3-2, to=3-3]
	\arrow["{\Omega \sigma_{q+s-2}}"', from=3-3, to=4-3]
	\arrow["{\Omega (- \cup f^W_sw)}", from=4-1, to=4-2]
	\arrow[from=4-2, to=3-3]
	\arrow[from=4-2, to=4-3]
\end{tikzcd}\]
    We are interested in the commutativity of the outer rectangle.
    The middle rectangle commutes up to sign
    $(-1)^{q-1}$
    by induction (on $q$).
    The side ``triangles'' commute by Lemma \ref{le:twoangles}. The lower rectangle commutes by Lemma \ref{le:29}.
    For the upper right square, we prove that the corresponding two maps of arrows are equal, see the $n$-inductive step below.

    It remains to show that the upper left square anticommutes. Note that $I_{0,0}^{-,F_s^W}$ is natural (Lemma \ref{le:i-nat}), so $\Omega ({e^Y_{q}}^{-1} \star F_s^W) \circ I_{0,0}^{F_{q-1}^{\Omega Y}, F_s^W}(-,w) \circ e_q^Y = I_{0,0}^{\Omega F^Y_q, F^W_s}(-, w)$. The remaining triangle anticommutes by Lemma \ref{I00-oout}.

    The $s$-inductive step when $q=1$ is analogous, with two differences to pay attention to in the lower row. First, this time we use the version of Lemma \ref{le:29} with the looping of the second factor, which adds the sign $(-1)^{q-1} = 1$. Second, this time we consider the pairing $Y \wedge \Omega W \to \Omega Z$. The induced pairing  $\pi_{1} Y \otimes \pi_{s-1} \Omega W \to \pi_{s}\Omega Z$ then applies $\oout_2$ $(s-1)$ times, then applies $\oout$ and then $\oout_2$ again. This differs by one transposition from $\pi_{1} Y \otimes \pi_{s} W \to \pi_{s+1} Z$, so an additional minus sign is added for every increase of $s$. In total, we get $(-1)^{q}$.

    Now suppose $n>1$. Note that we are only interested in cases $n\leq q$ (as $n=p+q$ and $p$ is negative). We may thus assume $q > 1$. The inductive step is analogous in spirit to the inductive step in the proof of the base case. Consider the diagram
\[\begin{tikzcd}
	{\Omega^{n-1} F^Y_q \times \Omega^{m-1}F^W_s} & {\Omega^{n+m-2}(\Omega F_q^Y \join F_s^W)} & {\Omega^{n+m-1} (F_q^Y \join F_s^W)} & {\Omega^{n+m-1} F^Z_{q+s}} \\
	{\Omega^{n-2} F_{q-1}^{\Omega Y} \times \Omega^{m-1} F_s^W} & {\Omega^{n+m-2}(F_{q-1}^{\Omega Y} \join F_s^W)} & {\Omega^{n+m-2} F^{\Omega Z}_{q+s-1}} \\
	{K_{q-n}\pi_{q-1}\Omega Y \times K_{s-m} \pi_s W} && {K_{q+s-n-m}\pi_{q+s-1}\Omega Z} \\
	{K_{q-n}\pi_{q}Y \times K_{s-m} \pi_s W} &&& {K_{q+s-n-m}\pi_{q+s}Z,}
	\arrow["{-I_{n-2,m-1}^{\Omega F_q^Y, F_s^W}}", from=1-1, to=1-2]
	\arrow["{\Omega^{n-2} e^Y_q \times \id}"', from=1-1, to=2-1]
	\arrow["{-(\sigma^{-1} \circ \Omega^{n-1}f^{Y}_{q}) \times \sigma^{-1}}"', curve={height=18pt}, from=1-1, to=4-1]
	\arrow["{\Omega^{n+m-2}\oout}", from=1-2, to=1-3]
	\arrow["{\Omega^{n+m-2} (e^Y_{q} \join \id)}"', from=1-2, to=2-2]
	\arrow[from=1-3, to=1-4]
	\arrow["{\Omega^{n+m-2}e^Z_{q+s}}", from=1-4, to=2-3]
	\arrow["{- \sigma^{-1} \circ (\Omega^{n+m-1}f^Z_{q+s})}"', from=1-4, to=4-4]
	\arrow["{I_{n-2,m-1}^{F_{q-1}^{\Omega Y}, F_s^W}}", from=2-1, to=2-2]
	\arrow["{(\sigma^{-1} \circ \Omega^{n-2}f^{\Omega Y}_{q-1}) \times \sigma^{-1}}", from=2-1, to=3-1]
	\arrow[from=2-2, to=2-3]
	\arrow["{\sigma^{-1} \circ \Omega^{n+m-2}f^{\Omega Z}_{q+s-1}}"', from=2-3, to=3-3]
	\arrow[from=3-1, to=3-3]
	\arrow["\equiv"', from=3-1, to=4-1]
	\arrow["\equiv", from=3-3, to=4-4]
	\arrow[from=4-1, to=4-4]
\end{tikzcd}\]
    where by $\sigma$ we mean appropriate compositions of $\Omega^\bullet \sigma_\bullet$.
    We are interested in the outer rectangle.
    The middle rectangle commutes up to sign
    $(-1)^{mq+mn+n-1}$
    by induction (on $n$ and $q$). The upper left rectangle anticommutes by the naturality of $I$. The ``triangles'' on the sides again commute by Lemma \ref{le:twoangles}.

    Before we finish with this diagram, let us note what happens with the signs when $n=1$ and we induct on $m$. We need to show commutativity up to $(-1)^{mq+m+1}$.
    In this case, the signs of the middle and the lower rows behave differently. First, the reduction will be to the case of the previous $m$ and $s$ in the middle row, which will then inductively commute up to sign $(-1)^{(m-1)q+m-1+1} = (-1)^{mq+m+1-(q+1)}$. Second, the pairing $\pi_{q} Y \otimes \pi_{s-1} \Omega W \to \pi_{q+s-1 } \Omega Z$ differs from $\pi_{q} Y \otimes \pi_s W \to \pi_{q+s} Z$ by a transposition of $\oout_2$ and $\ooutk{q}$. The upper left square still anticommutes, so we get in total $(-1)^{mq+m+1}$.

    To finish, it is enough to prove that the diagram
    \[\begin{tikzcd}
    	{\Omega F_q^Y \join F_s^W} & {\Omega (F_q^Y \join F_s^W)} & {\Omega F^Z_{q+s}} \\
    	{F_{q-1}^{\Omega Y} \join F_s^W} & {F^{\Omega Z}_{q+s-1}}
    	\arrow["\oout", from=1-1, to=1-2]
    	\arrow["{e^Y_{q} \join F^W_s}"', from=1-1, to=2-1]
    	\arrow[from=1-2, to=1-3]
    	\arrow["{e^Z_{q+s}}", from=1-3, to=2-2]
    	\arrow[from=2-1, to=2-2]
    \end{tikzcd}\]
    commutes. The ``down-right'' composite is a fiber of a known map. Unfortunately, $\oout$ does not seem to be induced from any obvious function. Note, however, that we have the following factorization:

    \[\begin{tikzcd}
    	& {\Omega F_q^Y \join F_s^W} & {\Omega (F_q^Y \join F_s^W)} & {\Omega F^Z_{q+s}} \\
    	& {\fib P^{\wedge}{(\Omega t_q^Y, t_s^W)}} & {\Omega \fib P^{\wedge}{(t_q^Y, t_s^W)}} \\
    	{F_{q-1}^{\Omega Y} \join F_s^W} & {\fib P^\wedge(t_{q-1}^{\Omega Y}, t_s^W)} & {F^{\Omega Z}_{q+s-1}.}
    	\arrow["\oout", from=1-2, to=1-3]
    	\arrow[from=1-2, to=2-2]
    	\arrow["{e^Y_{q} \join F^W_s}"', from=1-2, to=3-1]
    	\arrow[from=1-3, to=1-4]
    	\arrow[from=1-3, to=2-3]
    	\arrow["{e^Z_{q+s}}", from=1-4, to=3-3]
    	\arrow["i", from=2-2, to=2-3]
    	\arrow[from=2-2, to=3-2]
    	\arrow[from=2-3, to=1-4]
    	\arrow[from=3-1, to=3-2]
    	\arrow[from=3-2, to=3-3]
    \end{tikzcd}\]
    Here the map $i$ is the induced map of fibers. Now, all the maps coming from $\fib P^{\wedge}{(\Omega t_q^Y, t_s^W)}$ are the induced maps of fibers, hence it is enough to prove that the corresponding maps of arrows are equal. We reduce it to showing the commutativity of the diagram

    \[\begin{tikzcd}
    	{P^\wedge(\Omega t_q^Y, t_s^W)} & {\Omega P^\wedge(t_q^Y, t_s^W)} & {\Omega t^Z_{q+s}} \\
    	{P^\wedge(t_{q-1}^{\Omega Y}, t_s^W)} & {t^{\Omega Z}_{q+s-1}} & {\id_{\Omega Z}.}
    	\arrow["\oout", from=1-1, to=1-2]
    	\arrow["{P^\wedge({E}^Y_{q}, \operatorname{id}_{t^W_s})}"', from=1-1, to=2-1]
    	\arrow[from=1-2, to=1-3]
    	\arrow["{{E}^Z_{q+s}}"{description}, from=1-3, to=2-2]
    	\arrow[from=2-1, to=2-2]
    	\arrow[from=2-2, to=2-3]
    \end{tikzcd}\]
    (Here $E^\bullet_\bullet$ are maps of arrows whose induced maps of fibers are $e^\bullet_\bullet$.) To do this, denote $C: \equiv \Omega Z $ and $k:\equiv q+s$ and note that the domain of the domain is a pushout of $(k-1)$-connected spaces, while the domain of $t^{\Omega Z}_{q+s-1}$ is the $(k-1)$-connected cover of $C$. For any pointed space $X$, we have the fiber sequence
    $$\Omega||C||_{{k-1}}^X \to (C^{k})^X \to C^X$$
    where we use exponential notation for the function types. The fiber $\Omega||C||_{{k}-1}^X$ is contractible as soon as $\operatorname{conn}(X) \geq k-2$. Then to equate two elements of $(C^{k})^X$ it is sufficient to equate their images in $C^X$ (note that $(C^{k})^X$ is an H-space, so proving an equality of two of its elements is equivalent to proving that an element equals the basepoint). If we now have as $X$ a pushout of $(k-1)$-connected spaces $L \coprod^M N$, we get a map $(C^k)^{L \coprod^M N} \to C^{L \coprod^M N}$ which is equivalently ${(C^k)}^L \times_{{(C^k)^M}} (C^k)^N \to {C}^L \times_{{C^M}} C^N $. Then this map is an embedding since the maps of the corresponding components are embeddings.
    \\
    We conclude that to equate two maps between the domains, it is enough to equate them after postcomposing to $\Omega Z$.
    \\
    Next, we need to equate two trivializations of these two maps in the codomain of $t^{\Omega Z}_{q+s-1}$. Denote these arrows by $f$ and $g$. We are then equating two equalities $\alpha, \beta : f = g$, which is equivalent to proving $\alpha \cdot \beta^{-1} = \refl_{f}$. This is equivalent to proving a certain equality in $(\Omega C^{k-1})^{X}$, where $X$ is the domain of the domain of $P^\wedge(\Omega t_q^Y, t_s^W)$.
    We get a fiber sequence
    $$\Omega||\Omega C||_{{k-2}}^X \to (\Omega C^{k-1})^X \to (\Omega C)^X.$$
    Now the fiber is even more truncated, so the same argument ensures we can whisker our cells with $\Omega C^{k-1} \to \Omega C$.

    To prove commutativity, extend the diagram to
    \[\begin{tikzcd}
    	{P^\wedge(\Omega  t_q^Y, t_s^W)} & {\Omega P^\wedge(t_q^Y, t_s^W)} & {\Omega t^Z_{q+s}} \\
    	{P^\wedge(t_{q-1}^{\Omega Y}, t_s^W)} & {t^{\Omega Z}_{q+s-1}} & {\operatorname{id}_{\Omega Z}} \\
    	& {\operatorname{id}_{\Omega Y \wedge W}} & {\operatorname{id}_{\Omega (Y \wedge W)}.}
    	\arrow["\oout", from=1-1, to=1-2]
    	\arrow["{P^\wedge ({E}^Y_{q}, \operatorname{id}_{t^W_s})}"', from=1-1, to=2-1]
    	\arrow[from=1-1, to=3-2]
    	\arrow[from=1-2, to=1-3]
    	\arrow[from=1-2, to=3-3]
    	\arrow["{{E}^Z_{q+s}}"{description}, from=1-3, to=2-2]
    	\arrow[from=1-3, to=2-3]
    	\arrow[from=2-1, to=2-2]
    	\arrow[from=2-1, to=3-2]
    	\arrow[from=2-2, to=2-3]
    	\arrow[from=3-2, to=2-3]
    	\arrow[from=3-2, to=3-3]
    	\arrow[from=3-3, to=2-3]
    \end{tikzcd}\]
    The commutativities of the three triangles come from the definitions of the arrows involved. Then the lower trapezoid commutes by the definition of $P^\wedge(t_{q-1}^{\Omega Y}, t_s^W) \to t^{\Omega Z}_{q+s-1}$. The parallelogram is seen to commute for formal reasons, and the two maps $\Omega P^\wedge(t_q^Y, t_s^W) \to \operatorname{id}_{\Omega Z}$ are equal by $\Omega$ of the definition of $P^\wedge(t_q^Y, t_s^W) \to t^Z_{q+s}$.

    This finishes the inductive step and thus the proof of the theorem.
\end{proof}

Thus, after the change of coordinates $(p, q) \mapsto (-p, N-q)$, we get an isomorphism between $E_1^{p,q}(Y^{\bullet}) \otimes E_1^{t,s}(Y^{\bullet}) \to E_1^{p+t,q+s}(Z^{\bullet})$ and $H^{p}(X , H^q(F_x, G)) \otimes H^{t}(X , H^s(F_x, G)) \to H^{p+t}(X , H^{q+s}(F_x, G))$ up to sign $(-1)^{q(N+s)+(N+m)(N+q) + (N+m)(N+n) + N+n} = (-1)^{q(N+s)+(N+t+s)p + N+p+q}$.

Now, let us discuss the convergence. We know that the pairing of abstract spectral sequences converges to the pairing $\pi_{N-p-q}(Y, \sgtn) \otimes \pi_{N-t-s}(Y, \sgtn) \xrightarrow{\circ} \pi_{2N-p-q-t-s}(Y \wedge Y, \sgtn) \to \pi_{2N-p-q-t-s}(Z, \sgtn)$. By Proposition \ref{pr:smashing}, this is
$(-1)^{N+p+q}$
times the absolute homotopy group pairing $\pi_{N-p-q}(Y) \otimes \pi_{N-t-s}(Y) \to \pi_{2N-p-q-t-s}(Z)$, and by Lemma \ref{le:pairing-gives-inner-prod}, the latter is
$(-1)^{(p+q)(N+t+s)}$
times the cohomology product pairing. We conclude that
the pairing of our spectral sequences converges to the product $H^{n}(\tot{x:X}{F_x}, G) \otimes H^{m}(\tot{x:X}{F_x}, G) \to H^{n+m}(\tot{x:X}{F_x}, G)$ up to sign $(-1)^{N+n(N+m+1)}$.

Let us redefine the abstract spectral sequence product by twisting it by the coefficient $(-1)^{N+(n-N)+(n-N)N}$.
Then the product converges to
$(-1)^{mn}$
times the cohomology product, and the product on the starting page is
$(-1)^{q(N+s)+(N+t+s)p + N+p+q+N+p+q+(p+q)N} = (-1)^{qs+pt+ps} = (-1)^{nm+qt}$
times the cohomology product.

The graded-commutativity coefficient for the spectral sequence then becomes
\\
$(-1)^{(N-n)(N-m)+N-n+N-m+N+n+nN+N+m+mN+N^2} = (-1)^{nm}$ ($N^2$ here is from $K_NG \wedge K_NG \to K_{2N}G$).

The Leibniz formula becomes $(-1)^{N+n+nN}d(xy) = (-1)^{N+n+1+(n+1)N+1}dxy + (-1)^{N+n+nN+(N-n)}xdy$, which is equivalently
$$(-1)^N d(xy) = dx \cdot y + (-1)^{n}xdy.$$

Note now that if $N$ is even, we get the standard signs for these formulas. For the input and the output of the spectral sequence, we still have an undesirable summand $nm$ in the exponent.

We now make the final correction: redefine the identifications of $E_1$ and of the abutment with the cohomology groups, i.e. the isomorphisms $E_1^{p,q} \simeq H^p(X, \lambda x.H^q(F_x, G))$ and $\pi_{N-n}(Y) \simeq H^n(\tot{x:X}{F_x}, G)$, by twisting them by the sign $(-1)^{e(n)}$, where $n = p+q$ on $E_1^{p,q}$ and $e(n) :\equiv \frac{n(n-1)}{2}$. Note that $e(n) + e(m) + e(n+m) \equiv_2 mn$, so this rids us from the undesirable summand $mn$.

We have thus settled the signs:
\begin{thm}\label{th:main-result}
    Given a type $X$, a type family $F_{-}: X \to U$ and a ring $R$, there exists a multiplicative spectral sequence $E_\bullet^{\bullet,\bullet}$ such that:
    \begin{enumerate}
        \item $E_\bullet^{\bullet,\bullet}$ naturally depends on $X$ and $F_{-}$;
        \item $E_1^{p,q} \simeq H^{p}(X, \lambda x.H^{q}(F_x, R))$;
        \item the diagram
        \[\begin{tikzcd}
        	{E_1^{p,q} \times E_1^{t,s} } & {E_1^{p+t,q+s}} \\
        	{H^{p}(X, \lambda x.H^{q}(F_x, R)) \times H^{t}(X, \lambda x.H^{s}(F_x, R))} & {H^{p+t}(X, \lambda x.H^{q+s}(F_x, R))}
        	\arrow[from=1-1, to=1-2]
        	\arrow["\simeq"', from=1-1, to=2-1]
        	\arrow["\simeq", from=1-2, to=2-2]
        	\arrow[from=2-1, to=2-2]
        \end{tikzcd}\]
        commutes up to sign $(-1)^{qt}$;
        \item $E_\bullet^{\bullet,\bullet}$ converges multiplicatively to the ring $H^{\bullet}(\tot{x:X}{F_x}, R)$;
        \item the pairing $-\cdot-: E_r^{p,q} \times E_r^{t,s} \to E_r^{p+t,q+s}$ is distributive and satisfies
        \begin{enumerate}
            \item
            the Leibniz formula
            $$d(x \cdot y) = (dx)\cdot y + (-1)^{p+q}x \cdot dy;$$
            \item
            the graded-commutativity relation
            $$xy = (-1)^{(p+q)(t+s)}yx.$$
            \end{enumerate}

    \end{enumerate}
\end{thm}
\begin{re}
    Can we get rid of the sign $(-1)^{qt}$ in (3) without spoiling the properties?
    The graded cup-product commutes up to sign $(-1)^{pt+qs}$, while the spectral sequence product commutes up to sign $(-1)^{(p+q)(t+s)}$. Thus, if we are to connect them up to sign $(-1)^{f(p,q,t,s,N)},$ the function $f$ has to satisfy $f(p,q,t,s) + f(t,s,p,q) = qt+ps$. $f(p,q,t,s) =  qt$ is the simplest function with this property.
\end{re}
\begin{proof}
    We want to take the colimit of the spectral sequences from above over all even $N$, for which we need to assemble them into a diagram. We will use the definition of the spectral sequence $E_\bullet^{\bullet, \bullet}$ associated to a tower $W^\bullet$ as arising from an exact couple $(E, D) = (\bigoplus \pi_{p+q}(F_q^W), \bigoplus \pi_{p+q}(W^q))$. (Our definition is equivalent, but we also change the indexing by $(p,q) \mapsto (-p, N-q)$.) Let $W^\bullet :\equiv \depfun{x:X}{(F_x \to K_NR)^\bullet}$ and $\bar{W}^\bullet  :\equiv \depfun{x:X}{(F_x \to K_{N-1}R)^\bullet}$, with colimits $W$ and $\bar{W}$. The map $\sigma_{N-1}$ gives rise to an equivalence $\bar{W} \to \Omega W$. To construct a map between consecutive spectral sequences $E^{\bullet, \bullet}_{\bullet}(\bar{W}^{\bullet})$ and $E^{\bullet, \bullet}_{\bullet}({W}^{\bullet})$, we construct a map of the corresponding exact couples.

    Let $p+q < N$. We define a map between the groups $E$ of the exact couples as the inverse of the sequence of natural isomorphisms
    \begin{align*}
        E^{p,q}(W^\bullet) =
        \pi_{N-p-q-1}(F^W_{N-q}) \equiv\pi_{N-p-q-2}(\Omega F^W_{N-q}) \to \\
        \pi_{N-p-q-2}(F^{\Omega W}_{N-q-1}) \to \\
        \pi_{N-p-q-2}(F^{\bar{W}}_{N-q-1}) =
        \pi_{(N-1)-p-q-1}(F^{\bar{W}}_{(N-1)-q})
        = E^{p,q}(\bar{W}^\bullet).
    \end{align*}
    For $D$, there is similarly an isomorphism
    \begin{align*}
        D^{p,q}(W^\bullet) = \pi_{N-p-q}(W^{N-q}) \equiv \pi_{N-p-q-1}(\Omega W^{N-q}) \to \\
        \pi_{N-p-q-1}((\Omega W)^{N-q-1}) \to \\
        \pi_{N-p-q-1}(\bar{W}^{N-q-1}) =
        \pi_{(N-1)-p-q}(\bar{W}^{(N-1)-q}) = D^{p,q}(\bar{W}^\bullet).
    \end{align*}
    When $p+q \geq N$, the domains $E^{p,q}(\bar{W}^\bullet), D^{p,q}(\bar{W}^\bullet)$ are trivial, and so will be the maps for the colimit diagram we are constructing.

    First of all, we need to check that this is a morphism of exact couples. This follows from the naturality of the long exact sequence of homotopy groups. (Recall that the maps $i,j,k$ in the exact couple construction come from the long exact sequences associated to pairs $(W^q, W^{q+1})$. Our maps of groups $E$ and $D$ are obtained from the maps of pairs $(\Omega W^{N-q}, \Omega W^{N-q+1}) \to ( {(\Omega W)}^{N-q-1},  {(\Omega W)}^{N-q}) \to ({\bar{W}}^{N-q-1}, {\bar{W}}^{N-q})$.)

    This proves the existence of the limiting spectral sequence. With what we already proved, it is straightforward to check the announced properties. It converges to the cohomology groups because from the definitions it follows that for $p+q<N$ the triangle
    \[
    \begin{tikzcd}
    	{\pi_{N-1-p-q}(\bar{W})} && {\pi_{N-p-q}(W)} \\
    	& {H^{p+q}(\tot{x:X}{F_x}, R)}
    	\arrow[from=1-1, to=1-3]
    	\arrow[from=1-1, to=2-2]
    	\arrow[from=1-3, to=2-2]
    \end{tikzcd}
    \]
    commutes.

    Next, we'd like to show that the triangles
    \[\begin{tikzcd}
    	{E^{p,q}(\bar{W})} && {E^{p,q}({W})} \\
    	& {H^{p}(X, \lambda x.H^q(F_x, R))}
    	\arrow[from=1-1, to=1-3]
    	\arrow[from=1-1, to=2-2]
    	\arrow[from=1-3, to=2-2]
    \end{tikzcd}\]
    commute for all $p,q$. This boils down to diagrams
    \[\begin{tikzcd}
    	{\Omega F_q^W} & {F_{q-1}^{\Omega W}} \\
    	{\Omega K_{q-1} \pi_{q}W} & {K_{q-2}\pi_{q-1}\Omega W}
    	\arrow["{e_q^W}", from=1-1, to=1-2]
    	\arrow["{\Omega f^W_q}"', from=1-1, to=2-1]
    	\arrow["{f^{\Omega W}_{q-1}}", from=1-2, to=2-2]
    	\arrow["{\sigma^{-1}_{q-2}}"', from=2-1, to=2-2]
    \end{tikzcd}\]
    which turn out to \textit{anti}commute by Lemma \ref{le:twoangles}. This is not a problem for us as we are restricting to only even $N$, going two arrows at a time.

    Next, we show that the products in these spectral sequences are compatible. Let $Z^{\bullet} : \equiv \depfun{x:X}{(F_x \to K_{2N}R)^{\bullet}}$.
    \\
    Since the product on every $E_r^{\bullet, \bullet}$ is defined in terms of its action on representatives of classes of elements from $\pi_\bullet({W^{\bullet}, W^{\bullet+1}})$, it is enough to analyze the commutativity of the diagram
    \[\begin{tikzcd}
    	{\Omega^n F_q^W \times \Omega^m F_s^W} & {\Omega^{n+m+1}(F^W_q \star F_s^W)} & {\Omega^{n+m+1}F^Z_{q+s}} \\
    	{\Omega^{n-1} F^{\bar{W}}_{q-1} \times \Omega^{m-1} {F}^{\bar{W}}_{s-1}} & {\Omega^{n+m-1}({F}^{\bar{W}}_{q-1} \star {F}_{s-1}^{\bar{W}})} & {\Omega^{n+m-1}{F}^{\bar{Z}}_{q+s-2}}
    	\arrow[from=1-1, to=1-2]
    	\arrow[from=1-1, to=2-1]
    	\arrow[from=1-2, to=1-3]
    	\arrow[from=1-3, to=2-3]
    	\arrow[from=2-1, to=2-2]
    	\arrow[from=2-2, to=2-3]
    \end{tikzcd}\]
    for all $n,m,q,s$.

    From the proof of Theorem \ref{th:products-coincide} we know that the upper composite can be rewritten as
    \[\begin{tikzcd}
    	{K_{q-n-1}\pi_qW \times K_{s-m-1}\pi_sW} & {K_{q+s-m-n-2}(\pi_qW \otimes \pi_sW)} & {K_{q+s-m-n-2}\pi_{q+s}Z} \\
    	{\Omega^n K_{q-1}\pi_qW \times \Omega^m K_{s-1}\pi_sW} && {\Omega^{n+m+1}K_{q+s-1}\pi_{q+s}Z} \\
    	{\Omega^n F_q^W \times \Omega^m F_s^W} & {\Omega^{n+m+1}(F^W_q \star F_s^W)} & {\Omega^{n+m+1}F^Z_{q+s}}
    	\arrow[from=1-1, to=1-2]
    	\arrow[from=1-2, to=1-3]
    	\arrow[from=1-3, to=2-3]
    	\arrow[from=2-1, to=1-1]
    	\arrow[from=2-3, to=3-3]
    	\arrow[from=3-1, to=2-1]
    	\arrow[from=3-1, to=3-2]
    	\arrow[from=3-2, to=3-3]
    \end{tikzcd}\]
    up to sign $(-1)^{qt}$, where $t = m-s, q =n-p$.

    The lower composite can be rewritten similarly, up to sign $(-1)^{(q-1)t}$. Now we have to check two things.
    The first is to check the commutativity of
    \[\begin{tikzcd}
    	{\Omega^{n-1} K_{q-2}\pi_{q-1}{\bar{W}} \times \Omega^{m-1} K_{s-2}\pi_{s-1}{\bar{W}}} & {\Omega^{n+m-1}K_{q+s-3}\pi_{q+s-2}\bar{Z}} \\
    	{\Omega^n K_{q-1}\pi_qW \times \Omega^m K_{s-1}\pi_sW} & {\Omega^{n+m+1}K_{q+s-1}\pi_{q+s}Z.}
    	\arrow[from=1-1, to=1-2]
    	\arrow[from=2-1, to=1-1]
    	\arrow[from=2-1, to=2-2]
    	\arrow[from=2-2, to=1-2]
    \end{tikzcd}\]
    Here $\bar{Z}^{\bullet }  :\equiv \depfun{x:X}{(F_x \to K_{2N-2}R)^\bullet}$ and the horizontal maps are from the composites as above. This follows from the naturality of all of the arrows forming the horizontal composites.

    The second thing is that the diagram
    \[\begin{tikzcd}
    	{\Omega^n K_{q-1}\pi_qW } & {\Omega^{n-1} K_{q-2}\pi_{q-1}{\bar{W}}} \\
    	{\Omega^n F_q^W } & {\Omega^{n-1} F^{\bar{W}}_{q-1}}
    	\arrow[from=1-1, to=1-2]
    	\arrow[from=2-1, to=1-1]
    	\arrow[from=2-1, to=2-2]
    	\arrow[from=2-2, to=1-2]
    \end{tikzcd}\]
    (as well as the two similar diagrams) commutes. To use the definition of the vertical maps, we apply $\Omega^{q-n-1}$ to the diagram and postcompose with the maps to $\pi_q W$ and $\pi_{q-1}\bar{W}$ to reduce to the diagram
    \[\begin{tikzcd}
    	{\pi_q W} & {\pi_{q-1}{\bar{W}}} \\
    	{\pi_qW^q} & {\pi_{q-1}\bar{W}^{q-1}} \\
    	{\Omega^{q-1}{F_q^W}} & {\Omega^{q-2} F^{\bar{W}}_{q-1}}
    	\arrow[from=1-1, to=1-2]
    	\arrow[from=2-1, to=1-1]
    	\arrow[from=2-1, to=2-2]
    	\arrow[from=2-2, to=1-2]
    	\arrow[from=3-1, to=2-1]
    	\arrow[from=3-1, to=3-2]
    	\arrow[from=3-2, to=2-2]
    \end{tikzcd}\]
    whereupon the commutativity is clear.
    Thus, after two steps along the ``$N$-axis'', we get commutativity up to sign $(-1)^{t+t} = 1$ (our change of coordinates does not change the parity of $t$). The Leibniz formula follows at once.

    Similarly, from naturality it follows that the convergences to the product in the target are compatible since there is no sign involved in the isomorphism.

    Finally, there are two ways to go about the graded-commutativity. The first one is to deduce it from the graded-commutativity of the spectral sequence and the commutativity of $W \wedge W \to Z$, as described before the statement of the theorem (see also Remark \ref{re:graded-comm-SS}). A simpler way seems to be to note that we already have it on $E_1^{\bullet, \bullet}$ due to the isomorphism with cohomology. Then, it follows for the following pages as by construction the product respects the fact that $E_{r+1}$ is a subquotient of $E_r$.
    \end{proof}

\section{Applications}\label{se:appls}

As in the classical case, multiplicativity enables simple calculations for many cohomology rings. We enumerate some examples here. See \cite{FF}, Ch. 3, 4 and \cite{Hat5} for many examples of basic computations with the Serre spectral sequence.

\begin{itemize}
\item
    $H^\bullet(\CC P^n, R)$.

    There is a fiber sequence $S^1 \to S^{2n+1} \to \CC P^n$ (\cite{BR}, Section 6). We know (\cite{Bru}, 5.3.1) that $H^k(S^n, \ZZ)$ is $\ZZ$ for $k = 0, n$ and $0$ otherwise. Taking coefficients in a ring $R$, we get a spectral sequence with two nonzero rows, both isomorphic to $H^\bullet(\CC P^n, R).$ Since $H^1(S^{2n+1}, R) = 0$, all groups $E_{\infty}^{p,q}$ on the diagonal $p+q = 1$ must be trivial. There are no nonzero differentials hitting $E_{\bullet}^{1,0}$, so $E_1^{1,0} = H^1(\CC P^n, R)$ must be the trivial group. Next, for the copy of $H^1(S^1, R)$ sitting at $(0,1)$, which is a group generated by some element $x$, to be killed, the differential $d_1^{0,1}: E_1^{0,1} \to E_1^{2,0}$ must be an isomorphism. We conclude that $E_1^{2,0} = H^2(\CC P^n, R) = R$. Let $y$ be the generator of this group.

    We remark that $E_1^{2,1}$ is generated by the element $x \cdot y$. Indeed, let $X : \equiv \CC P^n, R_\bullet :\equiv H^\bullet(S^1, R)$. The product of $x$ and $y$ in $H^\bullet(X , H^{\bullet}(S^1, R))$ is represented by the map $X \xrightarrow{(x,y)} K_0 R_1 \wedge K_2 R_0 \xrightarrow{} K_2 (R_1 \otimes R_0) \to K_2 R_1$. The product $R_1 \otimes R_0 \to R_1$ defining the last map is by definition equal to the cup-product in $H^\bullet(S^1, R)$. Knowing this cohomology ring, we conclude that $R_1 \otimes R_0 \to R_1$ equals $R \otimes R \xrightarrow{\cdot} R$.
    Thus, the composite above is equal to the ordinary product in $H^\bullet(X, R)$ where we multiply the 2-dimensional generator by the 0-dimensional generator. We will implicitly use such arguments further.

    Continuing inductively, we get that $H^{2k+1}(\CC P^n, R) = 0$ and $d_1^{2k,1}: x \cdot y^k \mapsto y^{k+1} $ is an isomorphism for $k < n$ by the Leibniz formula. This stops at $k=n$, when the domain of $d_1^{2k,1}$ doesn't have to get killed. Instead, we have $R = H^{2n+1}(S^{2n+1}, R) = \ker d_1^{2k,1}$. Thus, if $R$ has no proper subgroups isomorphic to $R$ (e.g. if $R = \ZZ/k$), we conclude:
    $$H^\bullet(\CC P^n, R) \simeq R[x]/x^{n+1}, |x| = 2.$$

    For the same conclusion for a general ring $R$, note that finite-dimensional CW-complexes have finite cohomological dimension as their cohomology groups can be calculated from the cellular chain complex, as in \cite{BuchFav}. Taking Example 8 from \cite{MimOle} with $K_2 \ZZ$ in place of $B\ZZ_2$, we conclude that $\CC P^n$ is such a CW-complex. From the spectral sequence it follows that $H^{2n+2k}(\CC P^n, R) = \im d_1^{2n,1}$ for all $k>0$, and hence $d_1^{2n,1} = 0$.
\item
    $H^\bullet(K_2\mathbb{Z}, R)$.

    This is the same argument as the first half of the previous one. There is a fiber sequence $K_1\mathbb{Z} \to \sgtn \to K_2\mathbb{Z},$ giving rise to a Serre spectral sequence converging to $H^\bullet(\sgtn, R)$. Hence, the differential $d_1^{0,1}: E_{1}^{0, 1} \to E_{1}^{2,0}$ must be an isomorphism. Then $H^2(K_2\mathbb{Z}, R) \simeq E_{1}^{2,0} \simeq  R.$ Let $x$ generate $E_{1}^{0, 1}$ and $y$ generate $E_{1}^{2,0}$. We further calculate $d_1(x \cdot y) = d_1(x) \cdot y = y^2;$ $d_1(x \cdot y^k) = d_1(x) \cdot y^k = y^{k+1},$ so we conclude that $H^\bullet(K_2\mathbb{Z}, R) \simeq R[a]$ with $|a|=2.$
\item
    $H^\bullet(\Omega S^n, \mathbb{Z})$.

    There is a fiber sequence $\Omega S^n \to \sgtn \to S^n,$ giving rise to a Serre spectral sequence converging to $H^\bullet(\sgtn, \ZZ)$. Hence, the differential $E_{n-1}^{0, n-1} \to E_{n-1}^{n,0}$ must be an isomorphism, and so $E_1^{0, n-1} \simeq \ZZ$. Let $a_1$ be its generator and let $x : \equiv d_{n-1}^{0, n-1}(a_1)$. We similarly get, for all $k > 1$, that $E_1^{0, kn-k} \simeq \ZZ = \langle a_k \rangle$ (to kill the group located at $(n, (k-1)(n-1))$, which is generated by $a_{k-1} x$).
    \\
    When $n$ is odd, $d_{n-1}(a_1^2) = a_1x + (-1)^{n-1} a_1x = 2 a_1 x$, and so $a_1^2 =  2a_2$. Inductively, we get $a_1^k = k! a_{k}$. We conclude that $H^\bullet(\Omega S^n, \mathbb{Z}) \simeq \Gamma[x], |x| = n-1$, where $\Gamma[x]$ is the \textit{divided polynomial algebra}, defined as the subring of $\QQ[x]$ additively generated by the monomials $x^i/i!$.
    \\
    When $n$ is even, a simple algebraic argument along the same lines shows (see \cite{Hat5}, Example 5.17) that $H^\bullet(\Omega S^n, \mathbb{Z}) \simeq \Lambda [a] \otimes \Gamma[b], |a| = n-1, |b| = 2n - 2$ (where $\Lambda[x]$ is the exterior algebra on one generator).
\item
    $H^\bullet(({S^3})^4, \ZZ)$.

    Let $X$ be the 3-connected cover of $S^3$.
    There is a fiber sequence $K_2\ZZ \to X \to S^3$. Since the term $(0,2)$ does not survive in $E_\infty$, the differential $d_2^{0,2}$ is an isomorphism. It then follows from the Leibniz formula that for all positive $k$, $d_2^{0,2k}$ is multiplication by $k$. (Note that this is the first time when our sign $(-1)^{qt}$ enters. However, it doesn't have any importance.) Thus,
    $$H^k(X, \ZZ) = \begin{cases}
        \ZZ/m\ZZ, & k = 2m + 1,\ m = 2, 3, 4, \ldots \\
        0, & k \neq 5, 7, 9, \ldots
    \end{cases}.$$

    Note that if we had the universal coefficient theorem, or at least a nontrivial mapping $H^5(X) \to H_4(X),$ this would compute the Brunerie number.

\item
    $H^\bullet(K_n\ZZ, \mathbb{Q})$.

    This ring is isomorphic to $H^\bullet(S^n, \mathbb{Q})$ for $n$ odd and $\mathbb{Q}[x], |x|=n$ for $n$ even. This is already known to us for $n=1, 2$ (as $K_1 \ZZ$ is a circle and $K_2 \ZZ$ was analyzed above). We proceed by induction on $n$, using the fiber sequence $K_{n-1}\ZZ \to \sgtn \to K_{n}\ZZ$. If $n$ is even, the argument is the same as for $K_2 \ZZ$ above. If $n$ is odd, we have $\mathbb{Q}[x]$ in columns $0$ and $n$ and trivial groups in between.
    Let the isomorphism $d_{n-1}: E_{n-1}^{0,n-1} \to E_{n-1}^{n, 0}$ send the generator $x$ to the generator $y$. We get that $dx^2=(1+(-1)^{n-1})xdx = 2xdx = 2xy$ generates $E_{n-1}^{n,n-1}$, and likewise $dx^{k+1}$ generates $E_{n-1}^{n,k(n-1)}$. Thus, columns $0$ and $n$ annihilate each other on the $(n-1)$-th page, except for the $(0,0)$-th term. It follows that all differentials reaching $E_\bullet^{p,0}$ for $p \neq 0, n$ are zero, implying our claim.

\end{itemize}

\subsection{Future work}
Theorem \ref{th:main-result} opens more doors to the synthetic development of classical results. In work in progress, we study cohomology operations from a more computational point of view. This requires some form of the K\"unneth formula. Consider a trivial fiber sequence $X \times Y \to X$. Taking the coefficients to be a finite field $k$, one gets a spectral sequence whose differentials can be shown to all vanish. We then get $H^n(X \times Y, k) \simeq \bigoplus_i H^i(X, H^{n-i}(Y, k))$. Now, multiplicativity ensures that this isomorphism is compatible with the product structure.

\section{Appendix: technicalities}\label{se:lemmas}

\subsection{Category theory}

\begin{de}\label{de:nat-tr}
    Define a (wild) \textit{natural transformation} between wild functors $F, G : \mathcal{C} \to \mathcal{D}$ as a function $\eta: \depfun{c:\mathrm{ob}\mathcal{C}}{Fc \to Gc}$ together with a witness $\operatorname{nat}_\eta: \depfun{c,c':\mathrm{ob}\mathcal{C}}{\depfun{f: c \to c'}{Gf \circ \eta c = \eta c' \circ Ff}}$.
\end{de}

Let $f, g: c \to c'$ be two arrows in $\mathcal{C}$, $h: f=g$ be a homotopy, and $\eta$ be a natural transformation from $F$ to $G$. There is an equality of type $h_* (\operatorname{nat}_\eta(f)) = \operatorname{nat}_\eta(g)$.
This type is the type of fillers for the prism
\[\begin{tikzcd}
    	Fc & {Fc'} \\
    	Gc & {Gc'}
    	\arrow["Fg", curve={height=6pt}, from=1-1, to=1-2]
    	\arrow["Ff", curve={height=-6pt}, from=1-1, to=1-2]
    	\arrow[from=1-1, to=2-1]
    	\arrow[from=1-2, to=2-2]
    	\arrow["Gg", curve={height=6pt}, from=2-1, to=2-2]
    	\arrow["Gf", curve={height=-6pt}, from=2-1, to=2-2]
\end{tikzcd}\]
whose square sides are naturalities of $\eta$. In Lemma \ref{le:merid-smash-Sm} we face the case where $(\mathcal{C}, \mathcal{D}, F, G, \eta) \equiv (U_*, U_*, \Sigma (- \wedge B), \Sigma (-) \wedge B, \siin)$ and $f, g$ are constant maps. We will implicitly use there that $\Sigma \sgtn = \sgtn$ and $\sgtn \wedge B = \sgtn$. Transported along these equalities, the square faces of the prism disappear and then its filler $\ap_{\operatorname{nat}_\eta(-)}(h)$ says that $\siin^{-1}$ takes $\Sigma$ out from homotopies:
$$\siin^{-1} \whisk (\Sigma \alpha \wedge B) \whisk \siin =_* \Sigma (\alpha \wedge B).$$

\begin{de}\label{de:coherent-nat-tr}
    We call a natural transformation $\eta: F \implies G$ \textit{coherent} if it respects compositions; namely, if for any $a \xrightarrow{f} b \xrightarrow{g} c$, there is a filler in
    \[\begin{tikzcd}
    	Fa & Fb \\
    	Ga & Gb & Fc \\
    	&& {Gc.}
    	\arrow[from=1-1, to=1-2]
    	\arrow[from=1-1, to=2-1]
    	\arrow[from=1-1, to=2-3]
    	\arrow[from=1-2, to=2-2]
    	\arrow[from=1-2, to=2-3]
    	\arrow[from=2-1, to=2-2]
    	\arrow[from=2-1, to=3-3]
    	\arrow[from=2-2, to=3-3]
    	\arrow[from=2-3, to=3-3]
    \end{tikzcd}\]
\end{de}

\begin{lem}\label{le:nat-tr-cube}
Given a coherent natural transformation $\eta: F \implies G$ and a commutative square
\begin{tikzcd}
	a & b \\
	c & {d}
	\arrow[from=1-1, to=1-2]
	\arrow[from=1-1, to=2-1]
	\arrow["s"{description}, Rightarrow, from=1-2, to=2-1]
	\arrow[from=1-2, to=2-2]
	\arrow[from=2-1, to=2-2]
\end{tikzcd}, the cube
\[\begin{tikzcd}
	& Fa \\
	Fc & Ga & Fb \\
	Gc & Fd & Gb \\
	& Gd
	\arrow[from=1-2, to=2-1]
	\arrow[from=1-2, to=2-2]
	\arrow[from=1-2, to=2-3]
	\arrow[from=2-1, to=3-1]
	\arrow[from=2-1, to=3-2]
	\arrow["{\ }"{description}, from=2-2, to=3-1]
	\arrow["{\ }"{description}, from=2-2, to=3-3]
	\arrow[from=2-3, to=3-2]
	\arrow[from=2-3, to=3-3]
	\arrow[from=3-1, to=4-2]
	\arrow[from=3-2, to=4-2]
	\arrow[from=3-3, to=4-2]
\end{tikzcd}\]
has a filler.
\end{lem}
\begin{proof}
    Decompose this cube into three prisms along the following ``cuts'' of its top:
    \[\begin{tikzcd}
    	& Fa \\
    	Fc && Fb \\
    	& {Fd.}
    	\arrow[from=1-2, to=2-1]
    	\arrow[from=1-2, to=2-3]
    	\arrow[curve={height=-12pt}, from=1-2, to=3-2]
    	\arrow[curve={height=12pt}, from=1-2, to=3-2]
    	\arrow[from=2-1, to=3-2]
    	\arrow[from=2-3, to=3-2]
    \end{tikzcd}\]
    The middle prism is the naturality of the witness of the naturality of $\eta$ with respect to $s$. The two side prisms are proofs that $\eta$ is coherent.
\end{proof}

\begin{de}\label{de:pentagon}
    A wild category is said to satisfy the \textit{pentagon axiom} if for any composable arrows $X_1 \xrightarrow{d} X_2 \xrightarrow{c} X_3 \xrightarrow{b} X_4 \xrightarrow{a} X_5$, composing associators can only produce a trivial element of $\Omega [X_1 \to X_5, a(b(cd))]$. That is, the diagram
    \[\begin{tikzcd}
    	& {(ab)(cd)} & {((ab)c)d} \\
    	{a(b(cd))} & {a((bc)d)} & {(a(bc))d}
    	\arrow[from=1-2, to=1-3]
    	\arrow[from=1-3, to=2-3]
    	\arrow[from=2-1, to=1-2]
    	\arrow[from=2-1, to=2-2]
    	\arrow[from=2-2, to=2-3]
    \end{tikzcd}\]
    commutes.
\end{de}

\begin{re}\label{re:pentagon}
    One needs the pentagon identity in a wild category to hold to be able to paste 3-fillers. Let us illustrate how it arises. Consider first pasting two squares:
    \[\begin{tikzcd}
    	\bullet & \bullet & \bullet \\
    	\bullet & \bullet & {\bullet.}
    	\arrow["x", from=1-1, to=1-2]
    	\arrow["w"', from=1-1, to=2-1]
    	\arrow["y", from=1-2, to=1-3]
    	\arrow[from=1-2, to=2-2]
    	\arrow["z", from=1-3, to=2-3]
    	\arrow["k"', from=2-1, to=2-2]
    	\arrow["t"', from=2-2, to=2-3]
    \end{tikzcd}\]
    This has to produce an element of type $z(yx) = (tk)w$. To use the right square to construct such a path, one has to first use the associator $z(yx) = (zy)x$, since the right square is a path starting from $zy$.

    Going one dimension higher, imagine pasting two cubes:
    \[\begin{tikzcd}
    	\bullet &&&&& \bullet \\
    	& \bullet &&& \bullet \\
    	&& \bullet & \bullet \\
    	&& \bullet & \bullet \\
    	& \bullet &&& \bullet \\
    	\bullet &&&&& {\bullet.}
    	\arrow[from=1-1, to=1-6]
    	\arrow["a", from=1-1, to=2-2]
    	\arrow[from=1-1, to=6-1]
    	\arrow[from=1-6, to=2-5]
    	\arrow[from=1-6, to=6-6]
    	\arrow[from=2-2, to=2-5]
    	\arrow["b", from=2-2, to=3-3]
    	\arrow[from=2-2, to=5-2]
    	\arrow[from=2-5, to=3-4]
    	\arrow[from=2-5, to=5-5]
    	\arrow["c", from=3-3, to=3-4]
    	\arrow[from=3-3, to=4-3]
    	\arrow["d", from=3-4, to=4-4]
    	\arrow[from=4-3, to=4-4]
    	\arrow[from=5-2, to=4-3]
    	\arrow[from=5-2, to=5-5]
    	\arrow[from=5-5, to=4-4]
    	\arrow[from=6-1, to=5-2]
    	\arrow[from=6-1, to=6-6]
    	\arrow[from=6-6, to=5-5]
    \end{tikzcd}\]
    Each cube separately tells us that the two 2-cells identifying two paths between its opposite vertices are equal. Now we try to construct an equality between two 2-dimensional paths that start at, say, $d \circ (c \circ (b\circ a))$. The path ``$dcba$'' then appears in 5 arising 2-cells (3-cells for the initial diagram), the edges of which form a pentagon:
\[\begin{tikzcd}
	&& {} \\
	{} & {\ap_{\operatorname{assoc}}(B)} & {((dc)b)a} & {f \whisk a} & {} \\
	{F :\equiv \> ?} & {(dc)(ba)} && {(d(cb))a} & {\ap_{\operatorname{assoc}}(t)} \\
	& {d(c(ba))} && {d((cb)a)} \\
	{} && {\mathrm{def}(T)} && {}
	\arrow["{(B \whisk b) \whisk a}"{description}, from=2-3, to=1-3]
	\arrow["{B \whisk (ba)}"{description}, squiggly, from=3-2, to=2-1]
	\arrow[squiggly, from=3-2, to=2-3]
	\arrow[squiggly, from=3-4, to=2-3]
	\arrow["{(d \whisk t) \whisk a}"{description}, squiggly, from=3-4, to=2-5]
	\arrow[squiggly, from=4-2, to=3-2]
	\arrow[squiggly, from=4-2, to=4-4]
	\arrow["{ d \whisk T}"{description}, squiggly, from=4-2, to=5-1]
	\arrow[squiggly, from=4-4, to=3-4]
	\arrow["{d \whisk (t \whisk a)}"{description}, from=4-4, to=5-5]
\end{tikzcd}\]
    Here $T$ denotes the top 2-cell of the composite cube, $t$ the top 2-cell of the ``back'' (smaller) cube, and $B$ denotes the back 2-cell. We try to define the left cell $F$ of this diagram, using, in particular, the filler of the back cube $f$. The cell $F$ should equate two paths starting at $d(c(ba))$: one first applies $T$ and then continues along irrelevant sides of the cube, and the other first applies $B$ (for which it first uses an associator), and then continues along irrelevant sides of the cube. The bottom cell of the diagram is the definition of $T$ as a composite of the two smaller top 2-cells. We don't know anything about $B \whisk (ba)$, but the 3-filler $f$ of the back cube knows about $B \whisk b$. Hence we pass to this filler via an associator, which is the top two cells of the diagram. Finally, we note that another associator allows us to close the loop via the rightmost cell.
\end{re}

\subsection{Properties of -in and -out functions}\label{se:-in-out}
\begin{lem}\label{le:sout-w}
    For any pointed types $A, B: U_*$, the map
    $\siin: \Sigma (A \wedge B) \pto \Sigma A \wedge B$ is an equivalence.
\end{lem}
\begin{proof}
    An explicit inverse $\sout: \Sigma A \wedge B \pto \Sigma (A \wedge B) $ can be defined by
    \begin{flalign*}
        & \langle N, b \rangle, \auxr, \auxl \mapsto N \\
        & \langle S, b \rangle \mapsto S \\
        & \langle \merid_a, b \rangle \mapsto \merid_{\langle a,b \rangle} \\
        & \gluel: N \mapsto \refl_N \\
        & \gluel: S \mapsto \merid_{\bs{}}^{-1} \\
        & \gluel: \merid_a \mapsto byContr \\
        & \gluer: b \mapsto \refl_N.
    \end{flalign*}
\end{proof}

\begin{lem}\label{le:sout-j}
    For any pointed types $A, B: U_*$, the map
    $\siin: \Sigma (A \join B) \pto \Sigma A \join B$ is an equivalence.
\end{lem}
\begin{proof}
    An explicit inverse $\sout: \Sigma A \join B \pto \Sigma (A \join B) $ can be defined by
    \begin{flalign*}
        & \inl N \mapsto N \\
        & \inl S \mapsto S \\
        & \ap_{\inl}(\merid_a) \mapsto \merid_{\inl a} \\
        & \inr b \mapsto N \\
        & \glue(N, b) \mapsto \refl_N \\
        & \glue(S, b) \mapsto \merid_{\bs{A}}^{-1}  \\
        & \ap_{\glue}(\merid_a) \mapsto \ap_{\merid^{-1}}(\glue(a, b) \cdot \glue(\bs{A}, b)^{-1}),
    \end{flalign*}
    where for the last line we implicitly use that transport along $\merid_a$ precomposes with $\merid^{-1}_{\inl(a)}$.
\end{proof}

\begin{lem}\label{le:siin-nat-w}
    For any $B:U_*$, $\siin$ is a natural transformation from $\Sigma (- \wedge B)$ to $(\Sigma -) \wedge B$.
\end{lem}
\begin{proof}
    Straightforward.
\end{proof}

\begin{lem}\label{le:siin-nat-j}
    $\siin: \Sigma (A \join B) \to \Sigma A \join B$ is natural in $A$ and $B$: given arrows $f: A \to A', g: B \to B'$ we have $\siin \circ \Sigma (f \join g) = (\Sigma f \join g) \circ \siin$.
    \\
    The same is true in $U_*$.
\end{lem}
\begin{proof}
    The proof is straightforward. Equate the two functions on $N$ by $\refl_{\inl(N)}$ and on $S$ by $\refl_{\inl(S)}$. Now we need to equate the images of $\merid_x$ for $x: A \join B$. When $x \equiv \inl a$, both functions compute on $\merid_x$ as $\ap_{\inl}(\merid_{f(a)})$. For $x \equiv \inr b$, both functions take $\merid_x$ to $\glue(N, g(b)) \cdot \glue(S, g(b))^{-1}$. Finally, when $x$ moves along $\glue(a,b)$, we get two fillers for
    \[\begin{tikzcd}
    	{\inl(N)} & {\inr(gb)} \\
    	{\inl(S).}
    	\arrow["{\glue(N, gb)}", squiggly, from=1-1, to=1-2]
    	\arrow["{\ap_{\inl}(\merid_{f(a)})}"', curve={height=12pt}, squiggly, from=1-1, to=2-1]
    	\arrow["{\glue(S, gb)^{-1}}", squiggly, from=1-2, to=2-1]
    \end{tikzcd}\]
    One is $\cone_{g(b)}(N, S, \merid_{f(a)})$ and the other is $\ap_{\ap_{\Sigma f \join g}}(\cone_b(N, S, \merid_a))$. Replacing the former with $\cone_{g(b)}(ha, ha', \ap_hp)$ and the latter with $\ap_{\ap_{h \join g}}(\cone_b(a, a', p))$ for arbitrary $a, a': \Sigma A, p: a=a', h: \Sigma A \to \Sigma A'$, we get the equality by path induction.
    For the pointed case, note that both functions are pointed by $\refl_{\inl(N)}$, and our homotopy equals $\refl_{\inl(N)}$ on $N$.
\end{proof}

\begin{lem}\label{le:oout-nat-j}
    The join $\oout$ is a natural transformation: let $f: A \pto A'$ be a pointed map and $B$ be a type. Then $\oout \circ (\Omega f \join B) = \Omega (f \join B) \circ \oout$.
    \\
    The same is true in $U_*$.
\end{lem}
\begin{proof}
    As always, by path induction we can assume that $f$ is pointed by $\refl$. The 0-dimensional constructors are easy, e.g. send $\inl(p)$ to the obvious path
    $$\ap_{\inl} \Omega f(p)
    \rightsquigarrow \ap_{\inl} \ap_f(p)
    \rightsquigarrow \ap_{f \join B} \ap_{\inl}(p)
    \rightsquigarrow \Omega{(f \join B)} \oout(\inl p).$$

    For $\glue(p, b)$, we have to fill the diagram
    \[\begin{tikzcd}
    	{\ap_{\inl} \ap_f(p)} & {\ap_{f \join B} \ap_{\inl}(p)} \\
    	{\refl_{\bs{}}} & {\ap_{f \join B}(\refl_{\bs{}}).}
    	\arrow[squiggly, from=1-1, to=1-2]
    	\arrow["{c(\ap_fp, b)}"', squiggly, from=1-1, to=2-1]
    	\arrow["{\ap_{\ap_{f \join B}}(c(p,b))}", squiggly, from=1-2, to=2-2]
    	\arrow[squiggly, from=2-1, to=2-2]
    \end{tikzcd}\]

    Replace the loop $p$ with an arbitrary path $p: a=_Aa'$. Then this diagram generalizes to a diagram in $\inl(fa)= \inl (fa')$
    \[\begin{tikzcd}
    	{\ap_{\inl} \ap_f(p)} & {\ap_{f \join B} \ap_{\inl}(p)} \\
    	{\glue(fa, b) \cdot \glue(fa', b)^{-1}} & {\ap_{f \join B}(\glue(a, b) \cdot \glue(a', b)^{-1}).}
    	\arrow[squiggly, from=1-1, to=1-2]
    	\arrow["{\cone_b(fa, fa', \ap_fp)}"', squiggly, from=1-1, to=2-1]
    	\arrow["{\ap_{\ap_{f \join B}}(\cone_b(a, a', p))}", squiggly, from=1-2, to=2-2]
    	\arrow[squiggly, from=2-1, to=2-2]
    \end{tikzcd}\]
    A filler now exists by path induction.
\end{proof}

\begin{lem}\label{le:oout-nat-w}
    The smash $\oout$ is natural in $A$ and $B$: let $f: A \pto A', g: B \pto B'$ be pointed maps. Then $\oout \circ (\Omega f \wedge g) =_* \Omega (f \wedge g) \circ \oout$.
\end{lem}
\begin{proof}
    Straightforward; the proof only involves $byContr$ and obvious path algebra.
\end{proof}
By induction, we then get
\begin{lem}\label{le:ooutk-nat}
    $\ooutk{k}: \Omega^k A \wedge B \to \Omega^k(A \wedge B)$ is natural in $A$, $B$. The same is true for $\ooutk{k}_2$.
\end{lem}

\begin{lem}\label{le:i-nat}
    The map $I_{n,m}^{A,B}$ from Theorem \ref{th:defs-of-pairing} is natural in $A$ and $B$.
\end{lem}
\begin{proof}
    The base case $n=m=0$ follows immediately from the definition. The inductive step follows
    from the naturality of $\oout$ (Lemma \ref{le:oout-nat-j}).
\end{proof}

\begin{lem}\label{le:oout-anticom-w}
    The smash $\oout$ and $\oout_2$ anticommute: $\Omega \oout \circ \oout_2 =_{\Omega X \wedge \Omega Y \to \Omega^2(X \wedge Y)} -\Omega\oout_2 \circ \oout$.
\end{lem}
\begin{proof}
We will assume $X \wedge Y$ is pointed by $\langle \bs{X}, \bs{Y}\rangle$.
We have $\ap_{\oout_2}(\lps{p}{q}) = \ap_{\oout_2}\ap_{\ctos_{X, \Omega Y}}((p,q)) = \ap_{\lambda x.\rpsl{x}{q}}(p)$, which can be seen as the square $[p,q]^{-1}$
    \[\begin{tikzcd}
    	{\bs{}} & {\bs{}} \\
    	{\bs{}} & {\bs{}}
    	\arrow["{\lps{p}{\bs{Y}}}", squiggly, from=1-1, to=1-2]
    	\arrow[""{name=0, anchor=center, inner sep=0}, "{\rps{\bs{X}}{q}}", squiggly, from=2-1, to=1-1]
    	\arrow["{\lps{p}{\bs{Y}}}"', squiggly, from=2-1, to=2-2]
    	\arrow[""{name=1, anchor=center, inner sep=0}, "{\rps{\bs{X}}{q}}"', squiggly, from=2-2, to=1-2]
    	\arrow[between={0.2}{0.8}, Rightarrow, from=0, to=1]
    \end{tikzcd}\]
    where the left and the right sides are contracted by $byContr$.
Then $\ap_{\oout_2}(\oout(\langle p,q\rangle))$ is the same square whose top and bottom sides are also contracted, thus the inverse of $[p,q]^l$. Similarly, $\ap_{\oout}(\oout_2(\langle p,q\rangle))$ equals $[p,q]^l$, which provides the equality on inputs of the form $\langle p, q\rangle$. It is straightforward to see that it respects path constructors, as $\oout$ just acts on them using cancellations and contractibility.
\end{proof}

\begin{cor}\label{cor:oout-anticom-w}
    For the smash $\oout$, we have
    $$\Omega^q \ooutk{p} \circ \ooutk{q}_2 =_* (-1)^{pq}\Omega^p \ooutk{q}_2 \circ \ooutk{p}.$$
\end{cor}
\begin{proof}
    The cases when $p=0$ or $q=0$ are obvious, and the case $p=q= 1$ is Lemma \ref{le:oout-anticom-w}.

    When $pq >1$, we induct on $p$:
    \begin{align*}
        \Omega^q \ooutk{p} \circ \ooutk{q}_2 \equiv \Omega^q (\Omega^{p-1}\oout \circ \ooutk{p-1}) \circ \ooutk{q}_2
        \\ =_* \Omega^{q+p-1}\oout \circ \Omega^q \ooutk{p-1} \circ \ooutk{q}_2
        \\ =_* (-1)^{(p-1)q}\Omega^{q+p-1}\oout \circ \Omega^{p-1} \ooutk{q}_2 \circ \ooutk{p-1}
        \\ =_* (-1)^{(p-1)q}\Omega^{p-1} (\Omega^q \oout \circ \ooutk{q}_2) \circ \ooutk{p-1}
        \\ =_* (-1)^{pq}\Omega^{p-1} (\Omega \ooutk{q}_2 \circ \oout) \circ \ooutk{p-1}
        \\ =_* (-1)^{pq}\Omega^{p} \ooutk{q}_2 \circ \Omega^{p-1} \oout \circ \ooutk{p-1}
        \\ \equiv (-1)^{pq}\Omega^{p} \ooutk{q}_2 \circ \ooutk{p}.
    \end{align*}
\end{proof}

The same result can be repeated for the join $\oout$. The proof will be different only for the base case:
\begin{lem}\label{le:oout-anticom-j}
    The join $\oout$ and $\oout_2$ anticommute: for $A, B : U_*$, we have $\Omega \oout_2 \circ \oout =_{\Omega A \join \Omega B \to \Omega^2 (A \join B)} -\Omega \oout \circ \oout_2$.
\end{lem}
\begin{proof}
    We will leave out the details as we never use join anticommutativity except for reality check.

    Denote $w(a,a',b) : \equiv \glue(a, b) \glue(a',b)^{-1},$ $w'(a,b,b') : \equiv  \glue(a,b)^{-1} \glue(a,b')$.

    Define by path induction the type $\tetrahedron = \tetrahedron (p: a=_Aa', q: b=_Bb')$ as the tetrahedron obtained by joining $p$ and $q$. Technically, it is a 2-cell in the type $w(a,a',b) = w(a,a',b')$ between the path $\mathcolor{blue}{\cone_b(a,a',p)^{-1} \cone_{b'}(a, a', p)}$ and the horizontal composition (see Remark \ref{re:ap-composite}) of cells $\mathcolor[HTML]{007700}{\glue(a,b) \ap_{\inr}(q) \rightsquigarrow \glue(a,b')}$ and $\mathcolor[HTML]{007700}{(\glue(a',b) \ap_{\inr}(q))^{-1} \rightsquigarrow \glue(a',b')^{-1}}$ that come from $\cone'_a(b,b',q)$ and $\cone'_{a'}(b,b',q)$:
    \[\begin{tikzcd}
    	& {\inr(b)} \\
    	& {\inr(b')} && {\inl(a')} \\
    	{\inl(a)}
    	\arrow[squiggly, from=1-2, to=2-2]
    	\arrow[squiggly, from=2-4, to=1-2]
    	\arrow[""{name=0, anchor=center, inner sep=0}, squiggly, from=2-4, to=2-2]
    	\arrow[squiggly, from=3-1, to=1-2]
    	\arrow[""{name=1, anchor=center, inner sep=0}, squiggly, from=3-1, to=2-2]
    	\arrow[""{name=2, anchor=center, inner sep=0}, squiggly, from=3-1, to=2-4]
    	\arrow[draw={rgb,255:red,92;green,214;blue,92}, between={0.2}{1}, Rightarrow, from=1-2, to=0]
    	\arrow[draw={rgb,255:red,92;green,214;blue,92}, between={0.2}{1}, Rightarrow, from=1-2, to=1]
    	\arrow[draw={rgb,255:red,92;green,92;blue,214}, between={0.2}{1}, Rightarrow, from=1-2, to=2]
    	\arrow[draw={rgb,255:red,92;green,92;blue,214}, between={0.2}{1}, Rightarrow, from=2, to=2-2]
    \end{tikzcd}\]

    Let $a:A$. Generalize the map $\oout: \Omega A \join \Omega B \to \Omega (A \join \Omega B)$ to $({\bs{A} = a) \join \Omega B}$ using the same formula: $\inl p \mapsto \ap_{\inl}(p), \inr q \mapsto w(\bs{A}, a, q), \glue(p,q) \mapsto \cone_q(\bs{A},a,p)$. Now the codomain is the type $\inl \bs{A} = \inl a$; note that $\ap_{\oout_2}$ still takes this to $\Omega^2{(A \join B)}$ since both endpoints are sent to $\refl$. The composition then sends $\inl (a,p)$ to $\refl_{\refl}$, $\inr q$ to $c'(a, q) \cdot c'(a', q)^{-1}: \refl \rightsquigarrow \inr q \rightsquigarrow \refl$, and $\glue(a,p,q)$ to $t :\equiv \ap_{\ap_{\oout_2}} (\cone_q(\bs{A}, a, p))$. Note that it makes sense to compare the latter 3-cell $t$ with $\tetrahedron(p,q)$ after the following identifications. First, $w(a,a',b) = w(a,a',b')$ equals $\refl_{\bs{}} = \refl_{\bs{}}$ when $q$ is a loop. Then, transported along this equality, the endpoints of $\tetrahedron(p,q)$ can be connected to those of $t$ using path algebra. To equate $t$ and $\tetrahedron(p,q)$ now, it is enough to note that when $p = \refl_a$, they both are equal to $cancel$.

    Similarly, the second composite $\Omega \oout \circ \oout_2$ equals the map that sends $\inl p$ and $\inr q$ to $\refl_{\refl}$ and $\glue(p,q)$ to a whisker of $\tetrahedron'(p,q)^{-1}$ with cancellation, where $\tetrahedron'$ is a 2-path in the type $w'(a,b,b') = w'(a',b,b')$ defined symmetrically to $\tetrahedron$. Note the inversion that comes from the fact that now the endpoint of $\tetrahedron'$ (the green path that now is $\mathcolor[HTML]{007700}{\cone'_a(b,b',q)^{-1} \cone'_{a'}(b,b',q)}$) is now the starting point of the image of $\glue(p,q)$. We now use our standard argument: by path induction $\tetrahedron(p,q)$ and $\tetrahedron'(p,q)$ are equal up to path algebra that makes them of the same type. This path algebra commutes with contractions, and so the images of $\glue(p,q)$ under the two maps are equal.
\end{proof}

\begin{lem}\label{le:oout-mirrors-oout2-j}
    For $X, Y: U_*$, the diagram
    \[\begin{tikzcd}
    	{\Omega^p X \join Y} & {\Omega^p (X \join Y)} \\
    	{ Y \join \Omega^p X} & {\Omega^p (Y \join X)}
    	\arrow["{\ooutk{p}}", from=1-1, to=1-2]
    	\arrow["\mir"', from=1-1, to=2-1]
    	\arrow["{\Omega^p\mir}", from=1-2, to=2-2]
    	\arrow["{\ooutk{p}_2}"', from=2-1, to=2-2]
    \end{tikzcd}\]
    commutes in $U_*$.
\end{lem}
\begin{proof}
    The case $p=1$ follows from the definitions. Namely, we need to check that for any $p: \Omega X, y: Y$, there is a filler for
    \[\begin{tikzcd}
    	{\glue(\bs{Y}, \bs{X})\ap_{\inr}(p)\glue(\bs{Y}, \bs{X})^{-1}} \\
    	{\glue(\bs{Y}, \bs{X})\glue(\bs{Y}, \bs{X})^{-1}} \\
    	{\refl_{\inl(\bs{Y})},}
    	\arrow["{\ap_{\mathrm{\Omega \mir}}(c(p, y))}", curve={height=-12pt}, squiggly, from=1-1, to=2-1]
    	\arrow[curve={height=12pt}, squiggly, from=1-1, to=3-1]
    	\arrow[squiggly, from=2-1, to=3-1]
    \end{tikzcd}\]
    where the left arrow is the left-pointing cone $c_2(y,p)$. This generalizes to the diagram
    where on the right is $\glue(\bs{Y}, x) \whisk\ap_{\mathrm{\ap_{\mir}}}(\cone_{\bs{Y}}(x, x' ,p)) \whisk \glue(\bs{Y}, x')^{-1}$ composed with cancellations, and on the left is $\cone'_y(x,x',p)$ used to define $c_2$ (defined in the obvious way by path induction), whiskered with $\glue(\bs{Y}, x)^{-1}$ on the right and $\glue(\bs{Y}, x')$ on the left. Now, path induction is applicable.

    The case $p>1$ now follows by induction.
    \end{proof}

\begin{lem}\label{le:oout-mirrors-oout2-w}
    For $X, Y: U_*$, the diagram
    \[\begin{tikzcd}
    	{\Omega^p X \wedge Y} & {\Omega^p (X \wedge Y)} \\
    	{ Y \wedge \Omega^p X} & {\Omega^p (Y \wedge X)}
    	\arrow["{\ooutk{p}}", from=1-1, to=1-2]
    	\arrow["\swap"', from=1-1, to=2-1]
    	\arrow["{\Omega^p\swap}", from=1-2, to=2-2]
    	\arrow["{\ooutk{p}_2}"', from=2-1, to=2-2]
    \end{tikzcd}\]
    commutes in $U_*$.
\end{lem}
\begin{proof}
    The case $p=0$ is obvious and the case $p=1$ is immediate from the definitions. The case $p>1$ follows by induction.
\end{proof}

The transformations $\siin$ and $\oout$ are related via the adjunction $\Sigma \dashv \Omega$. Recall that the counit $\epsilon_A: \Sigma \Omega A \pto A$ sends $N,S$ to $\bs{A}$ and $\merid_p$ to $p$.
\begin{lem}\label{le:siin-oout-adj}
    For pointed spaces $X, Y$, the diagram
    \[\begin{tikzcd}
    	{\Sigma (\Omega X \join Y)} & {\Sigma \Omega (X \join Y)} \\
    	{\Sigma \Omega X \join Y} & {X \join Y}
    	\arrow["{\Sigma \oout}", from=1-1, to=1-2]
    	\arrow["{\siin }"', from=1-1, to=2-1]
    	\arrow["{\epsilon_{X \join Y}}", from=1-2, to=2-2]
    	\arrow["{\epsilon_X \join Y}", from=2-1, to=2-2]
    \end{tikzcd}\]
    commutes (in $U_*$).
\end{lem}
\begin{proof}
    Computing from the definitions, the first (top-right) composite sends $\ap_{\merid}(\glue(p,y))$ to the diagram
    \[\begin{tikzcd}
    	{\bs{}} \\
    	{\bs{},}
    	\arrow["{\ap_{\inl}(p)}"', curve={height=6pt}, squiggly, from=1-1, to=2-1]
    	\arrow["\refl", curve={height=-6pt}, squiggly, from=1-1, to=2-1]
    \end{tikzcd}\]
    where the 2-cell equals $\ap_{\ap_{\epsilon_{X \join Y}}}(c(p, y))$.
    The other composite sends it to
    \[\begin{tikzcd}
    	{\bs{}} \\
    	& {\inr(y)} \\
    	{\bs{},}
    	\arrow["{\glue(\bs{X}, y)}", squiggly, from=1-1, to=2-2]
    	\arrow["{\ap_{\inl}(p)}"', squiggly, from=1-1, to=3-1]
    	\arrow["{\glue^{-1}(\bs{X}, y)}", squiggly, from=2-2, to=3-1]
    \end{tikzcd}\]
    where the 2-cell equals $\cone_y(\bs{}, \bs{}, p)$.

    Now let us replace in these diagrams $p$ with an arbitrary $p:\bs{X}=_X x$. The coherence diagram for the second composite still makes sense (with $x$ in place of one of the lower $\bs{X}$); for the first one, instead of $c(p,y)$ we take $\cone_{y}(\bs{X}, x, p): \ap_{\inl}(p) =_{\inl(\bs{X}) = \inl(x)} \glue(\bs{X}, y) \glue(x, y)^{-1}$, and instead of $\epsilon_{X \join Y}$ -- the map
    $$\depfun{x:X}{\depfun{p:\bs{X} = x}{\Sigma (\inl(\bs{X}) = \inl x) \to X}} \join Y$$
    given by
    $N \mapsto \inl\bs{X}, S \mapsto \inl x, \merid_p \mapsto p$.
    Now the two diagrams are equal by path induction. Specializing to the case $x = \bs{X}$, we get the desired equality.
\end{proof}

\begin{lem}\label{le:siin-anticom}
    The following diagram in $U_*$ anticommutes:
    \[\begin{tikzcd}
    	{\Sigma \Sigma(A \wedge B)} & {\Sigma ( \Sigma A \wedge B)} \\
    	{\Sigma (A \wedge \Sigma B)} & {\Sigma A \wedge \Sigma B.}
    	\arrow["{\Sigma \siin}"', from=1-1, to=1-2]
    	\arrow["{\Sigma \siin_2}", from=1-1, to=2-1]
    	\arrow["{\siin_2}"', from=1-2, to=2-2]
    	\arrow["\siin", from=2-1, to=2-2]
    \end{tikzcd}\]
\end{lem}
We will restrict to a manageable number of glueings and prove this only for the case $A \equiv S^0$. This is the only case we will use. When $A \equiv S^0$, the map $B \pto A \wedge B$, $b \mapsto \langle 1_\s ,b \rangle,$ is an equivalence. Thus, we prove the lemma after precomposing with $\Sigma \Sigma B \pto \Sigma \Sigma (A \wedge B)$.
\begin{proof}

 Our $\Sigma \siin$ and $\Sigma \siin_2$ are assumed to be precomposed with $\Sigma \Sigma B \pto \Sigma \Sigma ( S^0 \wedge B)$, so let us write explicitly how the composites $\Sigma B \pto S^1 \wedge B$ and $\Sigma B \pto S^0 \wedge \Sigma B$ work. The first one acts as
$$N \mapsto \auxr, S \mapsto \auxr, \merid_b \mapsto [\mmerid_{{1_\s}}, b\rangle$$
and the second one as
$$N \mapsto \auxr, S \mapsto \auxr, \merid_b \mapsto \langle {1_\s}, \mmerid_b].$$
Both are pointed by reflexivity.

Denote the equality we are seeking by $\alpha: \siin_2 \circ \Sigma \siin =_{\Sigma \Sigma B \pto S^1 \wedge \Sigma B} \siin \circ (-\Sigma \siin_2)$.
For clarity, we add superscripts $0$ and $1$ to the outer and the inner meridians of $\Sigma \Sigma B$. For any $b:B$, there is a diagram in $\Sigma \Sigma B$ of the form

\[\begin{tikzcd}
	N \\
	{S,}
	\arrow[""{name=0, anchor=center, inner sep=0}, "{\merid^0_N}"', curve={height=12pt}, squiggly, from=1-1, to=2-1]
	\arrow[""{name=1, anchor=center, inner sep=0}, "{\merid^0_S}", curve={height=-12pt}, squiggly, from=1-1, to=2-1]
	\arrow[shorten <=5pt, shorten >=5pt, Rightarrow, from=0, to=1]
\end{tikzcd}\]
where the 2-cell is $\ap_{\merid_{(-)}^0}(\merid_b^1)$. We need to specify how $\alpha$ acts on each such diagram. First, set $\alpha(N) \equiv \alpha (S) :\equiv \refl_{\auxr}$.
Now we need to provide squares $\ap_{\alpha}(\merid^0_x)$ for $x \equiv N,S$. This is easy, because from the definitions both functions we equate compute to $\refl_{\auxr}$ on these meridians, e.g.
$$\ap_{\siin \circ \Sigma \siin_2}(\merid^0_N) = \ap_{\siin}(\merid^1_{\auxr}) = \refl_{\auxr}.$$

Next, fix $b:B$. We need the 3-filler $\ap_{\ap_{\alpha}(\merid^0_{(-)})}(\merid^1_b)$. We compute and compare the other two sides of the prism. We have
\begin{align*}
& \ap_{\ap_{\siin_2 \circ \Sigma \siin} \circ \merid^0_{(-)}} (\merid^1_b) = \\
& \ap_{\ap_{\siin_2} \circ (\ap_{\Sigma \siin} \circ \merid^0_{(-)})} (\merid^1_b)  = \\
& \ap_{\ap_{\siin_2} \circ {\merid^1_{(-)}  \circ \siin}} (\merid^1_b)  = \\
& \ap_{\ap_{\siin_2} \circ {\merid^1_{(-)}}} ([\mmerid_{1_{\s}}, b \rangle ) = \\
& \ap_{\ap_{\siin_2} \circ {\merid^1_{(-)}}} (\gluer_b^{-1} \cdot {\ap_{\ctos}}((\mmerid_{1_{\s}}, b)) \cdot \gluer_b).
\end{align*}
The factor without $\gluer_b$ computes as
\begin{align*}
& \ap_{\ap_{\siin_2} \circ {\merid^1_{(-)}}} ({\ap_{\ctos}}((\mmerid_{1_{\s}}, b ))) = \\
& \ap_{\ap_{\siin_2} \circ {\merid^1_{(-)}} \circ {\ctos}}((\mmerid_{1_{\s}}, b )) = \\
& \ap_{\lambda(s,b).\langle s, \mmerid_b]}((\mmerid_{1_{\s}}, b )).
\end{align*}
We conclude that
\begin{align*}
& \ap_{\ap_{\siin_2 \circ \Sigma \siin} \circ \merid^0_{(-)}} (\merid^1_b) = \\
& byContr \cdot \ap_{\ap_{\ctos}}(\ap_{(-, \mmerid_b)}(\mmerid_{1_{\s}})) \cdot byContr,
\end{align*}
where the two factors $byContr$ should be seen as contracting the sides of this square.

Similarly, $\ap_{\ap_{\siin \circ (-\Sigma \siin_2)} \circ \merid^0_{(-)}} (\merid^1_b)$ equals $
\ap_{\ap_{\ctos}}(\ap_{(\mmerid_{1_{\s}}^{-1}, -)}(\mmerid_b))$ up to paths $byContr$ that make the latter ``of the right type''.

The essence of our argument is that the cells $\ap_{\ap_{\ctos}}(\ap_{(-, \mmerid_b)}(\mmerid_{1_{\s}}))$ and $\ap_{\ap_{\ctos}}(\ap_{(\mmerid_{1_{\s}}^{-1}, -)}(\mmerid_b))$ are the same (equal up to what naturally makes them of the same type). Assume given types $X, Y$ and paths $p: x = x'$, $q: y = y'$. Then, intuitively, $\ap_{(-, q)}(p)$ and $\ap_{(p, -)}(q)$ represent the same square in $X \times Y$; the difference is that in the former case we sweep it horizontally, and in the latter case vertically. These two directions differ by a symmetry about the diagonal (which has degree $-1$). Such a symmetry is homotopic to the symmetry about the vertical axis; thus, to actually represent the same square, instead of $\ap_{(p, -)}(q)$ we can take $\ap_{(p^{-1}, -)}(q)$. We just need to make this simple idea ``well-typed''.

Consider $\ap_{(-, q)}(p): p_*((x,q)) = (x',q)$. The latter type is equal to $(x, q) \cdot (p, y') = (p, y) \cdot (x', q)$.
Similarly, there is an element $\ap_{(p^{-1}, -)}(q): q_*((p^{-1},y)) = (p^{-1},y')$ with the latter type equivalent to $(p^{-1},y) \cdot (x, q) = (x', q) \cdot (p^{-1},y'),$ or $(x, q) \cdot (p, y') = (p, y) \cdot (x', q)$. Thus, we get two elements of the same type. By path induction, they are equal. Importantly, we used only whiskers with sides, associators and cancellation to equate the two types.

Returning to our cells in $S^1 \wedge \Sigma B,$ we need to prove that $\ap_{\ap_{\ctos}}(\ap_{(-, \mmerid_b)}(\mmerid_{1_\s}))$ and $\ap_{\ap_{\ctos}}(\ap_{(\mmerid_{1_{\s}}^{-1}, -)}(\mmerid_b))$ are equal after composing with cells that contract the sides of these squares. By the argument above, they become equal after transporting along appropriate whiskering with their sides and canceling. To conclude, it only remains to note that by $byContr$, contracting the boundaries of the squares is the same as whiskering with some sides, canceling, and only then contracting the boundaries (compare with the base case proof of Corollary \ref{cor:distr1}).

It is easy to see that $\alpha$ is based, since by definition, $\alpha (N)$, as well as the pointedness proofs for the functions $\alpha$ equates, equals $\refl_{\auxr}$.
\end{proof}

Note that for any $A, B: U_*$ and for every $C$, there is a sequence of natural maps
    \begin{align*}
        &(A \wedge B \pto C) \xrightarrow{\sim} \\
        &(B \pto A \pto C) \to \\
        &(B \pto \Omega A \pto \Omega C) \xrightarrow{\sim} \\
        &(\Omega A \wedge B \pto \Omega C).
    \end{align*}
\begin{lem}\label{le:oout-desription}
    $\oout: \Omega A \wedge B \pto \Omega (A \wedge B)$ can be described as the image under the composition above
    of the universal map $\id: A \wedge B \pto A \wedge B$.
    \\
    $\oout_2$ can be described similarly.
\end{lem}
\begin{proof}
    It is straightforward to see that the image of $\id$ acts as $\oout$.
    \\
    The first arrow takes $\id$ to $\lambda b.\langle -, b \rangle$.
    The second map then takes this to $\lambda b. \lambda p. \mathrm{ptd}(\langle -, b \rangle)^{-1} \cdot \lps{p}{b} \cdot \mathrm{ptd}(\langle -, b \rangle) = [-, - \rangle$.
    The third arrow then takes this to the map that acts on 0-dimensional constructors as $\auxr, \auxl \mapsto \bs{}; \langle p, b \rangle \mapsto [p, b \rangle$ and on 1-dimensional constructors as
    \begin{align*}
    & \gluer_{b} \mapsto \mathrm{ptd}([-, b\rangle) = byContr \cdot cancel; \\
    & \gluel_{p} \mapsto \mathrm{ptd}([p, - \rangle) = byContr \cdot cancel.
    \end{align*}
\end{proof}

\begin{re}\label{re:oout-desription}
    It follows that for any $A, B, C : U_*$ the diagram
    \[\begin{tikzcd}
    	{A \wedge B \pto C} & {(\Omega (A \wedge B) \pto \Omega C)} & {((\Omega A \wedge B) \pto \Omega C)} \\
    	{(B \pto (A \pto C))} & {(B \pto (\Omega A \pto \Omega C))}
    	\arrow[from=1-1, to=1-2]
    	\arrow[from=1-1, to=2-1]
    	\arrow[from=1-2, to=1-3]
    	\arrow[from=1-3, to=2-2]
    	\arrow[from=2-1, to=2-2]
    \end{tikzcd}\]
commutes.
\end{re}

\begin{lem}\label{le:siin-compatible}
    Given pointed types $A, B,$ the diagram
    \[\begin{tikzcd}
    	{\Sigma (A \join B)} && {\Sigma A \join B} \\
    	{\Sigma \Sigma(A \wedge B)} & {\Sigma \Sigma(A \wedge B)} & {\Sigma(\Sigma A \wedge B)}
    	\arrow["\siin", from=1-1, to=1-3]
    	\arrow["{\Sigma\jtos}", from=1-1, to=2-1]
    	\arrow["\jtos", from=1-3, to=2-3]
    	\arrow["{-1}", from=2-1, to=2-2]
    	\arrow["{\Sigma \siin}", from=2-2, to=2-3]
    \end{tikzcd}\]
    commutes in $U_*$.
\\
 The same holds for $\siin_2$.
\end{lem}

\begin{proof}
     We construct a homotopy $\alpha: \jtos \circ \siin = \Sigma \siin \circ (-1) \circ \Sigma \jtos$.
     We are faced, for every $a: A, b: B$, with the prism
\[\begin{tikzcd}
	N & S \\
	N & {N.}
	\arrow[squiggly, from=1-1, to=1-2]
	\arrow["\refl"{description}, curve={height=-12pt}, squiggly, from=1-1, to=2-1]
	\arrow["{m_{\langle N, b \rangle} \cdot m^{-1}_{\langle S, b \rangle}}"', curve={height=12pt}, squiggly, from=1-1, to=2-1]
	\arrow["{m^{-1}_{\bs{}}}", curve={height=12pt}, squiggly, from=1-2, to=2-2]
	\arrow["{m^{-1}_{\bs{}}}", curve={height=-12pt}, squiggly, from=1-2, to=2-2]
	\arrow[squiggly, from=2-1, to=2-2]
\end{tikzcd}\]
    Here $m$ denotes the meridian of $\Sigma (\Sigma A \wedge B)$. The top and bottom edges are the values of our homotopy on $N$ and $S$, respectively. Set the former to be $\merid_{\bs{}}$ and the latter to be $\refl$.
    One square face corresponds to $\ap_\alpha(\merid_{\inr(b)})$. To provide it, we need to equate $\merid_{\langle N, b \rangle} \cdot \merid^{-1}_{\langle S, b \rangle}$ with $\merid_{\bs{}}^{-1}$. Equate the first pair of factors $\merid_{\langle N, b \rangle}$ and $\merid_{\bs{}}$ by $\ap_{\merid}(byContr)$, and the second pair of factors by $\ap_{\merid^{-1}_{\langle - , b \rangle}}(\merid^{-1}_{\bs{A}}) \cdot \ap_{\merid^{-1}}(byContr)$.
    Next, the other square face corresponds to $\ap_\alpha(\merid_{\inl(a)})$.
    To provide it, we need to equate $\refl_N$ with $\merid_{\bs{}}^{-1}$, which we do in the obvious way.
    \\
    The side faces of the prism correspond to $\ap_{\ap_{f}(\merid_{(-)})}(\glue(a,b))$, where $f$ is one of the two functions we are equating. We defined the action of $\siin$ on $\ap_{\merid} \glue(a,b)$ as $\cone_b(N, S, \merid_a)$, but by path induction, we can think of it as coming from $\ap_{\glue(-, b)}(\merid_a^{-1})$. Then concretely, the left side equals $\ap_{\merid_{\langle - , b \rangle}}(\merid_a)$  and the right side equals
    \begin{align*}
    \ap_{\ap_{\Sigma \siin \circ (-1) \circ \Sigma \jtos } (\merid)}(\glue(a,b)) = \ap_{\ap_{\Sigma \siin } (\merid)^{-1} \circ \jtos}(\glue(a,b)) \\
    = \ap_{(\merid)^{-1}  \circ \siin \circ \jtos}(\glue(a,b)) = \ap_{(\merid)^{-1}  \circ \siin}(\merid_{\langle a,b \rangle }) \\
    = \ap_{\merid^{-1}}([\merid_a , b \rangle).
    \end{align*}

    To provide $\ap_{\ap_\alpha(\merid_{(-)})}(\glue(a,b))$, we need to fill the prism. It is not hard to spot that all the sides of the prism cancel out. We rearrange them somewhat to make it easier to see:
    \[\begin{tikzcd}
    	&&& S \\
    	&&& S \\
    	N & N & N & N & N & N & {N.}
    	\arrow["{m^{-1}_{\bs{}}}"{description}, squiggly, from=1-4, to=3-1]
    	\arrow["{m^{-1}_{\bs{}}}"{description}, squiggly, from=1-4, to=3-7]
    	\arrow["\refl"{description}, squiggly, from=2-4, to=1-4]
    	\arrow["{m^{-1}_{\langle N, b \rangle}}"{description}, squiggly, from=2-4, to=3-1]
    	\arrow["{m_{\langle S, b \rangle}^{-1}}"{description}, squiggly, from=2-4, to=3-2]
    	\arrow["{m^{-1}_{\langle N, b \rangle}}"{description}, squiggly, from=2-4, to=3-3]
    	\arrow["{m^{-1}_{\bs{}}}"{description}, squiggly, from=2-4, to=3-4]
    	\arrow["{m^{-1}_{\langle N, b \rangle}}"{description}, squiggly, from=2-4, to=3-5]
    	\arrow["{m^{-1}_{\langle S, b \rangle}}"{description}, squiggly, from=2-4, to=3-6]
    	\arrow["{m^{-1}_{\langle N, b \rangle}}"{description}, squiggly, from=2-4, to=3-7]
    	\arrow[squiggly, from=3-1, to=3-2]
    	\arrow[squiggly, from=3-2, to=3-3]
    	\arrow[squiggly, from=3-3, to=3-4]
    	\arrow[squiggly, from=3-4, to=3-5]
    	\arrow[squiggly, from=3-5, to=3-6]
    	\arrow[squiggly, from=3-6, to=3-7]
    \end{tikzcd}\]
    All self-equalities of $N$ are $\refl$ in this picture.
    \\
    Here we subdivide the edge $\merid_{\langle N, b \rangle} \cdot \merid^{-1}_{\langle S, b \rangle}$ of the prism into two and set its middle $S$ to be the center of the diagram. The subdivided edge is then depicted as two sides of the leftmost small triangle, which itself corresponds to the left side of the prism. It is equal, as an equality from $\merid_{\langle S, b \rangle}^{-1} =_{S = N} \merid_{\langle N, b \rangle}^{-1}$, to $\ap_{\merid^{-1}_{\langle -, b \rangle}}(\merid_a)^{-1}$. The front of the prism then breaks: its upper part is the leftmost (big) triangle, and the lower part is the second and the third small triangles. Next, the four triangles to the right of the middle $\merid_{\bs{}}$ form the right side of the prism. Finally, the back of the prism is just cancellation (for this edge orientation) and it corresponds to the outer triangle (with sides $\refl$ and $\merid_{\bs{}}$). The picture is symmetric, and symmetric triangles are equal, so everything cancels out.

    Finally, to see that $\alpha$ is based we need to check that $\merid_{\bs{}} = \refl \cdot (\merid_{\bs{}}^{-1})^{-1}$, which is clear.
\end{proof}

\begin{lem}\label{I00-oout}
For pointed types $A, B$ and $b: B$, the diagram
\[\begin{tikzcd}
	{\Omega A} && {\Omega \Omega (A \star B)} \\
	&& {\Omega (\Omega A \star B)}
	\arrow["{-\Omega I_{0,0}^{A,B}(-,b)}", from=1-1, to=1-3]
	\arrow["{I_{0,0}^{\Omega A, B}(-,b)}"', from=1-1, to=2-3]
	\arrow["{\Omega \oout}"', from=2-3, to=1-3]
\end{tikzcd}\]
commutes.
\end{lem}
\begin{proof}
    Take $p: \Omega A$ in the domain. The horizontal map sends it to
    $$ cancel^{-1} \cdot  \ap_{\lambda a.\glue(\bs{A}, \bs{B})\cdot \glue(a, \bs{B})^{-1} \cdot \glue(a, b) \cdot \glue(\bs{A}, b)^{-1}}(p^{-1}) \cdot cancel.$$

    The composite of the other two arrows sends $p$ to
    $$cancel^{-1} \cdot c(\refl, \bs{B}) \cdot  c(p, \bs{B})^{-1} \cdot c(p, b) \cdot c(\refl, b)^{-1} \cdot cancel.$$
    Note that $cancel^{-1} \cdot c(\refl, b) = \refl_{\refl}$ for any $b:B$. Thus, we need to compare $c(p, \bs{B})^{-1} \cdot c(p,b)$ with the 2-loop above.
    This is formal, given that both just express the geometric idea of moving along the two triangles in
    \[\begin{tikzcd}
    	{\inr(\bs{B})} && {\inl(\bs{A})} && {\inr(b)} \\
    	&& {\inl(\bs{A})}
    	\arrow["{\glue(\bs{A}, \bs{B})}"', squiggly, from=1-3, to=1-1]
    	\arrow["{\glue(\bs{A}, b)}", squiggly, from=1-3, to=1-5]
    	\arrow["{\glue(\bs{A}, \bs{B})}"{description}, squiggly, from=2-3, to=1-1]
    	\arrow["{\ap_{\inl}(p)}"{description}, squiggly, from=2-3, to=1-3]
    	\arrow["{\glue(\bs{A}, b)}"{description}, squiggly, from=2-3, to=1-5]
    \end{tikzcd}\]
    but in slightly different ways. The path $c(p, \bs{B})^{-1} \cdot c(p,b)$ first writes $\refl$ as two sides of the left triangle, then applies its filler, then moves along the right triangle, and then cancels back to $\refl$. The other path writes $\refl$ as the boundary of the outer triangle, and then applies two of its sides (see Remark \ref{re:ap-composite}) before canceling. As always, a rigorous way to say this consists in making one of the $\bs{A}$ an arbitrary $a:A$ and applying path induction.
\end{proof}

\subsection{Lemmas and proofs}
\begin{lem}\label{le:LES-isom}
    Let $i: A \pto B$ be a pointed pair with fiber $F$. Consider the exact sequence of types
    $$\Omega^{n+1} B \xrightarrow{d} \Omega^n F \to \Omega^n A \to \ldots \to F \to A \to B.$$
    Let $e: \Omega^{n+1} B$. Then in the diagram
    \[
    \begin{tikzcd}
    	{S^n} & \sgtn & A \\
    	\sgtn & {S^{n+1}} & B & F \\
    	&&& \sgtn
    	\arrow[from=1-1, to=1-2]
    	\arrow[from=1-1, to=2-1]
    	\arrow[dashed, from=1-1, to=2-4]
    	\arrow[from=1-2, to=1-3]
    	\arrow[from=1-2, to=2-2]
    	\arrow[from=1-3, to=2-3]
    	\arrow[from=2-1, to=2-2]
    	\arrow["\lrcorner"{anchor=center, pos=0.125, rotate=180}, draw=none, from=2-2, to=1-1]
    	\arrow["{\adj{n+1}e}", from=2-2, to=2-3]
    	\arrow[from=2-4, to=1-3]
    	\arrow["\lrcorner"{anchor=center, pos=0.125, rotate=-135}, draw=none, from=2-4, to=2-3]
    	\arrow[from=2-4, to=3-4]
    	\arrow[from=3-4, to=2-3]
    \end{tikzcd}
    \]
    the dashed map equals $\adj{n}(de)$.
\end{lem}
\begin{proof}
    We induct on $n$. The base case $n=0$ follows from explicit descriptions for the universal properties of the fiber and the pushout.

    Assume now $n>0$ and consider the pointed pair $\Omega i$ with fiber $\Omega F$. By induction, in the diagram
    \[
    \begin{tikzcd}
    	{S^{n-1}} & \sgtn & {\Omega A} \\
    	\sgtn & {S^{n}} & {\Omega B} & {\Omega F} \\
    	&&& \sgtn
    	\arrow[from=1-1, to=1-2]
    	\arrow[from=1-1, to=2-1]
    	\arrow[dashed, from=1-1, to=2-4]
    	\arrow[from=1-2, to=1-3]
    	\arrow[from=1-2, to=2-2]
    	\arrow[from=1-3, to=2-3]
    	\arrow[from=2-1, to=2-2]
    	\arrow["\lrcorner"{anchor=center, pos=0.125, rotate=180}, draw=none, from=2-2, to=1-1]
    	\arrow["{\adj{n}e}", from=2-2, to=2-3]
    	\arrow[from=2-4, to=1-3]
    	\arrow["\lrcorner"{anchor=center, pos=0.125, rotate=-135}, draw=none, from=2-4, to=2-3]
    	\arrow[from=2-4, to=3-4]
    	\arrow[from=3-4, to=2-3]
    \end{tikzcd}
    \]
    the dashed arrow is $-\adj{n-1}(de): S^{n-1} \to \Omega F$, since the $m$-th differential of the fiber sequence of $i$ is minus the $(m-1)$-th differential of the fiber sequence of $\Omega i$ (see \cite{UF}, Lemma 8.4.4). Applying $\adj{}$ to this diagram, we get a diagram as in the statement of the lemma, except the square
    \[
    \begin{tikzcd}
    	{\Sigma S^{n-1}} & {\Sigma \sgtn} \\
    	{\Sigma \sgtn} & {\Sigma S^n}
    	\arrow[from=1-1, to=1-2]
    	\arrow[from=1-1, to=2-1]
    	\arrow["{\Sigma \merid}"{description}, Rightarrow, from=1-2, to=2-1]
    	\arrow[from=1-2, to=2-2]
    	\arrow[from=2-1, to=2-2]
    \end{tikzcd}
    \]
    is the opposite of the pushout square by Lemma \ref{le:susp-of-merid}. This provides the second minus sign, giving us what we want.
\end{proof}

We prove Proposition \ref{pr:smashing}.
\begin{prp}\label{pr:smashing-proof}
    For $n,m \geq 1$, the map $\pi_n(X,\sgtn) \times \pi_m(Y,\sgtn) \to \pi_{n+m}(X \wedge Y,\sgtn)$
    under the isomorphism $\pi_\bullet(-, \sgtn) \simeq \pi_{\bullet}(-)$
    equals the truncation of
    $$\Omega^n X \times \Omega^m Y \to \Omega^n X \wedge \Omega^m Y \xrightarrow{\ooutk{m}_2} \Omega^m (\Omega^n X \wedge Y) \xrightarrow{\Omega^m \ooutk{n}} \Omega^{n+m} (X \wedge Y)$$
    up to sign $(-1)^n$.
\end{prp}
\begin{proof}
    We will assume that $X \wedge Y$ is pointed at $\langle \bs{X}, \bs{Y} \rangle$.

    We prove the proposition by induction on $n$ and $m$.
    For the base case $n=m=1$, we need to show that the diagram
    \[\begin{tikzcd}
    	{\Omega X \times \Omega Y} & {\Omega (\Omega X \join\Omega Y)} \\
    	{\Omega X \wedge \Omega Y} & {\Omega^2 (X \wedge Y)}
    	\arrow["{I_{0,0}^{\Omega X, \Omega Y}}", from=1-1, to=1-2]
    	\arrow[from=1-1, to=2-1]
    	\arrow["{\Omega i}", from=1-2, to=2-2]
    	\arrow[from=2-1, to=2-2]
    \end{tikzcd}
    \]
    anticommutes. Here the lower map is $\Omega \oout \circ \oout_2$, and the map $i$ is $\Omega X \join \Omega Y \simeq \fib(P^\times(\sgtn \to X, \sgtn \to Y)) \simeq \fib(X \vee Y \to X \times Y) \to \fib(\sgtn \to X \wedge Y) \simeq \Omega (X \wedge Y)$ (see Lemma \ref{le:fiber-of-pp}).
    Concretely, $i$ equals the map that sends $\inl(p)$ to $\lps{p}{\bs{Y}}$, $\inr(q)$ to $\rps{\bs{X}}{q}$, and $\glue(p,q)$ to the cell shaped as $t \cdot [p,q]^{-1} \cdot l$, where the cells $t$ and $l$ identify two of $[p,q]$'s sides with $\refl$; that is, we take the square $[p,q]^{-1}$, contract a copy of $\lps{p}{\bs{Y}}$ and $\rps{\bs{X}}{q}$ on its boundary, and regard this as a path between the two remaining sides of $[p, q]$:
    \[\begin{tikzcd}
    	{\langle \bs{X}, \bs{Y} \rangle } & {\langle \bs{X}, \bs{Y} \rangle} \\
    	{\langle \bs{X}, \bs{Y} \rangle} & {\langle \bs{X}, \bs{Y} \rangle} & {\bs{}} \\
    	& {\bs{}.}
    	\arrow[""{name=0, anchor=center, inner sep=0}, "{\rps{\bs{X}}{q}}"', squiggly, from=1-2, to=1-1]
    	\arrow["\refl", squiggly, from=1-2, to=2-3]
    	\arrow[""{name=1, anchor=center, inner sep=0}, "{\lps{p}{\bs{Y}} }", squiggly, from=2-1, to=1-1]
    	\arrow[squiggly, from=2-1, to=3-2]
    	\arrow[squiggly, from=2-2, to=1-2]
    	\arrow[squiggly, from=2-2, to=2-1]
    	\arrow[""{name=2, anchor=center, inner sep=0}, "{\refl_{\bs{}}}"{description}, squiggly, from=3-2, to=2-3]
    	\arrow[shorten <=6pt, shorten >=6pt, Rightarrow, from=1, to=0]
    	\arrow[shorten >=2pt, Rightarrow, from=2-2, to=2]
    \end{tikzcd}\]
    We will denote it by $\underline{[p, q]}$.
    The 2-path $\ap_i({I_{0,0}^{\Omega X, \Omega Y}}(p,q))$ can be depicted as follows:
    \[\begin{tikzcd}
    	{\bs{}} & {\bs{}} & {\bs{}} \\
    	{\bs{}} & {\bs{}} & {\bs{}} \\
    	{\bs{}} & {\bs{}} & {\bs{}.}
    	\arrow[squiggly, from=1-1, to=1-2]
    	\arrow[squiggly, from=1-1, to=2-1]
    	\arrow["{\lps{p}{\bs{Y}}}", squiggly, from=1-2, to=2-2]
    	\arrow[squiggly, from=1-3, to=1-2]
    	\arrow[squiggly, from=1-3, to=2-3]
    	\arrow["{\rps{\bs{X}}{\refl}}", squiggly, from=2-1, to=2-2]
    	\arrow["{\rps{\bs{X}}{q}}", squiggly, from=2-3, to=2-2]
    	\arrow[squiggly, from=3-1, to=2-1]
    	\arrow[squiggly, from=3-1, to=3-2]
    	\arrow["{\lps{\refl}{\bs{Y}}}"', squiggly, from=3-2, to=2-2]
    	\arrow[squiggly, from=3-3, to=2-3]
    	\arrow[squiggly, from=3-3, to=3-2]
    \end{tikzcd}\]
    Here each square corresponds to $\underline{[r,t]}$, $r : \Omega X, t: \Omega Y$. Such a square is of the form $[r,t]^{-1}$ and half of it -- the outer part of the square in this diagram -- is trivialized by $byContr$. Thus the big square of the diagram is connected to $\refl_{\bs{}}$, and the diagram should be seen as the clockwise path
    $$\lps{\refl}{\bs{Y}} \rightsquigarrow \rps{\bs{X}}{\refl} \rightsquigarrow \lps{p}{\bs{Y}} \rightsquigarrow \rps{\bs{X}}{q} \rightsquigarrow \lps{\refl}{\bs{Y}}.$$
    The image of $(p,q)$ under the ``up-right'' composite $(\Omega i)({I_{0,0}^{\Omega X, \Omega Y}}(p,q))$ equals this loop conjugated with the obvious path $\mathrm{ptd}(i): \lps{\refl}{\bs{Y}} = \refl_{\bs{}}$.
    Note that the part of the diagram outside of the inner part of the upper right square is contractible by $byContr$, and so we conclude that this image is equal to the inverse of $[p,q]^l$.

    Now we check that the ``down-right'' composite takes $(p,q)$ to $[p,q]^l$. Of course, one would expect this to be a trivial consequence of path induction, but the technical difficulties of dealing with loop spaces make it somewhat awkward.

    We are interested in $\ap_{\oout}(\oout_2 \langle p,q \rangle) \equiv \ap_{\oout}(byContr^{-1} \cdot \rps{p}{q} \cdot byContr)$. Note that
    $$\ap_{\oout}(\rps{p}{q}) = \ap_{\oout}\ap_{\ctos}((\refl_p, q)) = \ap_{\lambda (p,y).\lpsl{p}{y}}((\refl_p, q)) = \ap_{\lpsl{p}{-}}(q) = \ap_{t(-)^{-1} \cdot \lps{p}{-} \cdot t(-)}(q),$$
    where $t(y) : \langle \bs{X}, y \rangle = \bs{}$ comes from contractibility. Now let $q: y =_Y y'$ be arbitrary. The above equality can be depicted in the diagram
    \[\begin{tikzcd}
    	& {\bs{}} \\
    	{\langle \bs{X}, y \rangle} && {\langle \bs{X}, y \rangle} \\
    	{\langle \bs{X}, y' \rangle} && {\langle \bs{X}, y' \rangle.}
    	\arrow[""{name=0, anchor=center, inner sep=0}, squiggly, from=2-1, to=1-2]
    	\arrow[""{name=1, anchor=center, inner sep=0}, "{\lps{p}{y}}", squiggly, from=2-1, to=2-3]
    	\arrow[""{name=2, anchor=center, inner sep=0}, squiggly, from=2-3, to=1-2]
    	\arrow[squiggly, from=2-3, to=3-3]
    	\arrow[""{name=3, anchor=center, inner sep=0}, squiggly, from=3-1, to=1-2]
    	\arrow["{\rps{\bs{X}}{q}^{-1}}", squiggly, from=3-1, to=2-1]
    	\arrow[""{name=4, anchor=center, inner sep=0}, squiggly, from=3-1, to=3-3]
    	\arrow[""{name=5, anchor=center, inner sep=0}, squiggly, from=3-3, to=1-2]
    	\arrow[shorten <=2pt, shorten >=2pt, Rightarrow, from=0, to=3]
    	\arrow[shorten <=4pt, shorten >=4pt, Rightarrow, from=1, to=4]
    	\arrow[shorten <=2pt, shorten >=2pt, Rightarrow, from=2, to=5]
    \end{tikzcd}\]
    Here the boundaries of the front and back triangles are $\oout(\langle p,y'\rangle)$ and $\oout(\langle p,y \rangle)$ and the side triangles are $\ap_t(q)$. The square corresponds to $\ap_{\lps{p}{-}}(q)$. Note that the side triangles are $byContr$.
    Then $\ap_{\oout}(\oout_2 \langle p,q \rangle) = \ap_{\oout}(byContr^{-1}) \cdot  \ap_{\oout}(\rps{p}{q}) \cdot \ap_{\oout}(byContr)$ is the above diagram where the front and back sides of the square are also contracted by $byContr$. Thus, it is enough to recognize the square as $[p, q]$. Our argument is like that at the end of the proof of Lemma \ref{le:siin-anticom}. The type of the square is $\rps{\bs{X}}{q}^{-1} \cdot \lps{p}{y} \cdot \rps{\bs{X}}{q} = \lps{p}{y'},$
    which by whiskering and cancellations is equivalent to $\lps{p}{y} \cdot \rps{\bs{X}}{q} = \rps{\bs{X}}{q} \cdot \lps{p}{y'}$  -- the type of $[p, q]$. Since we are interested in the square after the side contraction, it remains to note that whiskering with sides, canceling, and then contracting is the same as contracting right away.

    Now let us induct on $n$. The inductive step follows from the diagram
    \[\begin{tikzcd}
    	{\Omega^{n} X \times \Omega^{m} Y} & {\Omega^{n+m-2}(\Omega^{2} X \join\Omega Y)} && {\Omega^{n+m-1}(\Omega^{} X \join \Omega^{} Y)} \\
    	{\Omega^{n}X \wedge \Omega^{m}Y} & {\Omega^{n+m-1}(\Omega^{}X \wedge Y)} && {\Omega^{n+m}(X \wedge Y).}
    	\arrow["{-I_{n-2,m-1}^{\Omega^{2}X, \Omega Y}}", from=1-1, to=1-2]
    	\arrow["\ctos", from=1-1, to=2-1]
    	\arrow["{\Omega^{n+m-2}\oout}", from=1-2, to=1-4]
    	\arrow[from=1-2, to=2-2]
    	\arrow[from=1-4, to=2-4]
    	\arrow[from=2-1, to=2-2]
    	\arrow["{\Omega^{n+m-1}\oout}", from=2-2, to=2-4]
    \end{tikzcd}\]
    The left square commutes up to
    $-(-1)^{n-1}$
    by induction. We are left to show the commutativity of the right square.

    We show that for $X, Y$ pointed spaces, the diagram
        \[\begin{tikzcd}
        	{\Omega^2 X \join \Omega Y} & {\Omega(\Omega X \join \Omega Y)} \\
        	{\Omega(\Omega X \wedge Y)} & {\Omega^2(X \wedge Y)}
        	\arrow["\oout", from=1-1, to=1-2]
        	\arrow[from=1-1, to=2-1]
        	\arrow[from=1-2, to=2-2]
        	\arrow["{\Omega \oout}"', from=2-1, to=2-2]
        \end{tikzcd}\]
        commutes.

    Denote the rightmost map by $\Omega i$.
    To find the image of a generic path $\inl(\eta) \rightsquigarrow \inr(q)$ in the domain, we need to provide the upper and lower sides and the filler for the square
    \[\begin{tikzcd}
    	{\Omega i\ap_{\inl}(\eta)} && {\Omega \oout (\lps{\eta}{\bs{Y}})} \\
    	\refl && {\Omega \oout (\rps{\refl}{q}).}
    	\arrow[squiggly, from=1-1, to=1-3]
    	\arrow["{\ap_{\Omega i}c(\eta, q)}"', squiggly, from=1-1, to=2-1]
    	\arrow["{\ap_{\Omega \oout }(\underline{[\eta, q]})}"', squiggly, from=1-3, to=2-3]
    	\arrow[squiggly, from=2-1, to=2-3]
    \end{tikzcd}\]

    We generalize both composites to the case $\tot{p: \Omega X}{\refl=p}$ in place of $\Omega \Omega X$ to make path induction on $\eta$ applicable. For the down-right composite, send $\inl(\eta)$ to $\ap_{\oout}(\lps{\eta}{\bs{Y}}):\lpsl{\refl}{\bs{Y}} =_{\Omega (X \wedge Y)} \lpsl{p}{\bs{Y}}$, $\inr(q)$ to $byContr \cdot \ap_{\lpsl{p}{-}}(q^{-1}): \lpsl{\refl}{\bs{Y}} \rightsquigarrow \lpsl{p}{\bs{Y}} \rightsquigarrow \lpsl{p}{\bs{Y}}$. Note that the latter path equals $byContr$ when $p \equiv \refl$, so this generalizes $\Omega \oout(\rps{\refl}{q})$. For $\glue(\eta, q)$, we fill the diagram
    \[\begin{tikzcd}
    	{\lpsl{\refl}{\bs{Y}}} & {\lpsl{\refl}{\bs{Y}}} \\
    	{\lpsl{p}{\bs{Y}}} & {\lpsl{p}{\bs{Y}},}
    	\arrow["q", squiggly, from=1-1, to=1-2]
    	\arrow["\refl", curve={height=-12pt}, squiggly, from=1-1, to=1-2]
    	\arrow["\eta", squiggly, from=1-1, to=2-1]
    	\arrow["\eta"', squiggly, from=1-2, to=2-2]
    	\arrow["byContr", shift left=3, curve={height=-12pt}, squiggly, from=1-2, to=2-2]
    	\arrow["q", squiggly, from=2-1, to=2-2]
    \end{tikzcd}\]
    where the
    ``triangles'' are $byContr$ and the ``square'' is the obvious filler for the parallelepiped with two dimensions spanned by $\eta$ and the third by $q$ (a three-dimensional version of squares $\lambda p.\lambda q.[p,q]$). Clearly, this generalizes $\ap_{\Omega \oout}(\underline{[\eta, q]})$, since $ \ap_{\ap_{\oout}}\ap_{\ap_{{\ctos}}}([\eta, q]) = \ap_{\ap_{\lambda (p,y).\lpsl{p}{y}}}([\eta, q])$ is the same parallelepiped.

    We generalize the second composite as a composite $(\refl =_{\Omega X} p) \join \Omega Y \to (\inl(\refl) =_{\Omega X \join \Omega Y} \inl(p)) \xrightarrow{\ap_i} \lpsl{\refl}{\bs{Y}} =_{\Omega (X \wedge Y)} \lpsl{p}{\bs{Y}}$, with both maps defined by the same formulas. Namely, the first map sends $\inl(\eta)$ to $\ap_{\inl}(\eta): \inl(\refl) = \inl(p)$, $\inr(q)$ to $\glue(\refl, q) \cdot \glue^{-1}(p,q),$ and $\glue(\eta, q)$ to $\cone_q(\refl, p, \eta)$.

    Now it is enough to observe that the two functions are equal on $\inl(\eta)$ for $\eta \equiv \refl_{\refl}$ (note that $\sgtn \star \Omega Y$ is contractible).

    When $n=1$, we induct on $m$. Note that our map $\Omega X \times \Omega^m Y \to \Omega^{m+1}(X \wedge Y)$ still ends with $\oout$, so to induct we need to permute $\oout$ and $\oout_2$:
    \begin{align*}
        \Omega^mi \circ \Omega^{m-1} \oout_2 \circ -I_{0, m-2}^{\Omega X, \Omega^2Y} = \\
        -\Omega^m \oout_2 \circ \Omega^{m-1}i \circ I_{0, m-2}^{\Omega X, \Omega^2Y} = \\
        -\Omega^m \oout_2 \circ (\Omega^{m-1} \oout \circ \ooutk{m-1}_2 \circ \ctos) = \\
        \Omega^m \oout \circ \Omega^{m-1} \oout_2 \circ \ooutk{m-1}_2 \circ \ctos = \\
        \Omega^m \oout \circ  \ooutk{m}_2 \circ \ctos.
    \end{align*}
    Here the first step uses that $i$ and $\oout_2$ commute -- this can be proven like for $\oout$ above, or reduced to it by swapping left and right factors and checking that $i$ commutes with such a swap. The second step uses the induction hypothesis, and the third step permutes $\oout$ and $\oout_2$ (Lemma \ref{le:oout-anticom-w}).
\end{proof}

\begin{lem}\label{le:ap-composite}
    For any types $X, Y$, $x,x':X, r: x=x', p,q: X \to \Omega Y$,
    $$\ap_{\lambda x.p(x) \cdot q(x)}(r) = (\ap_{p}(r) \whisk q(x)) \cdot (p(x') \whisk \ap_{q}(r)). $$
    In particular, $\ap_{\lambda x.p \cdot q(x)}(r) = p \whisk \ap_q(r)$ and $\ap_{\lambda x.p(x) \cdot q}(r) = \ap_p(r) \whisk (q)$.
\end{lem}
\begin{proof}
    By path induction.
\end{proof}

\begin{re}\label{re:ap-composite}
    More generally, given functions $a,b,c: X \to Y$, a path $r: x = x'$, and paths with varying endpoints $p:\depfun{x:X}{ax = bx}$, $q:\depfun{x:X}{bx = cx}$, the cell $\ap_{\lambda x.p(x) \cdot q(x)}(r)$ can be computed from $\ap_{p(-)}(r)$ and $\ap_{q(-)}(r)$ via the same formula if we see transport along $r$ in the type of paths as pre- and post-composing with the values of the endpoints on $r$:
    \[\begin{tikzcd}
    	ax & bx & cx \\
    	{ax'} & {bx'} & {cx'.}
    	\arrow[""{name=0, anchor=center, inner sep=0}, "px", squiggly, from=1-1, to=1-2]
    	\arrow["{\ap_a r}"', squiggly, from=1-1, to=2-1]
    	\arrow[""{name=1, anchor=center, inner sep=0}, "qx", squiggly, from=1-2, to=1-3]
    	\arrow["{\ap_br}"{description}, squiggly, from=1-2, to=2-2]
    	\arrow["{\ap_c r}", squiggly, from=1-3, to=2-3]
    	\arrow[""{name=2, anchor=center, inner sep=0}, "{px'}"', squiggly, from=2-1, to=2-2]
    	\arrow[""{name=3, anchor=center, inner sep=0}, "{qx'}"', squiggly, from=2-2, to=2-3]
    	\arrow["{\ap_pr}"{description}, between={0.2}{0.8}, Rightarrow, from=0, to=2]
    	\arrow["{\ap_qr}"{description}, between={0.2}{0.8}, Rightarrow, from=1, to=3]
    \end{tikzcd}\]
\end{re}

Note that the homotopy $c_{\cdot,\cdot}$ from Definition \ref{de:c} can be seen as witnessing the naturality (Definition \ref{de:nat-tr}) of the transformation $A \wedge - \implies B \wedge -$.
\begin{lem}\label{le:c-natural}
    The natural transformation $A \wedge - \implies B \wedge -$ is coherent:
    given $f : X \pto Y, g: Y \pto Z$, $h : A \pto B$, there is a filler in
    \[\begin{tikzcd}
    	{A \wedge X} & {A \wedge Y} \\
    	{ B\wedge X} & {B \wedge Y} & {A \wedge Z} \\
    	&& {B \wedge Z,}
    	\arrow[from=1-1, to=1-2]
    	\arrow[from=1-1, to=2-1]
    	\arrow[from=1-1, to=2-3]
    	\arrow[from=1-2, to=2-2]
    	\arrow[from=1-2, to=2-3]
    	\arrow[from=2-1, to=2-2]
    	\arrow[from=2-1, to=3-3]
    	\arrow[from=2-2, to=3-3]
    	\arrow[from=2-3, to=3-3]
    \end{tikzcd}\]
    where the triangle cells witness that $A \wedge -$ and $B \wedge -$ preserve compositions.
\end{lem}
\begin{proof}
    By path induction, we can assume that $f,g,h$ are pointed by $\refl$. All 2-cells involved are equal to $\refl$ on 0-dimensional constructors, so the values of our filler on these constructors are just $\refl_{\refl}$. The action of 2-cells on the sides consists in trivial algebraic rewritings as in the proof of Lemma \ref{le:c-intro}, and they all cancel out for the juxtaposition of all 2-cells\footnote{Intuitively, because every edge belongs to exactly two sides.}.
\end{proof}

\begin{lem}\label{le:smashes-cube}
Suppose given pointed pairs
$f: A \pto B, g: C \pto D, f': A' \pto B', g': C' \pto D'$,
and maps of pairs
$(B, A) \to (B', A'), \eta: (D, C) \to (D',C'), (B,A) \to (D,C), (B',A') \to (D', C')$.
\\
Then in
\[\begin{tikzcd}
	& {A' \wedge C'} && {A' \wedge D'} \\
	{A \wedge C} && {A \wedge D} \\
	& {B' \wedge C'} && {B' \wedge D'} \\
	{B \wedge C} && {B \wedge D}
	\arrow[from=1-2, to=1-4]
	\arrow[from=1-4, to=3-4]
	\arrow[from=1-2, to=3-2]
	\arrow[from=3-2, to=3-4]
	\arrow[from=2-1, to=2-3]
	\arrow[from=2-3, to=4-3]
	\arrow[from=4-1, to=4-3]
	\arrow[from=2-1, to=4-1]
	\arrow[from=4-1, to=3-2]
	\arrow[from=4-3, to=3-4]
	\arrow[from=2-1, to=1-2]
	\arrow[from=2-3, to=1-4]
\end{tikzcd}\]
there are fillers in $U$ for the sides and the cube.
\end{lem}

\begin{proof}
This decomposes into two cubes:
\[\begin{tikzcd}
	&& {A' \wedge C'} &&& {A' \wedge D'} \\
	& {A \wedge C'} &&& {A \wedge D'} \\
	{A \wedge C} && {B' \wedge C'} & {A \wedge D} && {B' \wedge D'} \\
	& {B \wedge C'} &&& {B \wedge D'} \\
	{B \wedge C} &&& {B \wedge D.}
	\arrow[from=1-6, to=3-6]
	\arrow[dashed, from=3-3, to=3-6]
	\arrow[dashed, from=4-2, to=4-5]
	\arrow[from=2-2, to=4-2]
	\arrow[from=4-2, to=3-3]
	\arrow[from=4-5, to=3-6]
	\arrow[from=2-2, to=1-3]
	\arrow[from=2-5, to=1-6]
	\arrow[from=2-5, to=4-5]
	\arrow[from=1-3, to=3-3]
	\arrow[from=2-2, to=2-5]
	\arrow[from=1-3, to=1-6]
	\arrow[from=3-1, to=2-2]
	\arrow[from=5-1, to=4-2]
	\arrow[from=3-1, to=3-4]
	\arrow[from=5-1, to=5-4]
	\arrow[from=3-4, to=5-4]
	\arrow[from=3-1, to=5-1]
	\arrow[from=3-4, to=2-5]
	\arrow[from=5-4, to=4-5]
\end{tikzcd}\]

Denote the maps by $c: C \to C', d: D \to D'$. To find a filler for the first cube (which technically means proving $(c_{f, d} \circ (A \wedge g)) \cdot (\ap_{B \wedge d} \circ c_{f, g}) \cdot ((B \wedge \eta) \circ (f \wedge C)) =_{(f \wedge D') \circ (A \wedge d) \circ (A \wedge g) = (B \wedge g') \circ (B \wedge c) \circ (f \wedge C)} (\ap_{f \wedge D'} A \wedge \eta) \cdot (c_{f, g'} \circ A \wedge c) \cdot (\ap_{B \wedge g'} \circ c_{f, c})$), consider the natural transformation $A \wedge - \implies B \wedge -$. It is coherent by Lemma \ref{le:c-natural}. But then the cube is a special case of the one from Lemma \ref{le:nat-tr-cube}.

A similar proof fills the second cube.
\end{proof}

\begin{lem}\label{le:fiber-of-pp}
    Given pointed maps $f: A \pto B, g: C \pto D$, the fiber of $P^\times(f, g)$ is equivalent to the join of the fibers of $f$ and $g$.
\end{lem}
\begin{proof}
    By path induction, we can assume that $f$ and $g$ are pointed by reflexivity.

    We can construct equivalences directly. Denote the fibers of $f, g$ and $P^\times$ by $F$, $G$ and $X$. Then there is a map $F \join G \to X$ defined by $\inl(a, p) \mapsto (\inl(a, \bs{D}), (p, {\bs{D}})); \inr(c, q) \mapsto (\inr( \bs{B}, c), ({\bs{B}}, q)); \glue((a,p),(c,q)) \mapsto (\ap_{\inl}((a, q^{-1})) \cdot \glue(a,c) \cdot \ap_{\inr}((p, c)), s)$, where, if we see transport as precomposition, $s$ is the obvious filler in
    \[\begin{tikzcd}
    	{(\bs{B}, \bs{D})} & {(fa, \bs{D})} \\
    	{(\bs{B}, gc)} & {(fa, gc).}
    	\arrow["{(p, \bs{D})}"', squiggly, from=1-2, to=1-1]
    	\arrow[""{name=0, anchor=center, inner sep=0}, "{(\bs{B}, q)}", squiggly, from=2-1, to=1-1]
    	\arrow["{(p^{-1}, gc)}"', squiggly, from=2-1, to=2-2]
    	\arrow[""{name=1, anchor=center, inner sep=0}, "{(fa, q)}"', squiggly, from=2-2, to=1-2]
    	\arrow["s"{description}, between={0.2}{0.8}, Rightarrow, from=1, to=0]
    \end{tikzcd}\]
    Define $X \to F \join G$ by $\inl(a, d) \mapsto \lambda \operatorname{pair^=}(p, q). \inl(a, p)$; $\inr(b, c) \mapsto \lambda \operatorname{pair^=}(p, q).\inr(c, q)$; $\glue(a, c) \mapsto \funext(\lambda \operatorname{pair^=}(p,q).\glue((a,p),(c,q)))$. It is easy to show that these are inverses.
\end{proof}

\begin{re}
    The lemma also follows from the descent theorem (\cite{BvDR}, Theorem 2.2) applied to the cube
    \[\begin{tikzcd}
    	& {F \times G} \\
    	G & {A \times C} & F \\
    	{B \times C} & \sgtn & {A \times D} \\
    	& {B\times D,}
    	\arrow[from=1-2, to=2-1]
    	\arrow[from=1-2, to=2-2]
    	\arrow[from=1-2, to=2-3]
    	\arrow[from=2-1, to=3-1]
    	\arrow[from=2-1, to=3-2]
    	\arrow["{\ }"{description}, from=2-2, to=3-1]
    	\arrow["{\ }"{description}, from=2-2, to=3-3]
    	\arrow[from=2-3, to=3-2]
    	\arrow[from=2-3, to=3-3]
    	\arrow[from=3-1, to=4-2]
    	\arrow[from=3-2, to=4-2]
    	\arrow[from=3-3, to=4-2]
    \end{tikzcd}\]
    but we use the map $F \join G \to X$ explicitly in the proof of Proposition \ref{pr:smashing-proof}.
\end{re}

\subsubsection{Geometric lemmas about spheres}

In the following lemmas we extensively use algebraic properties of the minus sign in homotopy groups:
\begin{itemize}
    \item $\merid_{S^{n-1}}^{-1}$ induces $-[\id]: \pi_n(S^n)$;
    \item $-(-f) = f$;
    \item $\Sigma(-f) = -\Sigma f$;
    \item $(-f) \circ g = -(f \circ g)$;
    \item $(-f) \whisk \alpha = (f \whisk \alpha)^{-1}$.
\end{itemize}

Denote by $\merid^k$ the meridian of $S^k \equiv \Sigma S^{k-1}$.
\begin{lem}\label{le:Sn-smash-merid}
    For any $n \geq 1, m \geq 0,$ the based 2-cell
    \[\begin{tikzcd}
    	{S^{n+m-1}} & {S^{m,n-1}} & \sgtn \\
    	& \sgtn & {S^{m,n}} & {S^{n+m}}
    	\arrow["{j^{m,n-1}}", from=1-1, to=1-2]
    	\arrow[from=1-2, to=1-3]
    	\arrow[from=1-2, to=2-2]
    	\arrow["{S^m \wedge \merid^n}"{description}, Rightarrow, from=1-3, to=2-2]
    	\arrow[from=1-3, to=2-3]
    	\arrow[from=2-2, to=2-3]
    	\arrow["{({j^{m,n}})^{-1}}", from=2-3, to=2-4]
    \end{tikzcd}\]
    equals $(-1)^{m} \merid^{n+m}$.
\end{lem}
\begin{re}
    Technically, for $S^m \wedge \merid^n$ we have to write $S^m \wedge \sgtn$ in the corners. We tacitly transport along the equality $S^m \wedge \sgtn =_{U_*} \sgtn$.
\end{re}
\begin{proof}
    We induct on $m$ and $n$. Recall that the identification $j^{n,m}: S^{n+m} \to S^n \wedge S^m$ is by definition the composite $\siin \circ \Sigma j^{n-1,m}$ for $n>0$ and $\siin_2 \circ \Sigma j^{n,m-1}$ for $n=0$.

    We first directly check the case $m=0, n=1$:
    \begin{align*}
    (({j^{0,1}})^{-1} \whisk (S^0 \wedge \merid^1) \whisk j^{0,0})(1_{\s}) = \\
    \ap_{({j^{0,1}})^{-1}}(S^0 \wedge \merid^1)\langle 1_{\s}, 1_{\s} \rangle \equiv  \\
    \ap_{({j^{0,1}})^{-1}} \langle 1_{\s}, \merid^1_{1_{\s}}) = \\
    \ap_{(\Sigma {j^{0,0}})^{-1}} \ap_{\sout_2} \langle 1_{\s}, \merid^1_{1_{\s}} ) = \\
    \ap_{(\Sigma {j^{0,0}})^{-1}} \merid^1_{\langle 1_{\s},1_{\s} \rangle}
    = \\
    \merid^1_{1_{\s}}.
    \end{align*}
    Similarly, $0_{\s}$ gets sent to $\merid^1_{0_{\s}}$, and the check that this is in $U_*$ is formal.

    Assume now that $m=0$ and $n>1$. We have
    \begin{align*}
        ({j^{0,n}})^{-1} \whisk (S^0 \wedge \merid^n) \whisk j^{0,n-1} =_* \\
        ((\Sigma {j^{0,n-1}})^{-1} \circ \siin_2^{-1}) \whisk (S^0 \wedge \merid^n) \whisk (\siin_2 \circ \Sigma j^{0,n-2}) =_* \\
        (\Sigma {j^{0,n-1}})^{-1} \whisk (\siin_2^{-1} \whisk (S^0 \wedge \merid^n) \whisk \siin_2) \whisk \Sigma j^{0,n-2} =_* \\
        -(\Sigma {j^{0,n-1}})^{-1} \whisk (\siin_2^{-1} \whisk (S^0 \wedge \Sigma \merid^{n-1}) \whisk \siin_2) \whisk \Sigma j^{0,n-2} =_* && (*) \\
        -(\Sigma {j^{0,n-1}})^{-1} \whisk \Sigma (S^0 \wedge \merid^{n-1}) \whisk \Sigma j^{0,n-2} =_* \\
        -\Sigma (({j^{0,n-1}})^{-1} \whisk (S^0 \wedge \merid^{n-1}) \whisk j^{0,n-2}) =_* \\
        -\Sigma \merid^{n-1} =_* \merid^n.
    \end{align*}
    Here the equality $(*)$ follows from the naturality of $\siin_2$, see the discussion after Definition \ref{de:nat-tr}. Besides this, we use the fact that $\Sigma \merid^k =_* -\merid^{k+1}$ (Lemma \ref{le:susp-of-merid}), and induction on $n$.

    For the $m$-inductive step, we get essentially the same calculation, but this time to split $\Sigma$ from $S^m \wedge \merid^n$, we split it from $S^m$ rather than $\merid^n$. Since $S^m = \Sigma S^{m-1}$, this will involve one sign change less, hence the sign alternation every time we decrease $m$:
    \begin{align*}
    ({j^{m,n}})^{-1} \whisk (S^m \wedge \merid^n) \whisk j^{m,n-1} =_* \\
    ((\Sigma {j^{m-1,n}})^{-1} \circ \siin^{-1}) \whisk (S^m \wedge \merid^n) \whisk (\siin \circ \Sigma j^{m-1,n-1}) =_* \\
    (\Sigma {j^{m-1,n}})^{-1} \whisk (\siin^{-1} \whisk (S^m \wedge \merid^n) \whisk \siin) \whisk \Sigma j^{m-1,n-1} =_* \\
    (\Sigma {j^{m-1,n}})^{-1} \whisk (\siin^{-1} \whisk  (\Sigma S^{m-1} \wedge \merid^{n})  \whisk \siin) \whisk \Sigma j^{m-1,n-1} =_* \\
    (\Sigma {j^{m-1,n}})^{-1} \whisk \Sigma (S^{m-1} \wedge \merid^{n}) \whisk \Sigma j^{m-1,n-1} =_* \\
    \Sigma (({j^{m-1,n}})^{-1} \whisk (S^{m-1} \wedge \merid^{n}) \whisk j^{m-1,n-1}) =_* \\
    (-1)^{m-1} \Sigma \merid^{n+m-1} =_* (-1)^m\merid^{n+m}.
    \end{align*}
\end{proof}

\begin{lem}\label{le:j1-ind}
Let $m \geq 1$. Then $j^{1, m} =_* -\siin_2 \circ \Sigma j^{1, m-1}$.
\end{lem}
\begin{proof}
    \begin{align*}
        j^{1, m} \equiv \\
        \siin \circ \Sigma j^{0, m} \equiv \\
        \siin \circ \Sigma (\siin_2 \circ \Sigma j^{0, m-1}) =_* \\
        \siin \circ \Sigma \siin_2 \circ \Sigma^2 j^{0, m-1} =_* && (*) \\
        -\siin_2 \circ \Sigma \siin \circ \Sigma^2 j^{0, m-1} =_*  \\
        -\siin_2 \circ \Sigma (\siin \circ \Sigma j^{0, m-1}) \equiv \\
        -\siin_2 \circ \Sigma j^{1, m-1}.
    \end{align*}
    Here the only nonformal step is $(*)$, which uses Lemma \ref{le:siin-anticom}.
\end{proof}

\begin{lem}\label{le:merid-smash-Sm}
    For any $n \geq 1, m \geq 0,$ the based 2-cell
    \[\begin{tikzcd}
    	{S^{n+m-1}} & {S^{n-1,m}} & \sgtn \\
    	& \sgtn & {S^{n,m}} & {S^{n+m}}
    	\arrow["{j^{n-1,m}}", from=1-1, to=1-2]
    	\arrow[from=1-2, to=1-3]
    	\arrow[from=1-2, to=2-2]
    	\arrow["{\merid^n \wedge S^m}"{description}, Rightarrow, from=1-3, to=2-2]
    	\arrow[from=1-3, to=2-3]
    	\arrow[from=2-2, to=2-3]
    	\arrow["{(j^{n,m})^{-1}}", from=2-3, to=2-4]
    \end{tikzcd}\]
    equals $\merid^{n+m}$.
\end{lem}
\begin{proof}
    The proof proceeds essentially like that of Lemma \ref{le:Sn-smash-merid}.

    First, check the base case: assume $n = 1$ and $m=0$. We compute:
    \begin{align*}
    (({j^{1,0}})^{-1} \whisk (\merid^1 \wedge S^0) \whisk j^{0,0}) (1_{\s}) = \\
    \ap_{({j^{1,0}})^{-1}} (\merid^1 \wedge S^0)\langle 1_{\s}, 1_{\s} \rangle \equiv \\
    \ap_{({j^{1,0}})^{-1}} (\merid^1_{1_{\s}}, 1_{\s} \rangle = \\
    \ap_{(\Sigma {j^{0,0}})^{-1}} \ap_{\sout} (\merid^1_{1_{\s}}, 1_{\s} \rangle = \\
    \ap_{(\Sigma {j^{0,0}})^{-1}} \merid^1_{\langle 1_{\s}, 1_{\s} \rangle}
    = \\
    \merid^1_{1_{\s}}.
    \end{align*}
    Similarly, $0_{\s}$ gets sent to $\merid^1_{0_{\s}}$.

    Now use induction. Let $n>1$. Then
    \begin{align*}
    ({j^{n,m}})^{-1} \whisk (\merid^n \wedge S^m) \whisk j^{n-1, m} =_* \\
    ((\Sigma {j^{n-1,m}})^{-1} \circ \siin^{-1}) \whisk (\merid^n \wedge S^m) \whisk (\siin \circ \Sigma j^{n-2,m}) =_* \\
    (\Sigma {j^{n-1,m}})^{-1} \whisk (\siin^{-1} \whisk (\merid^n \wedge S^m) \whisk \siin) \whisk \Sigma j^{n-2,m} =_* \\
    - (\Sigma {j^{n-1,m}})^{-1} \whisk (\siin^{-1} \whisk (\Sigma \merid^{n-1} \wedge S^m) \whisk \siin) \whisk \Sigma j^{n-2,m} =_* \\
    - (\Sigma {j^{n-1,m}})^{-1} \whisk  \Sigma (\merid^{n-1} \wedge S^m) \whisk \Sigma j^{n-2,m} =_* \\
    - \Sigma (({j^{n-1,m}})^{-1} \whisk (\merid^{n-1} \wedge S^m) \whisk j^{n-2,m}) =_* \\
    - \Sigma \merid^{n+m-1} =_* \merid^{n+m}.
    \end{align*}

    Now assume $n=1$ and induct on $m$. We have
    \begin{align*}
    ({j^{1,m}})^{-1} \whisk (\merid^1 \wedge S^m) \whisk j^{0, m} =_* && (*) \\
    -((\Sigma {j^{1,m-1}})^{-1} \circ \siin_2^{-1}) \whisk (\merid^1 \wedge S^m) \whisk (\siin_2 \circ \Sigma j^{0,m-1}) =_* \\
    -(\Sigma {j^{1,m-1}})^{-1} \whisk (\siin_2^{-1} \whisk (\merid^1 \wedge \Sigma S^{m-1}) \whisk \siin_2) \whisk \Sigma j^{0,m-1} =_* \\
    -(\Sigma {j^{1,m-1}})^{-1} \whisk \Sigma(\merid^1 \wedge S^{m-1}) \whisk \Sigma j^{0,m-1} =_* \\
    -\Sigma (({j^{1,m-1}})^{-1} \whisk (\merid^1 \wedge S^{m-1}) \whisk j^{0,m-1}) =_* \\
    -\Sigma (\merid^{m}) =_* \\
    \merid^{m+1}.
    \end{align*}
    Here $(*)$ uses Lemma \ref{le:j1-ind}.
\end{proof}

\begin{lem}\label{le:susp-of-merid}
    For any $A: U_*$, there is a filler in $U_*$ for the cube
    \[\begin{tikzcd}
    	{\Sigma A} &&& {\Sigma \sgtn} \\
    	& {\Sigma A} & \sgtn \\
    	& \sgtn & {\Sigma \Sigma A} \\
    	{\Sigma \sgtn} &&& {\Sigma \Sigma A}
    	\arrow[from=1-1, to=1-4]
    	\arrow["\id", from=1-1, to=2-2]
    	\arrow[from=1-1, to=4-1]
    	\arrow["\sim"', from=1-4, to=2-3]
    	\arrow[from=1-4, to=4-4]
    	\arrow[from=2-2, to=2-3]
    	\arrow[from=2-2, to=3-2]
    	\arrow[from=2-3, to=3-3]
    	\arrow[from=3-2, to=3-3]
    	\arrow["\sim", from=4-1, to=3-2]
    	\arrow[from=4-1, to=4-4]
    	\arrow["\id"', from=4-4, to=3-3]
    \end{tikzcd}\]
    where the front is $\Sigma \merid$, the back is $-\merid: \lambda \_.S = \lambda \_ .N$, and the other 2-cells are self-evident (e.g. the bottom is $N \mapsto \refl_N, S \mapsto \merid_S, \merid_{*} \mapsto cancel(\merid_S)$).
\end{lem}
\begin{proof}
    Denote the meridian of $\Sigma (\Sigma A)$ by $\merid^0$ and the meridian of $\Sigma A$ by $\merid^1$.

    Denote by $f$ the composite $\Sigma A \to \Sigma A \to \sgtn \xrightarrow{S} \Sigma \Sigma A$ and by $g$ the composite $\Sigma A \xrightarrow{\mathrm{down}} \Sigma \sgtn \to \Sigma \Sigma A \xrightarrow{\id} \Sigma \Sigma A$. We have two proofs of $f = g$ and we have to prove their equality. Let $l$ be the composite of the back, left and bottom sides of the cube and $r$ be the composite of the other three sides. We have to provide equalities $l(N) = r(N)$ and $l(S) = r(S)$ and prove their equality as we move from $N$ to $S$ along $\merid^1_a$. This is encapsulated in the following cube:
    \[\begin{tikzcd}
    	N &&& N \\
    	& S & S \\
    	& S & S \\
    	S &&& {S.}
    	\arrow["\refl"{description}, squiggly, from=1-1, to=1-4]
    	\arrow["{\merid^0_S}"', squiggly, from=1-1, to=4-1]
    	\arrow["{\merid^0_S}", squiggly, from=1-4, to=4-4]
    	\arrow["{(\merid_N^0)^{-1}}"{description}, squiggly, from=2-2, to=1-1]
    	\arrow["\refl", squiggly, from=2-2, to=2-3]
    	\arrow["\refl", squiggly, from=2-2, to=3-2]
    	\arrow["{(\merid_N^0)^{-1}}"{description}, squiggly, from=2-3, to=1-4]
    	\arrow["\refl", squiggly, from=2-3, to=3-3]
    	\arrow["\refl"', squiggly, from=3-2, to=3-3]
    	\arrow["{(\merid^0_S)^{-1} \merid_S^0}"{description}, squiggly, from=3-2, to=4-1]
    	\arrow["\refl", squiggly, from=3-3, to=4-4]
    	\arrow["\refl"', squiggly, from=4-1, to=4-4]
    \end{tikzcd}\]
    Here the horizontal direction corresponds to the equality $l \rightsquigarrow r$, vertical to $\merid^1$, and the third dimension corresponds to the direction of $l$ and $r$. The top left edge is $l(N): f(N) = g(N)$; the top right edge is $r(N): f(N) = g(N)$; the top face is our choice and we take it to be $\refl$; the bottom left edge is $l(S): f(S) = g(S)$; the bottom right edge is $r(S): f(S) = g(S)$; the bottom face is our choice and we take it to be $cancel$;
    the back left and right edges are $\ap_{f}(\merid^1_a)$, the front left and right edges are $\ap_{g}(\merid^1_a)$;
    the left face is $\ap_{l}(\merid^1_a)$; the right face is $\ap_{r}(\merid^1_a)$;
    the back face is $\refl_{\ap_f({\merid^1_{(-)}})} = \refl_{\refl_S}$;
    the front is $\refl_{\ap_g({\merid^1_{(-)}})} = \refl_{\merid^0_S}$.

    Thus, to find a filler we need to prove that $\ap_{l}(\merid^1_a)$ composed with $\ap_{r}(\merid^1_a)$ gives $ \refl_{(\merid^0_N)^{-1}}$.

    Recall first how $\ap_{\Sigma \merid^1}(\merid^1_a)$ works. Given two maps $g,f: X \to Y$ and a homotopy $h: f = g$, the homotopy $\Sigma h: \Sigma f = \Sigma g$ is defined by $N \mapsto \refl_N, S \mapsto \refl_S$ and as we move along $\merid_x$, we have to fill the diagram

\[\begin{tikzcd}
	N \\
	{S.}
	\arrow["{\merid(f(x))}", curve={height=-12pt}, squiggly, from=1-1, to=2-1]
	\arrow["{\merid(g(x))}"', curve={height=12pt}, squiggly, from=1-1, to=2-1]
\end{tikzcd}\]
    This is done via $\ap_{\merid}(h(x)).$

    Going back to the cube, we deal with the cell
    \[\begin{tikzcd}
    	N & N & S & N & N \\
    	S & N & S & S & {S.}
    	\arrow["{\merid^0_S}"', squiggly, from=1-1, to=2-1]
    	\arrow["\refl"', squiggly, from=1-2, to=1-1]
    	\arrow["\refl"', squiggly, from=1-2, to=2-2]
    	\arrow[""{name=0, anchor=center, inner sep=0}, "{(\merid^0_N)^{-1}}"', squiggly, from=1-3, to=1-2]
    	\arrow["{(\merid^{0}_N)^{-1}}", squiggly, from=1-3, to=1-4]
    	\arrow["\refl", squiggly, from=1-3, to=2-3]
    	\arrow["\refl", squiggly, from=1-4, to=1-5]
    	\arrow[""{name=1, anchor=center, inner sep=0}, "{\merid^{0}_N}"{description}, squiggly, from=1-4, to=2-4]
    	\arrow[""{name=2, anchor=center, inner sep=0}, "{\merid^0_S}", squiggly, from=1-5, to=2-5]
    	\arrow["{\merid^0_S}"', squiggly, from=2-2, to=2-1]
    	\arrow[""{name=3, anchor=center, inner sep=0}, "{(\merid^0_S)^{-1}}", squiggly, from=2-3, to=2-2]
    	\arrow["\refl", squiggly, from=2-3, to=2-4]
    	\arrow["\refl"', squiggly, from=2-4, to=2-5]
    	\arrow["\eta"', between={0.2}{0.8}, Rightarrow, from=0, to=3]
    	\arrow["\beta"{description}, between={0.2}{0.8}, Rightarrow, from=1, to=2]
    \end{tikzcd}\]
         Each of the cells $l$ and $r$ is a juxtaposition of two cells (2 of the 6 cube sides are trivial); accordingly, $\ap_l(\merid^1_a)$ and $\ap_r(\merid^1_a)$ here start from the middle and go left and right, respectively (see Remark \ref{re:ap-composite}). The top row is the values at $N$ and the bottom is the values at $S$. Accordingly, $\beta \equiv \ap_{\merid^0_{}}(\merid^1_a)$ and $\eta \equiv \ap_{(\merid^0)^{-1}}(\merid^1_a)$, while the other two 2-cells are cancellations. Note that $\eta = \ap_{\lambda p.p^{-1}}(\beta)$. Now the statement follows from the case where $\beta$, $\merid^0_S$ and $\merid^0_N$ are arbitrary paths, which is immediate from path induction.

         It is also not hard to see that this 2-cell equality holds in $U_*$.
\end{proof}

\subsubsection{Lemmas about the pointed universe}

\begin{lem}
    The wild category $U_*$ satisfies the pentagon axiom (Definition \ref{de:pentagon}).
\end{lem}
\begin{proof}
    Assume given 4 composable arrows $a, b, c, d$. By path induction, we can assume that each is pointed by reflexivity. But then, each associator (such as $((ab)c, \ap_{ab}(\ptd{c})\ptd{ab})  =(a(bc), \ap_a(\ptd{bc}) \ptd{a})$) computes to reflexivity, and a filler for the pentagon exists.
\end{proof}

We prove that pushouts behave in $U_*$ as expected.
\begin{lem}\label{pushout-proof}
    The pointed pushout $(C \coprod^A B, \inl(\bs{C}))$ is the pushout of the span $C \xleftarrow{f} A \xrightarrow{g} B$ in $U_*$. This means that
    \begin{enumerate}
    \item the pushout square can be considered a square in $U_*$ and
    \item given a diagram in $U_*$ of the form
    \[\begin{tikzcd}
    	A & C \\
    	B & {D,}
    	\arrow["f", from=1-1, to=1-2]
    	\arrow["g"', from=1-1, to=2-1]
    	\arrow["\beta"', shorten <=4pt, shorten >=4pt, Rightarrow, from=1-2, to=2-1]
    	\arrow["l", from=1-2, to=2-2]
    	\arrow["r"', from=2-1, to=2-2]
    \end{tikzcd}\]
    one gets
    \begin{itemize}
        \item
        a map $I: C \coprod^A B \pto D$;
        \item
        equalities $T_L: I \circ \inl =_* l$, $T_R: I \circ \inr =_* r$;
        \item
        an equality between based homotopies $\beta$ and $(T_L \whisk f)^{-1} \cdot (I \whisk \glue) \cdot (T_R \whisk g)$.
    \end{itemize}
    Moreover, the type of the induced data above is contractible.
    \end{enumerate}
\end{lem}

\begin{re}
    More concisely, the above could be expressed by saying that in the wild category $U_*$, the pushout diagram forms a cone $\mathcal{C}$ over the span $S :\equiv C \xleftarrow{f} A \xrightarrow{g} B,$ and for any object $D$, there is an equivalence of types $ \operatorname{Hom}_{U_*}(C \coprod^A B, D) \simeq \operatorname{Cone}(S, D)$, given by postcomposition (whiskering) of $\mathcal{C}$.
\end{re}
\begin{proof}
    1) The map $\inl: C \to C \coprod^A B$ is pointed by $\refl_{\inl({\bs{C}})}$.
    The map $\inr: B \to C \coprod^A B$ is pointed so that $A, B, C, C \coprod^A B$ form a commutative square in $U_*$; that is, by $\ap_{\inr}{\ptd{g}^{-1}} \cdot \glue_{\bs{A}}^{-1} \cdot \ap_{\inl} \ptd{f}.$

    2) Given a square in $U_*$ as above,
    one gets a map (in $U$) $I: C \coprod^A B \to D$
    and a witness $w: \depfun{a:A}{\ap_I(\glue_a) = \beta(a)}$,
    and proofs (in $U$) $T_L: I \circ \inl = l$, $T_R: I \circ \inr = r$ which are reflexivity since the computation rules for pushouts are judgmental on 0-dimensional constructors. First, we need to show that $I$ is a pointed map and $T_R, T_L$ are based homotopies.

    The map $I$ is pointed by $\ptd{l}$. $T_L$ will be made based by an element $x: \ptd{I \circ \inl} \ptd{l}^{-1} = \refl_{l(\bs{C})}$ whose value we will determine later. For the basedness of $T_R$, we would need $\ap_I (\ptd{\inr}) \cdot \ptd{I} = \ptd{r}$.
    By $\bsd{\beta},$ we get:
    $\ptd{r} = \ap_r(\ptd{g})^{-1} \beta({\bs{A}})^{-1} \ap_ l({\ptd{f}}) \ptd{l};$
    by definition of $\ptd{\inr}$, naturality of $\ap$, witness $w(\bs{A})$, and basedness of $\beta$, we have
    $$\ap_I (\ptd{\inr}) = \ap_I (\ap_{\inr}{\ptd{g}^{-1}} \cdot \glue_{\bs{A}}^{-1} \cdot \ap_{\inl} \ptd{f}) = \ap_r(\ptd{g})^{-1} \beta(\bs{A})^{-1} \ap_l(\ptd{f}) = \ptd{r} \ptd{l}^{-1},$$
    giving us what we want.

    Now, we need to prove that $(T_L \whisk f)^{-1} \cdot (I \whisk \glue) \cdot (T_R \whisk g)$ equals $\beta$ as a pointed homotopy. Unwinding Remark \ref{re:equality-of-squares-simple-version} for our case, we need a 3-filler for the diagram

    \[\begin{tikzcd}
    	& {I(\inl f (\bs A))} && {I(\inr g (\bs A))} \\
    	{l(f (\bs A))} & {I(\inl(\bs{C}))} & {I(\inl(\bs{C}))} & {I(\inr(\bs B))} & {r(g(\bs A))} \\
    	& {l (\bs C)} & {\bs D} & {r(\bs B).}
    	\arrow[squiggly, from=1-2, to=1-4]
    	\arrow[squiggly, from=1-2, to=2-1]
    	\arrow[squiggly, from=1-2, to=2-2]
    	\arrow[squiggly, from=1-4, to=2-4]
    	\arrow[squiggly, from=1-4, to=2-5]
    	\arrow[curve={height=-12pt}, squiggly, from=2-1, to=2-5]
    	\arrow[squiggly, from=2-1, to=3-2]
    	\arrow[squiggly, from=2-2, to=2-3]
    	\arrow[squiggly, from=2-2, to=3-2]
    	\arrow[squiggly, from=2-3, to=3-3]
    	\arrow[squiggly, from=2-4, to=2-3]
    	\arrow[squiggly, from=2-4, to=3-4]
    	\arrow[squiggly, from=2-5, to=3-4]
    	\arrow[squiggly, from=3-2, to=3-3]
    	\arrow[squiggly, from=3-4, to=3-3]
    \end{tikzcd}\]

    Here the top front rectangle is $w(\bs{A})$, the bottom front rectangle is $\bsd{\beta}$, the ``triangles'' are $\ap_{T_L}(\ptd{f})$ and $\ap_{T_R}(\ptd{g})$, the top back rectangle is $\ap_{\ap_I}(\bsd{\glue})$, the right lower back square is $\bsd{T_R}$ and the left lower back square is $x$. Thus, if we set $x$ to be the composite of all other 2-cells of this diagram, we get the basedness for $T_L$ and a trivial 3-filler for the diagram, finishing with coherences.

    For the uniqueness, note that the argument of the previous paragraph essentially says that the last sigma-type in
    $$\tot{I: C \coprod^A B \to D}{\tot{T_L: I \circ \inl = l}{\tot{T_R: I \circ \inr = r}{\tot{w:\beta=T_L^{-1} \cdot \glue \cdot T_R}{\tot{\operatorname{ptd}(I)}{\tot{\operatorname{bsd}(T_L)}{\tot{\operatorname{bsd}(T_R)}{\operatorname{bsd}(w)}}}}}}}$$
    is contractible. Moreover, the first four together are contractible as a fiber over $\beta$ of the equivalence from \cite{UF}, Lemma 6.8.2. Finally, the following two sigmas are contractible by path induction. Thus, the total type is contractible.
    \end{proof}

\begin{lem}\label{le:octahedral}
    Let $f: A \pto B$ and $g: B \pto C$ be composable pointed maps. Then in the diagram
    \[\begin{tikzcd}
    	{\fib_{i}} & {\fib_f} \\
    	{\fib_{g \circ f}} & A & C \\
    	{\fib_g} & B & C
    	\arrow[dashed, from=1-1, to=1-2]
    	\arrow[from=1-1, to=2-1]
    	\arrow[from=1-2, to=2-2]
    	\arrow[from=2-1, to=2-2]
    	\arrow["i"', dashed, from=2-1, to=3-1]
    	\arrow["{g \circ f}", from=2-2, to=2-3]
    	\arrow["f"', from=2-2, to=3-2]
    	\arrow["\equiv", from=2-3, to=3-3]
    	\arrow[from=3-1, to=3-2]
    	\arrow["g"', from=3-2, to=3-3]
    \end{tikzcd}\]
    the upper dashed arrow (defined as the induced map of fibers) is an equivalence.
\end{lem}
\begin{re}
    This is a slightly strengthened ``unstable octahedral axiom'' from \cite{UF}, Exercise 4.4. In our case, the universe is pointed, and the equivalence is explicitly exhibited as an induced fiber map.
\end{re}
\begin{proof}
    The proof is a straightforward application of path induction. First note that we just need to show that the map in question is an equivalence in $U$ (by the remark in Definition \ref{de:pointed-equiv}). Then, by path induction we may assume that $f$ and $g$ are pointed by reflexivity. The map $i$ equals $(f \circ \pi_1, \pi_2): (a,p) \mapsto (f(a), p).$
    The fiber over $(\bs{B}, \refl_{\bs{C}})$ under $i$ is equal to
    $\tot{a:A}{\tot{p:gfa=\bs{C}}{\tot{q:fa=\bs{B}}{\ap_g(q)^{-1} \cdot p^{}=\refl_{\bs{C}}}}}$. This contracts to $\tot{a:A}{(fa=\bs{B})} \equiv \fib_f$. Note that this contraction is precisely the induced fiber map $\lambda (a,p,q,\eta).(a,q)$, so we are done.
\end{proof}

\begin{lem}\label{pasting_of_based}
    Pasting two based squares produces a based square: if in the diagram
\[\begin{tikzcd}
	A & B & C \\
	{A'} & {B'} & {C'}
	\arrow["f", from=1-1, to=1-2]
	\arrow["l"', from=1-1, to=2-1]
	\arrow["{f'}", from=1-2, to=1-3]
	\arrow["\alpha"', Rightarrow, from=1-2, to=2-1]
	\arrow["r", from=1-2, to=2-2]
	\arrow["{\alpha'}"', Rightarrow, from=1-3, to=2-2]
	\arrow["{r'}", from=1-3, to=2-3]
	\arrow["g", from=2-1, to=2-2]
	\arrow["{g'}", from=2-2, to=2-3]
\end{tikzcd}\]
the inner squares commute in $U_*$, then the whole rectangle commutes in $U_*$.
\end{lem}
\begin{proof}
    This can be easily shown explicitly, but note that pasting squares makes sense in any wild category as a composite of two whiskers.
\end{proof}

\begin{cor}\label{cor:comp_of_based}
    Assume $h: f_1 =_{A \to B} f_n$ decomposes as a sequence of homotopies $h_1 : f_1=f_2, h_2: f_2 = f_3, \ldots, h_{n-1}: f_{n-1} = f_n$ and all $h_i$ are based. Then $h$ is based.
\end{cor}

\begin{proof}
    Follows from the above, and can also be proven directly as
    $h(\bs{A}) = h_1(\bs{A}) \cdot h_2(\bs{A}) \cdot \ldots \cdot h_n(\bs{A}) = \ptd{f_1} \ptd{f_2}^{-1} \ptd{f_2} \ptd{f_3}^{-1} \ldots \ptd{f_n}^{-1} = \ptd{f_1} \ptd{f_n}^{-1}.$
\end{proof}

\begin{lem}\label{equality-of-squares}
    Given two commutative squares in $U_*$ which are the front and the back sides of the cube below, to prove their equality means to
    \begin{enumerate}
        \item Provide the sides (in $U_*$) and a filler for the cube

\[\begin{tikzcd}
	& {X'} && {Y'} \\
	& {Z'} && {W'} \\
	X && Y \\
	Z && {W;}
	\arrow["{g'}"{description}, from=1-2, to=1-4]
	\arrow["{f'}"', from=1-2, to=2-2]
	\arrow["{\beta'}"{description}, shorten <=14pt, shorten >=10pt, Rightarrow, from=1-4, to=2-2]
	\arrow["{r'}", from=1-4, to=2-4]
	\arrow["{l'}"{description}, from=2-2, to=2-4]
	\arrow[from=3-1, to=1-2]
	\arrow["g"{description}, from=3-1, to=3-3]
	\arrow["f"', from=3-1, to=4-1]
	\arrow[from=3-3, to=1-4]
	\arrow["\beta"{description}, shorten <=14pt, shorten >=14pt, Rightarrow, from=3-3, to=4-1]
	\arrow["r", from=3-3, to=4-3]
	\arrow[from=4-1, to=2-2]
	\arrow["l"{description}, from=4-1, to=4-3]
	\arrow[from=4-3, to=2-4]
\end{tikzcd}\]
        \item Prove $\ap_{\ap_{e_W}}\bsd{\beta} =_{\ap_{e_W}{\beta(\bs{X})} = \ap_{e_W} (\ptd{r \circ g} \ptd{l \circ f}^{-1})} c,$
    where $c$ is the 2-path

    \end{enumerate}

\[\makebox[\textwidth][c]{\begin{tikzcd}
{\ap_{e_W}\beta(\bs{X})} \\
	{R (g(\bs{X}))\cdot \ap_{r'} U(\bs{X}) \cdot \beta'(e_X(\bs{X})) \cdot \ap_{l'} L(\bs{X})^{-1} \cdot D(f(\bs{X}))^{-1}} \\
	{R (g(\bs{X}))\cdot \ap_{r'} U(\bs{X}) \cdot \ap_{r' \circ g'}(\ptd{e_X})  \beta'(\bs{X'}) \ap_{l' \circ f'}(\ptd{e_X})^{-1}  \cdot \ap_{l'} L(\bs{X})^{-1} \cdot D(f(\bs{X}))^{-1}} \\
	{R (g(\bs{X}))\cdot \ap_{r'} U(\bs{X}) \cdot \ap_{r' \circ g'}(\ptd{e_X})  \ptd{r' \circ g'} \ptd{l' \circ f'}^{-1} \ap_{l' \circ f'}(\ptd{e_X})^{-1}  \cdot \ap_{l'} L(\bs{X})^{-1} \cdot D(f(\bs{X}))^{-1}} \\
	{R(g(\bs{X})) \cdot \ap_{r'} (U(\bs{X}))\ap_{r'\circ g'}\ptd{ e_X}(\ap_{r'}(\ptd{g'}) \ptd{r'}) (D(f(\bs{X})) \ap_{l'} (L(\bs{X})) \ap_{l' \circ f'}\ptd{e_X}(\ap_{l'}(\ptd{f'})\ptd{l'}))^{-1}} \\
	{R(g(\bs{X})) \cdot \ap_{r'} (U(\bs{X})\ptd{g'\circ e_X}) \ptd{r'} (D(f(\bs{X})) \ap_{l'} (L(\bs{X}) \ptd{f' \circ e_X})\ptd{l'})^{-1}} \\
	{R(g(\bs{X})) \cdot \ap_{r'} (\ap_{e_Y}\ptd{g} \cdot \ptd{e_Y})\ptd{r'} (D(f(\bs{X})) \ap_{l'} (\ap_{e_Z}\ptd{f} \cdot \ptd{e_Z})\ptd{l'})^{-1}  } \\
	{R(g(\bs{X})) \cdot \ap_{r' \circ e_Y}\ptd{g} \cdot \ptd{r' \circ e_Y} (D(f(\bs{X})) \ap_{l' \circ e_Z}\ptd{f} \cdot \ptd{l' \circ e_Z})^{-1}  } \\
	{\ap_{e_W} \ap_{r} \ptd{g} \cdot R(\bs{Y}) \cdot \ptd{r' \circ e_Y}(\ap_{e_W} \ap_l \ptd{f} \cdot D(\bs{Z}) \ptd{l' \circ e_Z})^{-1}  } \\
	{\ap_{e_W} \ap_{r} \ptd{g} \cdot \ap_{e_W} \ptd{r} \cdot \ptd{e_W} (\ap_{e_W} \ap_l \ptd{f} \cdot \ap_{e_W} \ptd{l} \cdot \ptd{e_W})^{-1}} \\
	{\ap_{e_W} (\ptd{r \circ g} \ptd{l \circ f}^{-1})}
	\arrow["{\text{3-filler}}", squiggly, from=1-1, to=2-1]
	\arrow["{\ap_{\beta'}(\ptd{e_X})}", squiggly, from=2-1, to=3-1]
	\arrow["{\bsd{\beta'}}", squiggly, from=3-1, to=4-1]
	\arrow["{\bsd{U}, \bsd{L}}", squiggly, from=6-1, to=7-1]
	\arrow["\text{algebra}", squiggly, from=7-1, to=8-1]
	\arrow["{\ap_R{\ptd{g}}, \ap_D{\ptd{f}}}", squiggly, from=8-1, to=9-1]
	\arrow["{\bsd{R}, \bsd{D}}", squiggly, from=9-1, to=10-1]
	\arrow["\text{algebra}", squiggly, from=10-1, to=11-1]
	\arrow["\text{algebra}", squiggly, from=4-1, to=5-1]
	\arrow["\text{algebra}", squiggly, from=5-1, to=6-1]
\end{tikzcd}}\]
where $e_X, e_Y, e_Z, e_W$ are equivalences which are the sides of the cube and $L, R, U, D$ are the left, right, upper and lower faces of it.
\\
If we don't require that $e_\bullet$ be equivalences, this can be taken as a definition of a map of pointed squares.
\end{lem}

We only use a simpler version of the lemma:
\begin{re}\label{re:equality-of-squares-simple-version}
    In the special case when $e_X, e_Y, e_Z, e_W$ are the identity maps and $L$ and $R$ are reflexivity, item (2) above reduces to showing that the square

    \[\begin{tikzcd}
    	{{r \circ g}(\bs{X})} && {{r \circ g'}(\bs{X})} \\
    	& {\bs{W'}} \\
    	{{l \circ f}(\bs{X})} && {{l' \circ f}(\bs{X})}
    	\arrow[""{name=0, anchor=center, inner sep=0}, squiggly, from=1-1, to=1-3]
    	\arrow[squiggly, from=1-1, to=2-2]
    	\arrow[""{name=1, anchor=center, inner sep=0}, squiggly, from=1-1, to=3-1]
    	\arrow[squiggly, from=1-3, to=2-2]
    	\arrow[""{name=2, anchor=center, inner sep=0}, squiggly, from=1-3, to=3-3]
    	\arrow[squiggly, from=3-1, to=2-2]
    	\arrow[""{name=3, anchor=center, inner sep=0}, squiggly, from=3-1, to=3-3]
    	\arrow[squiggly, from=3-3, to=2-2]
    	\arrow["{\bsd{\beta}}", between={0}{0.8}, Rightarrow, from=1, to=2-2]
    	\arrow["{\bsd{r \whisk U}}", between={0}{0.8}, Rightarrow, from=0, to=2-2]
    	\arrow["{\bsd{\beta'}}"', between={0}{0.8}, Rightarrow, from=2, to=2-2]
    	\arrow["{\bsd{D \whisk f}}"', between={0}{0.8}, Rightarrow, from=3, to=2-2]
    \end{tikzcd}\]
equals the evaluation of the 3-filler on $\bs{X}$.
\end{re}

\end{document}